\pdfoutput=1

\documentclass[graybox]{svmult}

\usepackage{layout}

\usepackage{mathptmx}       %
\usepackage{helvet}         %
\usepackage{courier}        %

\usepackage{graphicx}        %
\usepackage{multicol}        %
\usepackage[bottom]{footmisc}%
\usepackage{mathrsfs}

\usepackage[numbers,sort,compress]{natbib}
\usepackage[fleqn]{amsmath}
\usepackage{ascmac}
\usepackage{psfrag}
\usepackage{amsmath}
\usepackage{amsbsy}
\usepackage{amssymb}
\usepackage{latexsym}
\usepackage{color}

\newtheorem{scenario}{\bf Scenario}
\newtheorem{fact}[theorem]{Fact}

\newtheorem{bold-example}[theorem]{\bf Example}
\newtheorem{bold-remark}[theorem]{\bf Remark}

\newtheorem{prop}[theorem]{\bf Proposition}
\newtheorem{cor}[theorem]{\bf Corollary}
\newtheorem{lem}[theorem]{\bf Lemma}

\newcommand{\bm}[1]{\mbox{\boldmath$#1$}}

\DeclareSymbolFont{largesymbolsA}{U}{txexa}{m}{n}
\DeclareMathSymbol{\varprod}{\mathop}{largesymbolsA}{16}

\newcommand{\reg}{\mathfrak{C}}
\newcommand{\olabel}[1]{\nonumber}

\begin{document}

\title{Hierarchical Convex Optimization \\ by the Hybrid Steepest Descent Method \\ with Proximal Splitting Operators \\ --- Enhancements of SVM and Lasso}
\titlerunning{Hierarchical Convex Optimization by the Hybrid Steepest Descent Method}
\author{Isao Yamada and Masao Yamagishi}
\authorrunning{}
\institute{Isao Yamada
  \at{Department of Information and Communications Engineering, Tokyo Institute of Technology, S3-60, Ookayama, Meguro-ku, Tokyo 152-8552, Japan,}
  \email{isao@ict.e.titech.ac.jp}
\and 
Masao Yamagishi
  \at{Department of Information and Communications Engineering, Tokyo Institute of Technology, S3-60, Ookayama, Meguro-ku, Tokyo 152-8552, Japan,}
  \email{myamagi@ict.e.titech.ac.jp \ \\
    \ \\
    {\footnotesize This manuscript is the \emph{accepted version} of ``Chapter 16 (by Isao Yamada and Masao Yamagishi, pp.~413--489) in: H.~H.~Bauschke, R.~S.~Burachik, D.~R.~Luke (eds.) Splitting Algorithms, Modern Operator Theory, and Applications, Springer, 2019.''
This version is not completely same as the published version in the book.
      For precise information, the authors would like to encourage to consult the published version in the book.}
  }
}
\maketitle

\abstract*{The breakthrough ideas in the modern proximal splitting  methodologies allow us to express the set of all minimizers of a superposition of multiple nonsmooth convex functions as the fixed point set of 
computable nonexpansive operators. In this paper, we present practical algorithmic strategies for the hierarchical convex optimization problems which require further strategic selection of a most desirable vector from the solution set of the standard convex optimization. The proposed algorithms are established by applying 
the hybrid steepest descent method to special nonexpansive operators designed through the art of proximal splitting. 
We also present applications of the proposed strategies to certain unexplored hierarchical enhancements of the support vector machine and the Lasso estimator. }

\abstract{The breakthrough ideas in the modern proximal splitting  methodologies allow us to express the set of all minimizers of a superposition of multiple nonsmooth convex functions as the fixed point set of 
computable nonexpansive operators. In this paper, we present practical algorithmic strategies for the hierarchical convex optimization problems which require further strategic selection of a most desirable vector from the solution set of the standard convex optimization. The proposed algorithms are established by applying 
the hybrid steepest descent method to special nonexpansive operators designed through the art of proximal splitting. 
We also present applications of the proposed strategies to certain unexplored hierarchical enhancements of the support vector machine and the Lasso estimator. }



\newcommand{\rev}[1]{#1}
\newcommand{\revb}[1]{#1}
\newcommand{\revc}[1]{#1}
\renewcommand{\theequation}{\arabic{equation}}
\section{Introduction}
\label{sec:1}
\rev{Convex} optimization has been playing \rev{a central role} in a broad range of mathematical sciences and engineering. Many optimization tasks in such applications can be interpreted as special instances of the following simple model:
\begin{equation}\label{typical-conv-opt}
\mbox{\rm minimize }f(x)+g(Ax)\mbox{ subject to }x\in {\cal X},
\end{equation}
where $({\cal X}, \langle \cdot,\cdot\rangle_{\cal X},\|\cdot\|_{\cal X})$, $({\cal K}, \langle \cdot,\cdot\rangle_{\cal K},\|\cdot\|_{\cal K})$ are real Hilbert spaces, $f:{\cal X}\rightarrow (-\infty,\infty]$
and $g:{\cal K}\rightarrow (-\infty,\infty]$ are proper lower semicontinuous convex functions, i.e., $f\in \Gamma_{0}({\cal X})$ and 
$g\in \Gamma_{0}({\cal K})$, and $A:{\cal X}\rightarrow {\cal K}$ is a bounded linear operator, 
i.e., $A\in {\cal B}({\cal X},{\cal K})$. Such a unified simplification is indebted entirely to the remarkable expressive ability of the abstract  Hilbert space. For example, a seemingly much more general model:  
\begin{equation}\label{multiple-typical-conv-opt}
\mbox{\rm find }x^{\star}\in {\cal S}:=\underset{x\in {\cal X}}{\operatorname{argmin}}
\left[
\Phi(x):=f(x)+\sum_{i=1}^{m}g_{i}(A_{i}x)
\right]\neq \varnothing,
\end{equation}
where $({\cal X}, \langle \cdot,\cdot\rangle_{\cal X},\|\cdot\|_{\cal X})$ and $({\cal K}_{i}, \langle \cdot,\cdot\rangle_{{\cal K}_{i}},\|\cdot\|_{{\cal K}_{i}})$ ($i=1,2,\ldots, m$) are real Hilbert spaces, $f\in \Gamma_{0}({\cal X})$,  
$g_{i}\in \Gamma_{0}({\cal K}_{i})$ ($i=1,2,\ldots, m$), 
and $A_{i}\in {\cal B}({\cal X},{\cal K}_{i})$ ($i=1,2,\ldots, m$), can also be translated into the problem in (\ref{typical-conv-opt}) by redefining a new Hilbert space  
\begin{equation}\label{prod-space0}
{\cal K}:={\cal K}_{1}\times \cdots \times {\cal K}_{m}=
\left\{
{\mathbf x}=(x_{1},\ldots, x_{m})\mid x_{i}\in {\cal K}_{i}\ (i=1,\ldots, m)\right\}
\end{equation}  
equipped with the addition $({\mathbf x},{\mathbf y})\mapsto 
(x_{1}+y_{1},\ldots, x_{m}+y_{m})$, the scalar multiplication 
$(\alpha, {\mathbf x})\mapsto (\alpha x_{1},\ldots, \alpha x_{m})$, and the inner product $({\mathbf x},{\mathbf y})\mapsto 
\langle {\mathbf x},{\mathbf y}\rangle_{\cal K}:=\sum_{i=1}^{m}\langle x_{i},y_{i}\rangle_{{\cal K}_{i}}$, a new convex function 
\begin{equation}\label{prod-space-g}
g:=\bigoplus_{i=1}^{m}g_{i}:  
{\cal K}\rightarrow (-\infty,\infty]: (x_{1},\ldots, x_{m})\mapsto 
\sum_{i=1}^{m}g_i(x_{i}), 
\end{equation}
and a new bounded linear operator 
\begin{equation}\label{prod-space-A}
A:{\cal X}\rightarrow {\cal K}: x\mapsto (A_{1}x,\ldots, A_{m}x). 
\end{equation}
Indeed, for many years, the model (\ref{multiple-typical-conv-opt}) 
has been accepted widely as a standard, where 
all players, $f, g_{i}\circ A_{i}\in \Gamma_{0}({\cal X})$ ($i=1,\ldots,m$) in (\ref{multiple-typical-conv-opt}) are 
designed strategically by users in order to achieve, after optimization, a valuable vector satisfying their requirements. 

The so-called \emph{proximal splitting methodology} has been built, on the rich mathematical foundations of convex analysis, monotone operator theory and fixed point theory of nonexpansive operators (see, e.g., \cite{bauschke11,Vu13,Condat13,combettes2015compositions}), in order to broaden the applicability of the proximity operators of convex functions \cite{Moreau62}, e.g., to the model (\ref{multiple-typical-conv-opt}). 
It is well-known that the solution set ${\cal S}$ in (\ref{multiple-typical-conv-opt}) can be characterized completely  
as the zero of the set-valued operator $\partial \Phi:{\cal X}\rightarrow 2^{\cal X}: x\mapsto 
\{u\in {\cal X}\mid \Phi(y)\geq \Phi(x)+\langle y-x,u\rangle_{\cal X}\ (\forall y\in {\cal X})\}$, called the \emph{subdifferential} of $\Phi$. The maximal monotonicity of $\partial \Phi$ provides us with 
further equivalent fixed point characterization in terms of a single valued operator 
${\rm prox}_{\Phi}:=\left({\rm I}+\partial \Phi\right)^{-1}:{\cal X}\rightarrow {\cal X}: 
x\mapsto \operatorname{argmin}_{y\in {\cal X}}\Phi(y)+\frac{1}{2}\|x-y\|_{\cal X}^{2}$ called the proximity operator of $\Phi$ (see Section~\ref{sec:Monotone}) \revb{[Note: The identity operator is denoted by ${\rm I}\colon \mathcal{X} \to \mathcal{X}$ but the common notation ${\rm I}$ is going to be used for the identity operator on any real Hilbert space in this paper]}. This fact is simply stated as 
\begin{equation}
\left(\forall z\in {\cal X}\right)\quad z\in {\cal S}\Leftrightarrow 0\in 
\partial \Phi(z)\Leftrightarrow z={\rm prox}_{\Phi}(z)
\olabel{1st-characterization}
\end{equation}
and some algorithms can generate, from any $x_{0}\in {\cal X}$, a weakly convergent sequence $x_{n}\in {\cal X}$ ($n\in {\mathbb N}$)  to 
a point in ${\cal S}$ in (\ref{multiple-typical-conv-opt}) if ${\rm prox}_{\Phi}$ is available 
as a computational tool (see, e.g., \cite{Martinet70,Martinet72,Rock76}). 
A simplest example of such algorithms generates a sequence 
$(x_{n})_{n\in {\mathbb N}}$ by    
\begin{equation}\label{prox-point-alg1}
x_{n+1}={\rm prox}_{\Phi}(x_{n})\quad (n=0,1,2,\ldots). 
\end{equation}
The algorithm (\ref{prox-point-alg1}) can be interpreted as a straightforward application of \emph{Krasnosel'ski\u{\i}-Mann Iterative Process} (see Fact \ref{fa:Mann} in Section \ref{sec:Monotone}) because ${\rm prox}_{\Phi}$ is known to be firmly nonexpansive, i.e., $2{\rm prox}_{\Phi}-{\rm I}:{\cal X}\rightarrow {\cal X}$ is a nonexpansive operator (see (\ref{def-nexp})).  
Although the above strategy in \eqref{prox-point-alg1} is  conceptually simple and elegant, its applicability has been very limited because 
the computation of ${\rm prox}_{\Phi}(x)$ requires to solve a regularized convex optimization problem 
$\min \Phi(\cdot)+\frac{1}{2}\|x-\cdot\|_{\cal X}^{2}$ whose unique solution is still hard to be computed for most scenarios of type 
(\ref{multiple-typical-conv-opt}) in many application areas.

On the other hand, there are many scenarios \rev{that fall} in the model (\ref{multiple-typical-conv-opt}) 
where the proximity operators of the all players, i.e., ${\rm prox}_f:{\cal X}\rightarrow {\cal X}$ and 
${\rm prox}_{g_{i}}:{\cal K}_{i}\rightarrow {\cal K}_{i}$ ($i=1,\ldots,m$), are available as computational tools while 
${\rm prox}_{\Phi}$ is not practically available (see, e.g., \cite{plc-pesq11,combettes2015compositions}). 
A major goal of recent active studies (see, e.g., \cite{plc-vrw05,plc-jcp07,Vu13,Condat13,YY2017}) on the \emph{proximal splitting methodology} has been the creation of  
more applicable iterative algorithms, for (\ref{multiple-typical-conv-opt}) and its variations, than (\ref{prox-point-alg1}) by utilizing 
computable tools ${\rm prox}_f$ and ${\rm prox}_{g_{i}}$ ($i=1,\ldots,m$) simultaneously.  
Such effort has culminated in many powerful algorithms which have been applied successfully to the broader classes of optimizations including the standard model (\ref{multiple-typical-conv-opt}).   

Usually, the standard model (\ref{multiple-typical-conv-opt}) is formulated in the form of a weighted average of multiple convex functions and 
the weights are designed in accordance with the level of importance of each convex function. However 
quantification of the level of importance is often challenging as well as influential to the final results of optimizations (see \revb{Section~\ref{subsec:5.1}} for a recent advanced strategy of such a parameter tuning for the  Lasso estimator which is a standard sparsity aware statistical estimation method).   
By keeping in mind (i) the remarkable flexibility of the standard model (\ref{multiple-typical-conv-opt}) 
proven extensively in many successful applications of the modern proximal splitting methodology, as well as (ii) 
the inherent difficulty in the weight design of multiple convex functions in (\ref{multiple-typical-conv-opt}), 
a question arises: \emph{Is there any alternative model of (\ref{multiple-typical-conv-opt}) which can also serve as a natural optimization strategy for multiple convex criteria} ? To see the light of the tunnel regarding this primitive question,
let us start to imagine important elements for us to consider in \emph{finding residence}. We may consider the house rent, 
the residential environment including living space and housing equipment, the
neighborhood environment, 
the accessibility to public 
transportation systems, and the commuting time, etc. We would prioritize the elements, e.g., firstly by narrowing down the candidates 
to the set $S_{1}$ of all residents of which the rents and commuting times are in your acceptable range. Next we may try to narrow 
down the candidates to the set $S_{2}(\subset S_{1})$ of all residents whose living spaces achieve maximum level among all in $S_{1}$. Further, 
we may probably like to select residents in $S_{2}$ as final choices by choosing the best ones, e.g., in the sense of the neighborhood environment or 
the housing equipment. This simple example suggests that we often optimize multiple criteria 
one by one hierarchically rather than  optimize the sum of different criteria at once 
certainly because there exists no universal justification for adding different criteria. 
In fact, many mathematicians and scientists have been challenging to pave the way for  
the so-called hierarchical convex optimization problems (see, e.g., \cite{Tikhonov63,Polyak79,Mang81,Dontchev-Zolezzi93,Com94,Attouch96,yama0,yama1,C-B99,Com2000,yama2001a,BC2001,Cominetti-Courdurier02,y-o-s2002,Cabot05,yamada2011minimizing,OYTSP}). Landmark theories toward $\mathcal{M}$-stage hierarchical convex optimization are found, e.g.,
in \cite{Attouch96,Cabot05} where, for given 
 $\Phi_{i}\in \Gamma_{0}({\mathcal X})$ ($i=0,1,\ldots,\mathcal{M}$) satisfying  
 $S_{i}:={\displaystyle \underset{x\in S_{i-1}}{\operatorname{argmin}}\Phi_{i}(x)\neq \varnothing}$ ($i=0,1,\ldots,\mathcal{M}$) with $S_{-1}:={\mathcal X}$, their major goals
 are set to establish computational strategies for iterative approximation of a point in $S_{\mathcal{M}}$. 
(Note: Every point in  $S_{i}$ of the hierarchical convex optimization is called a \emph{viscosity solution} of $S_{i-1}$ ($i=1,2,\ldots, \mathcal{M}$). To avoid confusion with "bilevel optimization" in the sense of 
\cite{Candler-Norton77,Vicent-Calamai94,Coloson-Marcotte-Savard07},
we do not use the designation \emph{bilevel optimization} for $S_{1}$ in our hierarchical convex optimization). 
Under the assumptions that $\dim({\cal X})<\infty$ and that $\Phi_{i}\in \Gamma_{0}({\mathcal X})$ ($i=1,2,\ldots,\mathcal{M}$) 
are real valued, Cabot \cite{Cabot05} showed that the sequence $(x_{n})_{n\in {\mathbb N}}$ defined by
\begin{eqnarray}
x_{n+1}&:=&{\rm prox}_{\left(\Phi_{0}+\varepsilon_{n}^{(1)}\Phi_{1}+\varepsilon_{n}^{(2)}\Phi_{2}+
\cdots +\varepsilon_{n}^{(\mathcal{M})}\Phi_{\mathcal{M}}\right)}(x_{n})\nonumber \\
&=&\left({\rm I}+\partial \left(\Phi_{0}+\varepsilon_{n}^{(1)}\Phi_{1}+\varepsilon_{n}^{(2)}\Phi_{2}+
\cdots +\varepsilon_{n}^{(\mathcal{M})}\Phi_{\mathcal{M}}\right)\right)^{-1}(x_{n})\label{Cabot05-alg}
\end{eqnarray}
satisfies (i) $\lim_{n\rightarrow \infty}d(x_{n}, S_{\mathcal{M}})=0$ and (ii) 
($\forall i\in \{0,1,\ldots, \mathcal{M}\}$) $\lim_{n\rightarrow \infty}\Phi_{i}(x_{n})=\min_{x\in S_{i-1}}\Phi_{i}(x)$ 
if positive number sequences $(\varepsilon_{n}^{(0)}:=1)_{n\in {\mathbb N}}$
and $(\varepsilon_{n}^{(i)})_{n\in {\mathbb N}}$ ($i\in \{1,\ldots,\mathcal{M}\}$) satisfy  
certain technical conditions including ${\displaystyle \lim_{n\rightarrow \infty}\varepsilon_{n}^{(i)}=0}$,  
${\displaystyle \lim_{n\rightarrow \infty}\frac{\varepsilon^{(i)}_{n}}{\varepsilon^{(i-1)}_{n}}=0}$ ($i=1,2,\ldots, \mathcal{M}$) and 
$\sum_{n=0}^{\infty}\varepsilon_{n}^{(\mathcal{M})}=\infty$ {\rm [Note: The scheme (\ref{Cabot05-alg}) is a simplified version of the  original scheme in \cite{Cabot05} by restricting to the case $\lambda_{n}=1$  and 
$\eta_{n}=0$ ($n\in {\mathbb N}$)]}. 

Clearly, the algorithms (\ref{Cabot05-alg}) and 
(\ref{prox-point-alg1}) have essentially a common limitation in their practical applicabilities 
because (\ref{Cabot05-alg}) requires ${\rm prox}_{\left(\Phi_{0}+\varepsilon_{n}^{(1)}\Phi_{1}+\varepsilon_{n}^{(2)}\Phi_{2}+
\cdots +\varepsilon_{n}^{(\mathcal{M})}\Phi_{\mathcal{M}}\right)}$, or its very good approximation, for every update in generation of $(x_{n})_{n\in {\mathbb N}}$. By recalling the breakthrough ideas developed in the recent \emph{proximal splitting methodology} for resolution of the inherent limitation in (\ref{prox-point-alg1}), an ideal as well as possibly realistic assumption to be imposed upon each player $\Phi_{i}:{\mathcal X}\rightarrow (-\infty,\infty]$ ($i=0,1,\ldots,\mathcal{M}$) in the above $\mathcal{M}$-stage hierarchical convex optimization seems to be certain differentiability assumptions or proximal decomposability assumptions, e.g.  
\begin{equation}\olabel{assump-conmbination-multiple-conv-opt}
\Phi_{i}(x):=f_{i}(x)+\sum_{\iota_{(i)}=1}^{M_{i}}g_{\iota_{(i)}}(A_{\iota_{(i)}}x),
\end{equation} 
with real Hilbert spaces $({\cal K}_{\iota_{(i)}}, \langle \cdot,\cdot\rangle_{{\cal K}_{\iota_{(i)}}},\|\cdot\|_{{\cal K}_{\iota_{(i)}}})$ ($\iota_{(i)}=1,2,\ldots, M_{i}$), $f_{i}\in \Gamma_{0}({\cal X})$,  
$g_{\iota_{(i)}}\in \Gamma_{0}({\cal K}_{i})$ ($i=0,1,\ldots, \mathcal{M}$), 
and bounded linear operators $A_{\iota_{(i)}}:{\cal X}\rightarrow {\cal K}_{\iota_{(i)}}$ ($\iota_{(i)}=1,2,\ldots, M_{i}$), where ${\rm prox}_{f_{i}}:{\cal X}\rightarrow {\cal X}$ and 
${\rm prox}_{g_{\iota_{(i)}}}:{\cal K}_{\iota_{(i)}}\rightarrow {\cal K}_{\iota_{(i)}}$ ($\iota_{(i)}=1,\ldots,M_{i}$), are available as computational tools while ${\rm prox}_{\Phi_{i}}$ is not necessarily available. 

In this paper, we choose to cast our primary target in the   
iterative approximation of a solution of 
\begin{equation}
\label{viscosity-conv-prob-1}
\mbox{\rm minimize }\Psi(x^{\star})\mbox{ \rm subject to }x^{\star}\in \underset{x\in {\cal X}}{\operatorname{argmin}}
\left[
\Phi(x):=f(x)+\sum_{i=1}^{m}g_{i}(A_{i}x)
\right]\neq \varnothing,
\end{equation}
i.e., a viscosity solution of the convex optimization problem 
(\ref{multiple-typical-conv-opt}), where we assume that $\Psi\in \Gamma_{0}({\cal X})$  is  
{\rm G\^{a}teaux} differentiable with Lipschitzian gradient $\nabla \Psi:{\cal X}\rightarrow {\cal X}$, i.e., 
\begin{equation}\olabel{grad-lip}
(\exists \kappa>0,\  \forall x,y\in {\cal X})\quad \|\nabla \Psi(x)-\nabla \Psi(y)\|\leq \kappa \|x-y\|, 
\end{equation}
and that  ${\rm prox}_f:{\cal X}\rightarrow {\cal X}$ and 
${\rm prox}_{g_{i}}:{\cal K}_{i}\rightarrow {\cal K}_{i}$ ($i=1,\ldots,m$) are available as computational tools. 

Although the application of such iterative algorithms is certainly restrictive compared to the overwhelming potential of the general hierarchical convex optimization, our target is realistic and still allows us to cover many applications of interest to practitioners who are searching for a step ahead 
optimization strategy and yet to be able to exploit maximally the central ideas in the modern proximal splitting methodologies. Especially for practitioners, we remark that if the suppression of $\sum_{k=1}^{L}\psi_{k}\circ B_{k}\in \Gamma_{0}({\cal X})$ 
over $\operatorname{argmin}_{x\in {\cal X}}\Phi(x)$ is required, where, for each $k\in \{1,2,\ldots, L\}$, ${\cal Y}_{k}$ is a real Hilbert space, $\psi_{k}\in \Gamma_{0}({\cal Y}_{k})$, ${\rm prox}_{\gamma \psi_{k}}:{\cal Y}_{k}\rightarrow {\cal Y}_{k}$ ($\gamma>0$) is available as computational tools, and $B_{k}\in {\cal B}({\cal X},{\cal Y}_{k})$, 
such a mission could be achieved satisfactorily by considering an alternative problem below of type (\ref{viscosity-conv-prob-1}):     
\begin{eqnarray}
&&\mbox{\rm minimize }\Psi(x^{\star}):=\sum_{k=1}^{L}{}^{\gamma}\psi_{k}(B_{k}x^{\star})\nonumber \\
&&\mbox{ \rm subject to }x^{\star}\in \underset{x\in {\cal X}}{\operatorname{argmin}}
\left[
\Phi(x):=f(x)+\sum_{i=1}^{m}g_{i}(A_{i}x)
\right]\neq \varnothing,  \label{viscosity-conv-prob-2}
\end{eqnarray}
where (i) ${}^{\gamma}\psi_{k}:{\cal Y}_{k}\rightarrow {\mathbb R}:
y_{k}\mapsto \min_{y\in {\cal Y}_{k}} \psi_{k}(y)+\frac{1}{2\gamma}\|y-y_{k}\|^{2}_{{\cal Y}_{k}}
$ ($k=1,2,\ldots, L$) are \emph{the Moreau envelopes (or the Moreau-\revb{Yosida} regularizations)} of a sufficiently small index $\gamma>0$ (see Fact \ref{fa:MoreauEnvelope} in Section~\ref{sec:Monotone} and \cite{yamada2011minimizing}).
This is because $\lim_{\gamma\downarrow 0}{}^{\gamma}\psi_{k}(y_{k})=\psi_{k}(y_{k})$ ($\forall y_{k}\in 
{\rm dom} \psi_{k}:=\{y\in {\cal Y}_{k}\mid \psi_{k}(y)<\infty\}$) and 
${}^{\gamma}\psi_{k}$ is {\rm G\^{a}teaux} differentiable with $\frac{1}{\gamma}$-Lipschitzian 
$\nabla {}^{\gamma}\psi_{k}:{\cal Y}_{k}\rightarrow {\cal Y}_{k}:y_{k}\mapsto \frac
{y_{k}-{\rm prox_{\gamma \psi_{k}}(y_{k})}}{\gamma}$ and therefore 
\begin{eqnarray*}
&&(\forall x_{1},x_{2}\in {\cal X})\quad 
\left\|\nabla \sum_{k=1}^{L}\left({}^{\gamma}\psi_{k}\circ B_{k}\right)(x_{1})-
\nabla \sum_{k=1}^{L}\left({}^{\gamma}\psi_{k}\circ B_{k}\right)(x_{2})\right\|_{\cal X}
\\
&=&
\left\|\sum_{k=1}^{L}
B_{k}^{*}\nabla {}^{\gamma}\psi_{k}\left(B_{k}x_{1}\right)
-
\sum_{k=1}^{L}
B_{k}^{*}\nabla {}^{\gamma}\psi_{k}\left(B_{k}x_{2}\right)
\right\|_{\cal X}
\leq \sum_{k=1}^{L}\frac{\|B_{k}\|_{\rm op}^{2}}{\gamma}\|x_{1}-x_{2}\|_{\cal X},
\end{eqnarray*}
where $B_{k}^{*}\in {\cal B}({\cal Y}_{k},{\cal X})$ is 
the conjugate of $B_{k}\in {\cal B}_{k}({\cal X},{\cal Y}_{k})$ and $\|\cdot\|_{\rm op}$ stands for the operator norm.

Fortunately, by introducing the exactly same translation used in the reformulation of 
Problem (\ref{multiple-typical-conv-opt}) as an instance of Problem (\ref{typical-conv-opt}), 
our problem (\ref{viscosity-conv-prob-1}) can also be simplified as 
\begin{equation}
\label{viscosity-conv-prob}
\mbox{\rm minimize }\Psi(x^{\star})\mbox{ \rm subject to }x^{\star}\in {\cal S}_{p}:=\underset{x\in {\cal X}}{\operatorname{argmin}}
\left[
f(x)+g(Ax)
\right]\neq \varnothing, 
\end{equation}
where ${\cal K}$, $g:{\cal K}\rightarrow (-\infty, \infty]$, and $A:{\cal X}\rightarrow {\cal K}$ are defined respectively\footnote{
  There are many practical conditions for $(f,g,A)$ to guarantee ${\cal S}_p \not = \varnothing$, see, e.g., \cite{Zeid3,bauschke11} and Fact~\ref{fa:uobqegbuqeh} in Section \ref{subsec:2.1}.
}
 by (\ref{prod-space0}), (\ref{prod-space-g}), and (\ref{prod-space-A}), and we can assume 
that (i) $\Psi\in \Gamma_{0}({\cal X})$ is {\rm G\^{a}teaux} differentiable with Lipschitzian gradient $\nabla \Psi:{\cal X}\rightarrow {\cal X}$, and 
that (ii) ${\rm prox}_f:{\cal X}\rightarrow {\cal X}$ and 
${\rm prox}_{g}:{\cal K}\rightarrow {\cal K}$ are available as computational tools because 
\begin{eqnarray}
&&{\rm prox}_{g}({\mathbf x}):=\underset{{\mathbf y}\in {\cal K}}{\operatorname{argmin}}
\left[
g({\mathbf y})+\frac{1}{2}\|{\mathbf y}-{\mathbf x}\|_{\cal K}^{2}
\right] \nonumber \\
&=&\underset{(y_{1},\ldots, y_{m})\in {\cal K}_{1}\times \cdots \times  {\cal K}_{m}}{\operatorname{argmin}} \ 
\sum_{i=1}^{m}\left[g_{i}(y_{i})+\frac{1}{2}\|y_{i}-x_{i}\|_{{\cal K}_{i}}^{2}\right] \nonumber\\
&=&\left({\rm prox}_{g_{1}}(x_{1}), \ldots, {\rm prox}_{g_{m}}(x_{m})\right).
    \label{eq:intro_proxdecompose}
\end{eqnarray}

The following two scenarios suggest the remarkable advantage achieved by algorithmic solutions to (\ref{viscosity-conv-prob-1}).   

\begin{scenario}\label{sce1} (Unification of Conditional Optimization Models)  
\emph{Let $f_{\langle {\cal D}\rangle}\in \Gamma_{0}({\cal X})$ and $g_{i\langle {\cal D}\rangle}\in \Gamma_{0}({\cal K}_{i})$ ($i=1,2,\ldots, m$) be  nonnegative valued functions which are defined with observed data ${\cal D}$.
Suppose that there exists a well-established data analytic strategy which utilizes with $\Psi\in \Gamma_{0}({\cal X})$ as  
\begin{eqnarray}
&&\mbox{\rm find }x^{\star}\in \underset{x\in {\cal S}_{0}}{\operatorname{argmin}} \ \Psi(x),\nonumber \\
&&\mbox{ where }{\cal S}_{0}:=
\left\{
x\in {\cal X}\mid f_{\langle {\cal D}\rangle}(x)=g_{i\langle {\cal D}\rangle}(A_ix)=0\ (i=1,2,\ldots,m)
\right\}, \label{optimization-in-ideal-case}
\end{eqnarray}
provided that the data ${\cal D}$ is consistent, i.e., it satisfies ${\cal S}_{0}\neq \varnothing$. }

\emph{
However, to deal with more general data ${\cal D}$, it is important 
to establish a mathematically sound extension of the above data analytic strategy to be applicable even to inconsistent data ${\cal D}$ s.t.  
${\cal S}_{0}=\varnothing$. One of the most natural extensions of (\ref{optimization-in-ideal-case}) would be 
the following hierarchical formulation:  
\begin{eqnarray}
&&\mbox{\rm find }x^{\star\star}\in \underset{x^{\star}\in {\cal S}_{\langle {\cal D}\rangle}}{\operatorname{argmin}} \ \Psi(x^{\star}),\nonumber \\
&&\mbox{ where }{\cal S}_{\langle {\cal D}\rangle}:=
\underset{x\in {\cal X}}{\operatorname{argmin}} 
\left[
f_{\langle {\cal D}\rangle}(x)+\sum_{i=1}^{m}g_{i \langle {\cal D}\rangle}(A_{i}x)
   \right],
   \olabel{optimization-in-general-case}
\end{eqnarray}
because  ${\cal S}_{\langle {\cal D}\rangle}\neq \varnothing$ holds under weaker assumption than ${\cal S}_{0}\neq \varnothing$, and 
${\cal S}_{\langle {\cal D}\rangle}={\cal S}_{0}$ holds true if ${\cal S}_{0}\neq \varnothing$.
However, the well-established data analytic strategies only for consistent data 
${\cal D}$ in the form of (\ref{optimization-in-ideal-case}) have often been modified, with the so-called \emph{tuning parameter} $\reg>0$, to 
\begin{equation}
\mbox{\rm find }\tilde{x}^{\star}\in 
\underset{x\in {\cal X}}{\operatorname{argmin}}
\left[
\frac{1}{\reg} \Psi(x)+f_{\langle {\cal D}\rangle}(x)+\sum_{i=1}^{m}g_{i{\langle {\cal D}\rangle}}(A_{i}x)
\right], \label{optimization-in-general-case-fake}
\end{equation}
which is not really an extension of (\ref{optimization-in-ideal-case}) because the model (\ref{optimization-in-general-case-fake})
unfortunately has no guarantee to produce $x^{\star}$ in (\ref{optimization-in-ideal-case}) even if ${\cal D}$ satisfies ${\cal S}_{0}\neq \varnothing$.
}  
\end{scenario} 

\begin{scenario}\label{sce2}  
\emph{Suppose that we are interested in an estimation problem of a desired vector in ${\cal X}$ 
at which the functions $f, g_{i}\circ A_{i}\in \Gamma_{0}({\cal X})$ ($i=1,\ldots,m$) in (\ref{multiple-typical-conv-opt})
 are known to achieve small values and therefore 
the model (\ref{multiple-typical-conv-opt}) has been employed as 
an estimation strategy. Suppose also that we newly found another effective criterion 
$\Psi\in \Gamma_{0}({\cal X})$ which likely to achieve small 
values
around the desired vector to be estimated.   
In such a case, our common utilization of $\Psi$, for improvement of the previous strategy, has often been modeled as a new optimization problem: 
\begin{equation}\label{multiple-typical-conv-opt2}
\mbox{\rm find }\tilde{x}^{\star}\in \widetilde{\cal S}:=\underset{x\in {\cal X}}{\operatorname{argmin}}
\left[
f(x)+\sum_{i=1}^{m}g_{i}(A_{i}x)+\Psi(x)
\right]\neq \varnothing.
\end{equation}
However, it is essentially hard to tell which is better between the estimation strategies (\ref{multiple-typical-conv-opt}) and 
(\ref{multiple-typical-conv-opt2}) because the criteria in these optimizations are different. Indeed, 
$\tilde{x}^{\star}$ does not necessarily achieve best in the sense of the model (\ref{multiple-typical-conv-opt}) while $x^{\star}$ certainly achieves best in the sense of the model (\ref{multiple-typical-conv-opt}). 
}
\emph{
On the other hand, if we formulate a new optimization problem, from a hierarchical optimization point of view, 
e.g., as   
\begin{equation}\label{multiple-typical-conv-opt3}
\mbox{\rm find }x^{\star\star}\in \underset{x^{\star}\in {\cal S}}{\operatorname{argmin}} \ \Psi(x^{\star}), \mbox{ where }{\cal S}:=\underset{x\in {\cal X}}{\operatorname{argmin}}
\left[
f(x)+\sum_{i=1}^{m}g_{i}(A_{i}x)
\right]\neq \varnothing, 
\end{equation}
its solution $x^{\star\star}$ certainly meets more faithfully all the requirements than $x^{\star}\in {\cal S}$ because  
both $x^{\star\star}, x^{\star}\in {\cal S}$ and $\Psi(x^{\star\star})\leq \Psi(x^{\star})$ are achieved. 
}
\end{scenario}

The following examples suggest that the hierarchical optimization   
has been offering well-grounded direction for advancement of computational strategies in inverse problems and data sciences.  

\begin{bold-example} (Hierarchical Convex Optimizations in Real World Applications\footnote{To the best of the authors' knowledge, little has been reported on the hierarchical \emph{nonconvex} optimization. We remark that the 
MV-PURE (Minimum-variance pseudo-unbiased reduced-rank estimator) (see, e.g., \cite{Yamada-E06,P-Yamada08,P-C-Yamada09}), for the unknown vector possibly subjected to linear constraints, 
is defined by a closed form solution of a certain hierarchical nonconvex optimization problem which characterizes 
a natural reduced rank extension of the Gauss-Markov (BLUE) estimator \cite{Luenberger69,Kailath95} to the case of reduced-rank estimator.
It was shown in \cite{P-Yamada08} that specializations of the MV-PURE include Marquardt's reduced rank estimator \cite{Marquardt70}, Chipman-Rao estimator \cite{ChipmanRao64}, and Chipman's reduced rank estimator \cite{Chipman99}.
In Section \ref{subsec:5.2} of this paper, we newly present a special instance of a hierarchical \emph{nonconvex} optimization problem which can be solved through multiple  hierarchical \emph{convex} optimization subproblems.
})
\label{ex:into}
\begin{enumerate}
\item[\rm (a)] (Generalized inverse / Moore-Penrose inverse \cite{Moore20,Penrose55,Rao-Mitra71,Ben-Israel,bauschke11})
Let $\cal X$ and ${\cal K}$ be real Hilbert spaces, let  $A\in {\cal B}({\cal X},{\cal K})$ be 
such that $\operatorname{ran}(A):=\{A(x)\in {\cal K}\mid x\in {\cal X}\}$ is closed. 
Then for every $y\in {\cal K}$, $C_{y}:=\{x\in {\cal X}\mid \|Ax-y\|_{\cal K}=\min_{z\in {\cal X}}\|Az-y\|_{\cal K}\}=
\{x\in {\cal X}\mid A^{*}A(x)=A^{*}(y)\}\neq \varnothing$.  
The generalized inverse  (in the sense of Moore-Penrose) $A^{\dagger}\in {\cal B}({\cal K},{\cal X})$ is 
defined as $A^{\dagger}:{\cal K}\rightarrow {\cal X}:y\mapsto P_{C_{y}}(0)$,where $P_{C_{y}}$ is 
the orthogonal projection onto $C_{y}$. $A^{\dagger}(y)$ can be seen as the unique solution to the hierarchical convex optimization problem (\ref{multiple-typical-conv-opt3}) for 
$f(z):=\|A(z)-y\|_{\cal K}$, $g_{i}(z):=0$ ($i=1,2,\ldots, m$) and $\Psi(z)=\frac{1}{2}\|z\|_{\cal X}^{2}$. 
The \emph{Moore-Penrose inverse} $A^{\dagger}\in {\cal B}({\cal K},{\cal X})$ of $A\in {\cal B}({\cal X},{\cal K})$ 
has been serving as one of the most natural generalizations of the inverse of $A$, typically in Scenario \ref{sce1}, under the situations 
where the existence of $A^{-1}\in {\cal B}({\cal K},{\cal X})$ is not guaranteed. 
In particular, for finite dimensional settings, there are many ways to express $A^{\dagger}$. 
These include the singular value decomposition of $A^{\dagger}$ in terms of the singular value decomposition 
of $A$.
\item[\rm (b)] (Tikhonov approximation \cite{Tikhonov63,Dontchev-Zolezzi93,Attouch96,bauschke11})  
Let $\Psi,f\in \Gamma_{0}({\cal X})$ and ${\rm argmin}(f)\cap {\rm dom}(\Psi)\neq \varnothing$ where $\Psi$ is coercive and strictly convex. Then $\Psi$ admits 
a unique minimizer $x_{0}$ over ${\rm argmin}(f)$. 
This $x_{0}$ can be seen as the solution of the hierarchical convex optimization in (\ref{multiple-typical-conv-opt3})
for $g_{i}(z):=0$ ($i=1,2,\ldots, m$).   
Moreover, if we define $x_{\varepsilon}\in {\cal X}$ as the unique minimizer of the regularized problem
\begin{equation}\label{reg-problem}
\mbox{\rm miminize }f(x)+\varepsilon \Psi(x)\mbox{ {\rm subject to }}x\in {\cal X}  
\end{equation}
for every $\varepsilon>0$, the desired $x_{0}$ can be approximated as  
(i) $x_{\varepsilon}\rightharpoonup x_{0}$ (as $\varepsilon \downarrow 0$) and 
(ii) $\Psi(x_{\varepsilon})\rightarrow \Psi(x_{0})$  (as $\varepsilon \downarrow 0$).
This fact suggests a strategy for approximating $x_{0}$ if we have a practical way of computing $x_{\varepsilon_{n}}$ for positive sequence $(\varepsilon_{n})_{n=1}^{\infty}$ satisfying $\varepsilon_{n}\downarrow 0$ (as $n\rightarrow \infty$). Many computational approaches to the hierarchical convex optimization seem to have been designed along this strategy.\rev{\footnote{\rev{The behavior of  $(x_{\varepsilon})_{\varepsilon\in(0,1)} \subset \mathcal{X}$ can be analyzed in the context of \emph{approximating curve} for monotone inclusion problem. For recent results combined with Yosida regularization, see \cite{plc-sah06}.}}}
We remark that many formulations of type \eqref{optimization-in-general-case-fake} in Scenario 1 can be seen as instances of \eqref{reg-problem} with $\epsilon = \frac{1}{\reg}$. However, in general, the hierarchical optimality can never be guaranteed by the solution of \eqref{reg-problem} for a fixed constant $\epsilon >0$.
\item[\rm (c)] Assuming differentiability, the iteration of (\ref{Cabot05-alg}) for $\mathcal{M}=1$ can also be interpreted as 
an implicit discretization of the continuous dynamical system:
\begin{equation}\label{dynamical-system1}
\dot{x}(t)+\nabla \Phi_{0}\left(x(t)\right)+\varepsilon(t)\nabla \Phi_{1}\left(x(t)\right)=0,\quad t\geq 0, 
\end{equation}
where $\varepsilon:{\mathbb R}_{+}\rightarrow {\mathbb R}$ is a control parameter tending to $0$ when $t\rightarrow \infty$. 
This observation has been motivating explicit discretization of (\ref{dynamical-system1}) for iterative approximation of point in $S_{1}$, e.g. by       
\begin{equation}
 \olabel{explicit-step0}
x_{n+1}\in x_{n}+{\lambda_{n}} \partial (\Phi_{0}+\varepsilon_{n}\Phi_{0})(x_{n}),
\end{equation}
and its variations (see, e.g., \cite{Solodov07,Solodov08,Helou-DePierro11,Helou-Simoes17}), where $\lambda_{n}$ is a nonnegative stepsize. However this class of algorithms cannot exploit recent advanced proximal splitting techniques for dealing with the constrained set $S_{0}$.  
\item[\rm (d)] Under the assumption that (i) $\Psi$ is {\rm G\^{a}teaux}   
differentiable with Lipschitzian gradient $\nabla \Psi:{\cal X}\rightarrow {\cal X}$, 
and (ii) ${\rm prox}_{f}:{\cal X}\rightarrow {\cal X}$ is available as a computable tool, the \emph{inertial forward-backward algorithm with vanishing Tikhonov regularization} was proposed \cite{Attouch17}, along in the frame of 
accelerated forward-backward methods\footnote{See \cite{Attouch17} on the stream of research, to name but a few, \cite{beck08:_fast_gradien_based_algor_for,Chambolle-Dossal15}, originated from Nesterov's seminal paper \cite{Nesterov83}.}, for an iterative approximation of the solution of a hierarchical convex optimization in (\ref{multiple-typical-conv-opt3})
for $g_{i}=0$ ($i=1,2,\ldots, m$).
\item[\rm (e)] 
In general, the convex optimization problems, especially in the convex feasibility problems \cite{BB96,Com93,C-Z97}, have infinitely many solutions that could be considerably different in terms of other criteria. However most iterative algorithms for  
convex optimization can approximate an anonymous solution of the problem. 
For pursuing a better solution in some other aspects, \emph{superiorization} \cite{Penfold10,Censor-Davidi-Herman10,Herman-Garduno-Davidi-Censor12,Nikazad12} introduces proactively designed perturbations into the original algorithms with preserving preferable convergence properties. 
Essentially, by adopting another criterion $\Psi$, these methods aim to lower the value of $\Psi$ with incorporating a perturbation involving the descent direction of $\Psi$. Apparently, as reported in \cite{YY2017}, 
the hierarchical convex optimization can serve as one of the ideal formulations for the superiorization.  
\item[\rm (f)] 
Let $\Psi\in \Gamma_{0}({\cal X})$ be {\rm G\^{a}teaux}   
differentiable and its gradient $\nabla \Psi:{\cal X}\rightarrow {\cal X}$ is 
Lipschitzian. Suppose that $f\in \Gamma_{0}({\cal X})$ is also {\rm G\^{a}teaux}   
differentiable with Lipschitzian gradient $\nabla f:{\cal X}\rightarrow {\cal X}$
and admits ${\rm argmin} (f+\iota_{K})\neq \varnothing$ for a nonempty closed convex set $K\subset {\cal X}$, 
where $\iota_{K}$ is the indicator function , i.e., 
\begin{equation}
  \olabel{def-ind-func1}
\iota_{K}(x):=\left\{
\begin{array}{ll}
0&\mbox{ if }x\in K,\\
\infty&\mbox{ otherwise. }
\end{array}
\right.
\end{equation}
Then 
\begin{equation}\label{contenp-math}
\mbox{\rm minimize }\Psi(x)\mbox{ {\rm  subject to }}x\in {\rm argmin} (f+\iota_{K})
\end{equation}
can be seen as an instance of the hierarchical convex optimization in (\ref{multiple-typical-conv-opt3})
for $g_{1}:=\iota_{K}$ and $g_{i}=0$ ($i=2,3,\ldots, m$).   
By applying the hybrid steepest descent method \cite{yama0,yama1,D-Y98,ccpip99,yama2001a,y-o-s2002} to several expressions of the set ${\rm argmin} (f+\iota_{K})$ as the fixed point set of certain computable nonexpansive operators 
$T:{\cal X}\rightarrow {\cal X}$ (see, e.g., \cite[Proposition 2.5]{y-o-s2002}, \cite[Example 17.6(b)]{yamada2011minimizing}), practical algorithms have been established to produce a sequence $x_{n}\in {\cal X}$ ($n=0,1,2,\ldots$) which is guaranteed to converge to a solution to Problem (\ref{contenp-math}).
These cover a version of {Projected Landweber method} \cite{Eicke,sabpott,potarun93} for $\Psi(x):=\frac{1}{2}\|x\|^{2}_{\cal X}$ 
and $f(x):=\|A(x)-b\|_{\cal K}$, where $A\in {\cal B}({\cal X},{\cal K})$ and the metric projection ${P}_{K}:{\cal X}\rightarrow K$ is assumed available as a computational tool. As will be discussed below, the main idea of the present paper specialized for Problem (\ref{viscosity-conv-prob}) (or equivalently Problem (\ref{viscosity-conv-prob-1})) is along this simple hierarchical optimization strategy \cite{y-o-s2002,yamada2011minimizing,OYTSP,YY2017,KitaharaYamada16} of applying the hybrid steepest descent method (HSDM: see Section \ref{subsec:2.3}) to the precise expressions of 
the solution sets of the convex optimization problems in terms of fixed point sets of computable nonexpansive operators defined on a certain real Hilbert space ${\cal H}$ which is not necessarily same as the original Hilbert space ${\cal X}$.  
\end{enumerate}
\end{bold-example}

Apparently, to tackle Problem (\ref{viscosity-conv-prob}) (or equivalently Problem (\ref{viscosity-conv-prob-1})), we need to exploit full information on ${\cal S}_{p}$ which is an infinite set in general. 
Moreover, even by using the recently developed powerful proximal splitting algorithms, specially designed for (\ref{typical-conv-opt}), we can produce only some vector sequence that converges to just an anonymous point in ${\cal S}_{p}$,  which implies that we need to add further a new twist to the well-known strategies applicable to Problem (\ref{typical-conv-opt}).  

Fortunately, the unified perspective from the viewpoint of convex analysis and monotone operator theory (see, e.g., \cite{bauschke11}) often enables us to enjoy notable characterizations of the solution set ${\cal S}_{p}$ in terms of  the set of all fixed points of a computable nonexpansive operator defined on certain real Hilbert spaces. Indeed, almost all existing proximal splitting algorithms for Problem (\ref{typical-conv-opt})
 more or less rely on the following type of characterizations of ${\cal S}_{p}$:
\begin{eqnarray}
&&{\cal S}_{p}=\underset{x\in {\cal X}}{\operatorname{argmin}} \  f(x)+g(Ax)=\Xi\left({\rm Fix}(T)\right):=\bigcup_{z\in {\rm Fix}(T)}\Xi(z)\subset {\cal X}
, \label{generic-view1}\\
&&{\rm Fix}(T):=\{z\in {\cal H}\mid T(z)=z\}\qquad \mbox{\rm (Fixed point set of $T$)}, \label{generic-view2} 
\end{eqnarray}
where $({\cal H}, \langle \cdot,\cdot\rangle_{\cal H},\|\cdot\|_{\cal H})$ is a certain real Hilbert space (not necessarily ${\cal H}={\cal X}$), $T:{\cal H}\rightarrow {\cal H}$ is a computable nonexpansive operator, i.e., an operator satisfying   
\begin{equation}\label{def-nexp}
(\forall z_{1},z_{2}\in {\cal H})\quad \|T(z_{1})-T(z_{2})\|_{\cal H}\leq \|z_{1}-z_{2}\|_{\cal H},  
\end{equation}
and $\Xi:{\cal H}\rightarrow 2^{\cal X}$ is a certain set valued operator.
Examples of such characterizations are found
in \cite{E-Y15,YY2017} for the augmented Lagrangian method \cite{hestenes1969multiplier,powell69},
in \cite{plc-vrw05,combettes2015compositions} for the forward-backward splitting approach \cite{Passty79,Gabay83,Tseng91},
in \cite[Proposition 18(iii)]{plc-jcp07} for the Douglas-Rachford splitting approach (see Section \ref{subsec:2.2}) \cite{L-M},
in \cite{E-B92,YY2017} for the alternating direction method of multipliers (ADMM) \cite{L-M,Gabay83,HeYuan12}, 
in \cite{Condat13,Vu13} for the primal-dual splitting method, and
in \cite{YY2017} for a generalized version (see Section \ref{subsec:2.2}) of the linearized augmented Lagrangian method \cite{yang2013linearized}.

If we find a computable nonexpansive operator $T$ satisfying (\ref{generic-view1}) as well as  
a computationally tractable way to extract a point in $\Xi(z)(\subset {\cal X})$ for a given $z\in {\rm Fix}(T)$, 
we can realize an algorithmic solution to Problem (\ref{typical-conv-opt}) by applying 
the so-called \emph{Krasnosel'ski\u{\i}-Mann Iterative Process} (see Fact \ref{fa:Mann} in Section \ref{sec:Monotone}) to $T$, and can produce a weak convergent sequence to a fixed point $z\in {\rm Fix}(T)$, followed by a point extraction from $\Xi(z)$. 
Indeed, the powerful proximal splitting methodologies for Problem (\ref{multiple-typical-conv-opt}) 
seem to have been built more or less along this strategy through innovative designs of computable nonexpansive operators by using ${\rm prox}_f:{\cal X}\rightarrow {\cal X}$ and 
${\rm prox}_{g_{i}}:{\cal K}_{i}\rightarrow {\cal K}_{i}$ ($i=1,\ldots,m$) as computational tools. 

On the other hand, every nonexpansive operator $T:{\cal H}\rightarrow {\cal H}$ can also be plugged into the hybrid steepest descent method for minimizing $\Theta\in \Gamma_{0}(\cal H)$, whose gradient $\nabla \Theta\colon {\cal H} \to {\cal H}$ is Lipschitz continuous, over the fixed point set ${\rm Fix}(T)\neq \varnothing$ (see Section\ref{subsec:2.3}). 
\footnote{By extending the idea in \cite{Haug}, another algorithm,
which we refer to as the \emph{generalized Haugazeau's algorithm},
was developed for minimizing a \emph{strictly convex} function
in $\Gamma_{0}({\cal H})$ over the fixed point set of a certain
quasi-nonexpansive operator \cite{Com2000}.
In particular, this algorithm was specialized in a clear way for finding
the nearest fixed point of a certain quasi-nonexpansive operator
\cite{BC2001} and applied successfully to an image recovery problem
\cite{CP2004}.
If we focus on the case of a nonstrictly convex function,
the generalized Haugazeau's algorithm is not applicable,
while some convergence theorems of the hybrid steepest descent method
suggest its sound applicability {\it provided that the gradient of the function is Lipschitzian}.}
Moreover, for Problem (\ref{typical-conv-opt}), if such a computable nonexpansive operator $T$ can be used to express 
${\cal S}_{p}$ as in (\ref{generic-view1}) but more nicely with some computable bounded linear operator 
$\Xi\in {\cal B}({\cal H},{\cal X})$, we can apply the hybrid steepest descent method to 
Problem (\ref{viscosity-conv-prob}) after translating it into 
\begin{equation}
\mbox{\rm find }z^{\star}\in \underset{z\in {\rm Fix}(T)}{\operatorname{argmin}} \ \Theta(z),
\label{trans-to-optfix}
\end{equation}
where $\Theta:=\Psi\circ \Xi$, because $\Theta\in \Gamma_{0}({\cal H})$ is certainly {\rm G\^{a}teaux}   
differentiable with Lipschitzian gradient $\nabla\Theta:z\mapsto \Xi^{*}\nabla \Psi \left(\Xi z\right)$ and 
$\Xi(z^{\star})\in {\cal X}$ is a solution of (\ref{viscosity-conv-prob}).

The goal of this paper is to demonstrate that plugging the modern proximal splitting operators into the hybrid steepest descent method is a powerful computational strategy for solving 
highly valuable hierarchical convex optimization problems (\ref{viscosity-conv-prob-1}) in Scenario 1 and Scenario 2. 
The remainder of the paper is organized as follows. 
In the next section, as preliminaries, we introduce elements of convex analysis and fixed point theoretic view of the modern proximal splitting algorithms. These include key ideas behind fixed point characterizations of ${\cal S}_{p}$ 
in Problem (\ref{viscosity-conv-prob}) as well as the hybrid steepest descent method for nonexpansive operators. 
Section \ref{sec:3} contains the main idea of the hierarchical convex optimization based 
on the hybrid steepest descent method applied to modern proximal splitting operators. In Section \ref{sec:4}, as a typical example of Scenario \ref{sce1}, we present an application of the proposed strategies to a hierarchical enhancement of \emph{the support vector machine} \cite{VL63,Cortes-Vapnik95,V98} where we demonstrate how we can compute \emph{the best linear classifier which achieves the maximal margin among all linear classifiers having least empirical hinge loss}. 
The proposed best linear classifier can be applied to general training data whether it is linearly separable or not.
In particular, for linearly separable data, the proposed best linear classifier, \emph{which does not require any parameter tuning}, is guaranteed to reproduce successfully the original support vector machine specially defined in \cite{VL63}. To the best of the authors' knowledge, such a unified generalization of original support vector machine for linearly separable data has not been achieved by previously reported SVMs
(see, e.g., \cite{Cortes-Vapnik95,Scholkopf,Bishop06,HastieTibshiraniFreidman09,Chang-Lin11,ScholkopfLV13,sergios15} and Section \ref{subsec:4.2}). In Section \ref{sec:5}, as a typical example along Scenario \ref{sce2}, 
we present an application of the proposed strategy to a hierarchical enhancement of Lasso \cite{Ro-Tibshirani96,HastieTibshiraniFreidman09}.
This enhancement is achieved by utilizing maximally the Douglas-Rachford splitting applied to a recently established proximity operator \cite{plc-perspective,plc-muller18} of a perspective function for the TREX problem \cite{Lederer15} which is certainly \emph{the state-of-the-art nonconvex formulation} for automatic sparsity control of Lasso. The proposed application can optimize further an additional convex criterion over the all solutions of the TREX problem.
 Finally, in Section \ref{sec:ConcludingRemarks}, we conclude this paper with some remarks on other possible advanced applications of the hybrid steepest descent method.

\section{Preliminary}
\label{sec:2}
 Let $\mathcal{X}$ be a real Hilbert space equipped with\footnote{
Often $\langle \cdot, \cdot \rangle_{\mathcal{X}}$ denotes $\langle \cdot, \cdot \rangle$ to explicitly describe its domain. 
} an inner product $\langle \cdot, \cdot \rangle$ and its induced norm $\| \cdot \| = \sqrt{\langle \cdot, \cdot \rangle}$, which is denoted by $(\mathcal{X}, \langle \cdot, \cdot \rangle, \|\cdot \|)$.
 Let $(\mathcal{K},\langle \cdot, \cdot \rangle_{\mathcal{K}}, \| \cdot \|_{\mathcal{K}})$ be another real Hilbert space. Let $A\colon \mathcal{X} \to \mathcal{K}$ be a bounded linear operator of which the norm is defined by $\|A\|_{\rm op}:=\sup_{x \in \mathcal{X}\colon \|x\|\leq 1} \|Ax\|_{\mathcal{K}}$. For a bounded linear operator $A\colon \mathcal{X} \to \mathcal{K}$, $A^*\colon \mathcal{K} \to \mathcal{X}$ denotes its adjoint or conjugate, i.e.,
\begin{equation*}
  (\forall (x,u) \in \mathcal{X} \times \mathcal{K})\quad \langle x,A^*u \rangle = \langle Ax,u \rangle_{\mathcal{K}}.
\end{equation*}

\subsection{Selected Elements of Convex Analysis and Optimization}\label{subsec:2.1}
For readers' convenience, we list minimum elements, in convex analysis, which will be used in the later sections (for their detailed accounts, see, e.g., \cite{BB96,bauschke11,plc-vrw05,plc-perspective,plc-muller18,Ekeland,hiriart,Rock-Wets98,Yamada09,Zalinescu}).

\noindent {\bf (Convex Set)}
A set $C \subset \mathcal{X}$ \revb{is said to be} convex if $\lambda x + (1-\lambda) y \in C$ for all $\lambda \in (0,1)$ and for all $x, y \in C$. 

\noindent {\bf (Proper Lower Semicontinuous Convex Function; See, e.g., \cite[Chapter 9]{bauschke11}} A function $f\colon \mathcal{X} \to (-\infty,\infty]$ is said to be proper if its effective domain $\mathrm{dom}(f):=\{{x}\in  \mathcal{X}\mid f({x})<\infty\}$ is nonempty. A function $f\colon \mathcal{X} \rightarrow (-\infty,\infty]$ \revb{is said to be} lower semicontinuous if its lower level set $\text{lev}_{\leq \alpha}f:=\{x \in \mathcal{X} \mid f(x) \leq \alpha \} (\subset \mathcal{X})$ is closed for every $\alpha \in \mathbb{R}$. 
A function $f\colon \mathcal{X} \rightarrow (-\infty,\infty]$ \revb{is said to be} convex if $f(\lambda{x}+(1-\lambda){y})\leq\lambda f({x})+(1-\lambda)f({y})$ for all $\lambda\in(0,1)$ and for all ${x}, {y}\in \operatorname{dom}(f)$.
In particular, $f$ \revb{is said to be} strictly convex if $f(\lambda{x}+(1-\lambda){y})<\lambda f({x})+(1-\lambda)f({y})$ for all $\lambda\in(0,1)$ and for all ${x}, {y}\in \operatorname{dom}(f)$ such that $x \not = y$.
The set of all proper lower-semicontinuous convex functions defined over the real Hilbert space $\mathcal{X}$ is denoted by $\Gamma_0(\mathcal{X})$. 

\noindent {\bf (Coercivity and Supercoercivity; See, e.g., \cite[Chapter 11]{bauschke11})}
A function $f\colon \mathcal{X} \to (-\infty, \infty]$ is said to be coercive if
\begin{equation*}
  \|x\| \to \infty \ \Rightarrow \ f(x) \to \infty
\end{equation*}
and supercoercive if
\begin{equation*}
  \|x\| \to \infty \ \Rightarrow \ \frac{f(x)}{\|x\|} \to \infty.
\end{equation*}
Obviously, supercoercivity of $f$ implies coercivity of $f$.
Coercivity of $f \in \Gamma_0(\mathcal{X})$ implies that 
$\operatorname{lev}_{\leq \alpha} f=\{x \in \mathcal{X} \mid f(x) \leq \alpha \}$ is bounded
for every $\alpha \in \mathbb{R}$ as well as $\operatorname{argmin}_{x \in \mathcal{X}} f(x) \not = \varnothing$.
Strict convexity of $f \in \Gamma_0(\mathcal{X})$ implies that the set of minimizers is at most singleton.

\begin{fact}\label{fa:uobqegbuqeh}
  {\bf (See, e.g., \cite[Section 11.4]{bauschke11})}
    Let $f \in \Gamma_0(\mathcal{X})$, $g \in \Gamma_0(\mathcal{K})$ and $A \in \mathcal{B}(\mathcal{X}, \mathcal{K})$ such that $\operatorname{dom}(f)\cap \operatorname{dom}( g \circ A) \not = \varnothing$. Then the following conditions
  \begin{align*}
    (\text{a}) & \text{ $\operatorname{argmin}(f+g\circ A)(\mathcal{X})$ is nonempty, closed, and bounded;} \\
    (\text{b}) & \text{ $f+g\circ A$ is coercive;} \\
    (\text{c}) & \text{ $f$ is coercive, and $g$ is bounded below;} \\
    (\text{d}) & \text{ $f$ is super-coercive;}
  \end{align*}
  satisfy that $((\text{d})\text{ or }(\text{c})) \Rightarrow (\text{b})  \Rightarrow (\text{a})$.
\end{fact}

\noindent {\bf (G\^{a}teaux Differential; See, e.g., \cite[Section 2.6]{bauschke11})} Let $U$ be an open subset of $\mathcal{X}$. Then a function $f\colon U \to \mathbb{R}$ \revb{is said to be} G\^{a}teaux differentiable at $x \in U$ if there exists $a(x) \in \mathcal{X}$ such that
\begin{eqnarray}
  \olabel{eq:GTA}
  \lim_{\delta \to 0}\frac{f(x+ \delta h) - f(x)}{\delta}=\langle a(x), h\rangle \quad  (\forall h \in \mathcal{X}).
\end{eqnarray} In this case, $\nabla f(x):= a(x)$ is called G\^{a}teaux gradient (or gradient) of $f$ at $x$.
Let $f \in \Gamma_{0}(\mathcal{X})$ be G\^{a}teaux differentiable at $x_{\star} \in \mathcal{X}$. Then $x_{\star}$ is a minimizer of $f$ if and only if $\nabla f(x_{\star})=0$. 

\noindent {\bf (Subdifferential; See, e.g., \cite[Chapter 16]{bauschke11})} For a function $f\in \Gamma_0(\mathcal{X})$, the subdifferential of $f$ is defined as the set valued operator
\begin{eqnarray*}
  \partial f\colon \mathcal{X} \to 2^{\mathcal{X}}:
  x \mapsto \{u \in \mathcal{X} \mid \langle y - x, u \rangle + f(x) \leq f(y), \forall y \in \mathcal{X} \}. 
\end{eqnarray*}
Every element $u \in \partial f(x)$ is called a subgradient of $f$ at $x$.
For a given function $f \in \Gamma_0(\mathcal{X})$, $x_{\star} \in \mathcal{X}$ is a minimizer of $f$ if and only if $0 \in \partial f(x_{\star})$. Note that if $f \in \Gamma_0(\mathcal{X})$ is G\^{a}teaux differentiable at $x \in \mathcal{X}$, then $\partial f(x):=\{\nabla f(x)\}$. 

\noindent {\bf (Conjugate Function; See, e.g., \cite[Chapter 13 and Chapter 16]{bauschke11})} For a function $f \in \Gamma_0(\mathcal{X})$, the conjugate of $f$ is defined by
\begin{equation*}
  f^*\colon \mathcal{X} \to [-\infty, \infty]: u \mapsto \sup_{x \in \mathcal{X}} ( \langle x, u \rangle - f(x)) = \sup_{x \in \operatorname{dom}(f)} ( \langle x, u \rangle - f(x)).
\end{equation*}
\noindent Let $f \in \Gamma_0(\mathcal{X})$. Then $f^* \in \Gamma_0(\mathcal{X})$ and $f^{**}=f$ are guaranteed. Moreover, we have 
\begin{equation*}
  (\forall (x, u)  \in \mathcal{X} \times \mathcal{X}) \quad  u \in \partial f(x)
  \Leftrightarrow
  f(x)+f^*(u) = \langle x,u \rangle
  \Leftrightarrow x \in \partial f^*(u),
\end{equation*}
which implies that $(\partial f)^{-1}(u) := \{x \in \mathcal{X} \mid u \in \partial f(x) \} = \partial f^*(u)$ and 
$(\partial f^*)^{-1}(x) := \{u \in \mathcal{X} \mid x \in \partial f^*(u) \} = \partial f(x)$. Often
\begin{align}
  \label{eq:FYidentity}
  (\forall x \in \operatorname{dom}(\partial f))(\forall u \in \partial f(x)) \quad  f(x)+f^*(u) = \langle x,u \rangle
\end{align}
is referred to as Fenchel-Young identity.

For hierarchical enhancement of Lasso in Section \ref{sec:5}, we exploit the following nontrivial example.
\begin{bold-example}\label{def-prop-perspective}(\rev{Subdifferential of Perspective; See \cite[Lemma 2.3]{plc-muller18}}) \\
  For given supercoercive $\varphi\in \Gamma_{0}({\mathbb R}^{N})$,
  the function
\begin{equation}
\widetilde{\varphi}:{\mathbb R}\times {\mathbb R}^{N}\rightarrow (-\infty, \infty]: (\eta,{\mathbf y})\mapsto 
\left\{
\begin{array}{ll}
\eta \varphi({\mathbf y}/\eta),
&\mbox{ if }\eta>0;\\
{\displaystyle \sup_{{\mathbf x}\in {\rm dom} (\varphi)}
\left[
\varphi({\mathbf x}+{\mathbf y})-\varphi({\mathbf x})
\right]},&\mbox{ if }\eta=0;\\
+\infty,&\mbox{ otherwise.} 
\end{array}
\right.\label{def-perspective}
\end{equation}
satisfies $\widetilde{\varphi}\in \Gamma_{0}\left({\mathbb R}\times {\mathbb R}^{N}\right)$ and is called the perspective of $\varphi$.

The subdifferential of $\widetilde{\varphi}$ is given by
\begin{equation}
\partial \widetilde{\varphi}(\eta, {\mathbf y})
=
\left\{
\begin{array}{ll}
\left\{
\left(
\varphi({\mathbf y}/\eta)-\langle {\mathbf y}/\eta, {\mathbf u}\rangle, {\mathbf u} 
\right)\in {\mathbb R}\times {\mathbb R}^{N}\mid {\mathbf u}\in \partial \varphi({\mathbf y}/\eta)
\right\}, &\mbox{ if }\eta>0;\\
\{(\mu, {\mathbf u})\in {\mathbb R}\times {\mathbb R}^{N}\mid \mu+\varphi^{*}({\mathbf u})\leq 0\},&\mbox{ if }\eta=0\mbox{ and }{\mathbf y}={\mathbf 0};\\
\varnothing,&\mbox{ otherwise.} 
\end{array}
\right.\label{subdiff-perspective2}
\end{equation}
\end{bold-example}

\noindent {\bf (Conical Hull, Span, Convex Sets; See, e.g., \cite[Chapter 6]{bauschke11})}
For a given nonempty set $C \subset \mathcal{X}$, $\operatorname{cone}(C):=\{\lambda x \mid \lambda >0, x \in C \}$ is called the conical hull of $C$, and ${\operatorname{span}}(C)$ denotes the intersection of all the linear subspaces of $\mathcal{X}$ containing $C$. The closure of ${\operatorname{span}}(C)$ is denoted by $\overline{\operatorname{span}}(C)$.
The strong relative interior of a convex set $C\subset \mathcal{X}$ is defined by
\begin{equation*}
  \operatorname{sri}(C):=\{x \in C \mid \operatorname{cone}(C -x) = \overline{\operatorname{span}}(C-x)\},
\end{equation*}
where
$C-x:=\{y-x \in \mathcal{X} \mid y \in C\}$. \\
\rev{Similarly, the relative interior of a convex set $C\subset \mathcal{X}$ is defined by
\begin{equation*}
  \operatorname{ri}(C):=\{x \in C \mid \operatorname{cone}(C -x) = \operatorname{span}(C-x)\}.
\end{equation*}
By $\operatorname{cone}(C -x) \subset \operatorname{span}(C-x) \subset \overline{\operatorname{span}}(C-x)$ for every $x \in C$, we have $\operatorname{sri}(C) \subset \operatorname{ri}(C)$. Moreover, $\operatorname{sri}(C) =\operatorname{ri}(C)$
if $\operatorname{span}(C-x) = \overline{\operatorname{span}}(C-x)$ for every $x \in C$, which implies 
\begin{align}
  \label{eq:srieqri}
  \operatorname{dim}(\mathcal{X})< \infty \ \Rightarrow \ \operatorname{sri}(C) =\operatorname{ri}(C).
\end{align}}
\noindent {\bf (Indicator Function)}
For a nonempty closed convex set $C \subset \mathcal{X}$, the indicator function of $C$ is defined by
\begin{equation*}
  \iota_{C}\colon \mathcal{X} \to (-\infty,\infty]: x \mapsto
  \left\{
    \begin{array}{ll}
      0, & \text{if }\  x \in C; \\
      +\infty, & \text{otherwise, }
    \end{array}
    \right.
\end{equation*}
which 
belongs to $\Gamma_0(\mathcal{X})$. 
In particular, for a closed subspace $V \subset \mathcal{X}$, 
\begin{align}
  u \in \partial \iota_V(x)
    \ \Leftrightarrow & \  \revb{x \in V \text{ and } u \in V^{\perp}:=\{y \in \mathcal{X}\mid (\forall v \in V) \ \langle v,y\rangle=0 \}}. \label{eq:subdiffwgdobueqg3}
\end{align}
Furthermore, the indicator function $\iota_{\{0\}}\in \Gamma_0(\mathcal{X})$ of $\{0\} \subset \mathcal{X}$ has the following properties: for all $x,u \in \mathcal{X}$
\begin{align}
  \label{eq:partialiota}
  \partial \iota_{\{0\}}(x)= & \ \left\{
    \begin{array}{cl}
      \mathcal{X}, & \text{if } \ x = 0; \\
      \varnothing, & \text{otherwise,}
    \end{array}
    \right. \\
  \iota_{\{0\}}^*(u)= & \ \sup_{y \in \mathcal{X}}(\langle y, u \rangle - \iota_{\{0\}}(y)) = 0, \label{eq:conjugatezero} \\
  \partial \iota_{\{0\}}^*(u)= & \ \{0\}.  
  \olabel{eq:iiiyuyvqet}
\end{align}

\noindent {\bf (Fenchel-Rockafellar Duality for Convex Optimization Problem Involving Linear Operator; See, e.g., \cite[Definition 15.19]{bauschke11})}
Let $f\in \Gamma_0(\mathcal{X})$, $g \in\Gamma_0(\mathcal{K})$, and $A \in {\cal B}(\mathcal{X},\mathcal{K})$.
The primal problem associated with the composite function $f + g \circ A$ is
\begin{equation}
  \label{eq:primal0}
\operatorname{minimize}_{x \in \mathcal{X}} \ f(x) + g(Ax),
\end{equation}
its dual problem is
\begin{equation}
  \label{eq:dual0}
  \operatorname{minimize}_{u \in \mathcal{K}} \ f^*(A^*u) + g^*(-u),
\end{equation}
$\mu := \inf_{x \in \mathcal{X}} (f(x) + g(Ax))$ is called the primal optimal value, 
and $\mu^* := \inf_{u \in \mathcal{K}} (f^*(A^*u) + g^*(-u))$ the dual optimal value.

\begin{fact}[\mbox{See, e.g., \cite[Theorem 15.23, Theorem 16.47, Corollary 16.53]{bauschke11}}]
  \label{fa:sc}
The condition 
\begin{align}
  \label{eq:qualification0}
  \left.
  \begin{array}{l}
    0 \in \text{\rm sri}(\operatorname{dom}(g) - A\operatorname{dom}(f)) \\
                   \text{($\operatorname{sri}$ can be replaced by $\operatorname{ri}$ in the case of $\operatorname{dim}(\mathcal{K})<\infty$, see \eqref{eq:srieqri})}
  \end{array}
                   \right\}
\end{align}
is the so-called \emph{qualification condition} for problem \eqref{eq:primal0}.
\begin{enumerate}
\item[\rm (a)]
  The condition \eqref{eq:qualification0} guarantees that
  the dual problem \eqref{eq:dual0} has a minimizer and satisfies
  \begin{align}
    \olabel{eq:pdprimaldualgap}
    \mu =\inf_{x \in \mathcal{X}} (f(x) + g(Ax))=-\min_{u \in \mathcal{K}}(f^*(A^* u) +g^*(-u))=-\mu^*;
  \end{align}
\item[\rm (b)]
  The condition \eqref{eq:qualification0} guarantees that
  the subdifferential of $f + g\circ A$ can be decomposed as
  \begin{equation}
    \olabel{eq:subdiffdecompose}
  \partial (f + g\circ A) = \partial f + A^* \circ (\partial g) \circ A;
\end{equation}
\item[\rm (c)] The qualification condition \eqref{eq:qualification0} with $f\equiv 0$ becomes $0 \in \operatorname{sri}(\operatorname{dom}(g) - \operatorname{ran}(A))$, where $\operatorname{ran}(A):= A(\mathcal{X}):=\{Ax \in {\cal K} \mid x \in {\cal X}\}$. Under this condition, (a), (b), and \eqref{eq:conjugatezero} guarantee
  \begin{align*}
    \left[
    \begin{array}{l}
      \mu =\inf_{x \in \mathcal{X}} g(Ax)=\inf_{x \in \mathcal{X}} \iota_{\{0\}}^*(x)+g(Ax)
      =-\min_{u \in \mathcal{K}}(\iota_{\{0\}}(A^* u) +g^*(-u)) \\
      \phantom{\mu =\inf_{x \in \mathcal{X}} g(Ax)=\inf_{x \in \mathcal{X}} \iota_{\{0\}}^*(x)+g(Ax)}
      =-\min_{u \in \mathcal{N}(A^*)}g^*(-u)
      =-\mu^*  \\    
      \partial (g \circ A) = A^* \circ \partial g \circ A.
    \end{array}
    \right.
\end{align*}
\end{enumerate}
\end{fact}

\begin{fact}[\mbox{\cite[Theorem 19.1]{bauschke11}}]
  \label{fa:pdsc}
  Suppose that $0 \in \operatorname{dom}(g)- A\operatorname{dom}(f)$
  (Note: This condition is
  not sufficient for \eqref{eq:qualification0}). Let $(x,u) \in \mathcal{X} \times \mathcal{K}$. Then the following are equivalent:
  \begin{itemize}
    \item[\rm (i)] \ \ $x$ is a solution of the primal problem \eqref{eq:primal0}, $u$ is a solution of the dual problem \eqref{eq:dual0}, and $\mu = -\mu^*$.
    \item[\rm (ii)] \ \ $A^*u \in \partial f(x)$ and $-u \in \partial g(Ax)$.
      \item[\rm (iii)] \ \ $x \in \partial f^*(A^*u) \cap A^{-1}(\partial g^*(-u))$.
  \end{itemize}
\end{fact}

\subsection{Selected Elements of Fixed Point Theory of Nonexpansive Operators for Application to Hierarchical Convex Optimization}
\label{sec:Monotone}
\rev{
  For readers' convenience, we list minimum elements in fixed point theory of nonexpansive mapping specially for application to hierarchical convex optimization in this paper
  (for their detailed accounts, see, e.g., \cite{BB96,bauschke11,plc-vrw05,plc-pesq11,plc-muller18,Frank2001,hiriart,Rock-Wets98,Takaha2,Yamada09}).
}

\noindent {\bf (Monotone Operator; See, e.g., \cite[Section 20.1]{bauschke11})} A set-valued operator $T\colon\mathcal{X} \to 2^{\mathcal{X}}$ \revb{is said to be} monotone over $S (\subset \mathcal{X})$ if 
\begin{equation*}
  (\forall x,y \in S)(\forall u \in Tx)(\forall v \in Ty) \quad \langle u - v, x - y\rangle \geq 0.
\end{equation*}
In particular, it \revb{is said to be} $\eta$-strongly monotone over $S$ if 
\begin{eqnarray*}
  (\exists \eta >0)(\forall x,y \in S)(\forall u \in Tx)(\forall v \in Ty) \quad  \langle u - v, x - y\rangle \geq \eta \|x - y\|^2.
\end{eqnarray*}

\noindent {\bf (Nonexpansive Operator; See, e.g., \cite{BB96} and \cite[Chapter 4]{bauschke11})}
An operator $T\colon \mathcal{X} \to \mathcal{X}$ \revb{is said to be} Lipschitz continuous with Lipschitz constant $\kappa>0$ (or $\kappa$-Lipschitzian) if 
\begin{eqnarray*}
  (\forall x,y \in \mathcal{X}) \quad \| Tx- Ty\| \leq \kappa \|x - y\|.
\end{eqnarray*}
In particular, an operator $T\colon \mathcal{X} \to \mathcal{X}$ \revb{is said to be} nonexpansive if it is $1$-Lipschitzian, i.e.,
\begin{eqnarray*}
  (\forall x,y \in \mathcal{X}) \quad \| Tx- Ty\| \leq \|x - y\|.
\end{eqnarray*}
A nonexpansive operator $T$ is said to be $\alpha$-averaged (or averaged with constant $\alpha$) \cite{bauschke11,BBR78} if there exist $\alpha \in (0,1)$ and a nonexpansive operator $\widehat{T}\colon \mathcal{X}\to\mathcal{X}$ such that
\begin{equation}
  \label{eq:averaged0}
  T=(1-\alpha) {\rm I} + \alpha \widehat{T},
\end{equation}
\revb{i.e., $T$ is an \emph{average} of the identity operator ${\rm I}$ and some nonexpansive operator $\widehat{T}$.}
If \eqref{eq:averaged0} holds for $\alpha=1/2$, $T$ is said to be firmly nonexpansive.
A nonexpansive operator $T$ is $\alpha$-averaged if and only if
\begin{equation}
  \label{eq:averaged1}
(\forall x,y \in \mathcal{X}) \quad \|Tx - Ty\|^2 \leq \|x - y\|^2 - \frac{1- \alpha}{\alpha}\|(x- Tx)-(y - Ty) \|^2.
\end{equation}
Suppose that a nonexpansive operator $T$ has the fixed point set $\operatorname{Fix}(T):=\{x \in \mathcal{X} \mid Tx =x \} \not = \varnothing$. Then $\operatorname{Fix}(T)$ can be expressed as the intersection of closed halfspaces:
  \begin{eqnarray*}
{\rm Fix}(T)={\displaystyle
\bigcap_{{\bm y}\in {\cal X}}
\left\{
{\bm x}\in {\cal X}\mid \langle {\bm y}-T({\bm y}), {\bm x}\rangle\leq
\frac{\|{\bm y}\|^{2}-\|T({\bm y})\|^{2}}{2}
\right\}
}
\end{eqnarray*}
and therefore $\operatorname{Fix}(T)$  is closed and convex (see, e.g., \cite[Proposition 5.3]{G-R}, \cite[Fact 2.1(a)]{yama2001a}, and \cite[Corollary 4.24]{bauschke11}). 
In addition, a nonexpansive operator $T$ with $\operatorname{Fix}(T)\not = \varnothing$ is said to be attracting \cite{BB96} if 
\begin{eqnarray*}
  (\forall x \not \in \operatorname{Fix}(T))(\forall z \in \operatorname{Fix}(T)) \quad  \| Tx- z\| < \|x - z\|.
\end{eqnarray*}
The condition \eqref{eq:averaged1} implies that $\alpha$-averaged nonexpansive operator $T$ is attracting if $\operatorname{Fix}(T) \not = \varnothing$.
Note that other useful properties on $\alpha$-averaged nonexpansive operators are found, e.g., in \cite{OY2001,Cegielski,combettes2015compositions}.

\begin{fact}{\bf (Krasnosel'ski\u{\i}-Mann (KM) Iteration \cite{Groetsch72} (See Also \cite[Section 5.2]{bauschke11},\cite{Mann,krasnoselskii55:_two,Dotson,reich1979weak,Cegielski}))}
\label{fa:Mann}
For a nonexpansive operator $T \colon \mathcal{X} \to \mathcal{X}$ with $\operatorname{Fix}(T) \not = \varnothing$ and any initial point $x_{0} \in \mathcal{X}$, the sequence $(x_n)_{n \in \mathbb{N}}$ generated by
\begin{equation*}
  x_{n+1} = (1-\alpha_n)x_n+\alpha_nTx_n
\end{equation*}
converges weakly\footnote{\label{fn:cluster}({Strong and weak convergences}) A sequence $(x_n)_{n \in \mathbb{N}} \subset \mathcal{X}$ is said to converge strongly to a point $x \in \mathcal{X}$ if the real number sequence $(\|x_n - x\|)_{n \in \mathbb{N}}$ converges to $0$, and to converge weakly to $x \in \mathcal{X}$ if for every $y \in \mathcal{X}$ the real number sequence $(\langle x_n - x, y \rangle)_{n \in \mathbb{N}}$ converges to $0$. If $(x_n)_{n \in \mathbb{N}}$ converges strongly to $x$, then $(x_n)_{n \in \mathbb{N}}$ converges weakly to $x$. The converse is true if $\mathcal{X}$ is finite dimensional, hence in finite dimensional case we do not need to distinguish these convergences. \\ ({Sequential cluster point}) If a sequence $(x_n)_{n \in \mathbb{N}} \subset \mathcal{X}$ possesses a subsequence that strongly (weakly) converges to a point $x \in \mathcal{X}$, then $x$ is a strong (weak) sequential cluster point of $(x_n)_{n \in \mathbb{N}}$. For weak topology of real Hilbert space in the context of Hausdorff space, see \cite[Lemma 2.30]{bauschke11}.}
 to a point in $\operatorname{Fix}(T)$ if $(\alpha_n)_{n \in \mathbb{N}} \subset [0,1]$ satisfies $\sum_{n\in \mathbb{N}} \alpha_n(1-\alpha_n) = \infty$ (Note: The weak limit of $(x_{n})_{n \in \mathbb{N}}$ depends on the choices of $x_0$ and $(\alpha_n)_{n \in \mathbb{N}}$).\footnote{
Some extensions to uniformly convex Banach spaces are found in \cite{Groetsch72,reich1979weak}.
} In particular, if $T$ is $\alpha$-averaged for some $\alpha \in (0,1)$ (see \eqref{eq:averaged0}), a simple iteration 
\begin{equation}
  \label{eq:Mannsimple}
  x_{n+1} = Tx_n = (1-\alpha) x_n + \alpha \widehat{T}x_n
\end{equation}
converges weakly to a point in $\operatorname{Fix}(T)=\operatorname{Fix}(\widehat{T})$. 
\end{fact}

\noindent {\bf (Proximity Operator \cite{Moreau62,Moreau65} (\rev{See Also} \cite[Chapter 24]{bauschke11}))} The proximity operator of $f \in \Gamma_0(\mathcal{X})$
is defined by 
\begin{equation*}
  {\rm prox}_f\colon \mathcal{X} \to \mathcal{X}: x \mapsto \underset{y \in \mathcal{X}}{\operatorname{argmin}} \ f(y) + \frac{1}{2}\|y - x \|^2.
\end{equation*}
Note that ${\rm prox}_f(x) \in \mathcal{X}$ is well defined for all $x \in \mathcal{X}$ due to the coercivity and the strict convexity of $f(\cdot)+\frac{1}{2}\|\cdot - x\|^2 \in \Gamma_{0}(\mathcal{X})$. 
It is also well known that ${\rm prox}_f$ is nothing but the resolvent of $\partial f$, i.e., ${\rm prox}_f = ({\rm I}+\partial f)^{-1}=:J_{\partial f}$, which implies that 
\begin{eqnarray}
  x \in \operatorname{Fix}({\rm prox}_f)
  \ \Leftrightarrow \
  {\rm prox}_f(x) = x \ \Leftrightarrow \ ({\rm I}+\partial f)^{-1}(x) = x \nonumber \\
  \quad \ \hspace{1.3mm}
  \ \Leftrightarrow \ x \in ({\rm I}+\partial f)(x)
  \ \Leftrightarrow \ 0 \in \partial f(x)
  \ \Leftrightarrow \ x \in \underset{y \in \mathcal{X}}{\operatorname{argmin}} \  f(y). \label{eq:eqiv}
\end{eqnarray}
Thanks to this fact, the set of all minimizers of $f \in \Gamma_0(\mathcal{X})$ can be characterized in terms of a single-valued map, i.e., $\operatorname{prox}_f$. Moreover, since the proximity operator is $1/2$-averaged nonexpansive, i.e., $\operatorname{rprox}_f:=2{\rm prox}_f - {\rm I}$ is nonexpansive, the iteration
\begin{equation}
  \label{eq:proximaliteration}
  x_{n+1} = {\rm prox}_f(x_n)
\end{equation}
converges weakly to a point in $\operatorname{argmin}_{x\in \mathcal{X}} f(x)=\operatorname{Fix}({\rm prox}_f)$ by \eqref{eq:Mannsimple} in Fact~\ref{fa:Mann}. The iterative algorithm \eqref{eq:proximaliteration} is known as \emph{proximal point algorithm} \cite{Rock76} (see \eqref{prox-point-alg1}).

\noindent 
In this paper, $f \in \Gamma_0(\mathcal{X})$ is said to be \emph{proximable} if  $\operatorname{prox}_{f}$ is available as a computable operator. 
Note that if $f \in \Gamma_0(\mathcal{X})$ is proximable, so is $f^*\in\Gamma_0(\mathcal{X})$. This is verified by 
\begin{equation*}
  {\rm prox}_{f^*} = J_{\partial f^*}  = J_{(\partial f)^{-1}}= {\rm I}-J_{\partial f} = {\rm I} - {\rm prox}_f,
\end{equation*}
which is a special example of \emph{the inverse resolvent identity} \cite[Proposition 23.20]{bauschke11}.
Note that the sum of two proximable convex functions is not necessarily proximable. Moreover, for $A \in {\cal B}(\mathcal{X},\mathcal{K})$, the composition $g \circ A \in \Gamma_0(\mathcal{X})$ for a proximable function $g \in \Gamma_0(\mathcal{K})$ is not necessarily proximable.
There are many useful formula to compute the proximity operator (see, e.g., \rev{\cite[Chapter 24]{bauschke11}}, \cite{plc-pesq11}).

\begin{bold-example}
  \label{ex:proximablefunctions}
  \
  \begin{enumerate}
  \item[\rm (a)] {\bf (Indicator function; see, e.g., \cite[Example 12.25]{bauschke11})} %
    For a nonempty closed convex set $C \subset \mathcal{X}$, 
    \begin{equation*}
      (\forall x \in \mathcal{X}) \quad {\rm prox}_{\iota_C}(x)
      = \underset{y \in \mathcal{X}}{\operatorname{argmin}} \left(\iota_C(y) + \frac{1}{2}\|y - x \|^2\right) = \underset{y \in C}{\operatorname{argmin}} \ 
      \frac{1}{2}\|y - x\|^2=:{P}_C(x)
    \end{equation*}
    holds, which implies that ${\rm prox}_{\iota_C}$ is identical to the \emph{metric projection} onto $C$. 
    In particular, if $\iota_{C}$ is proximable, $C$ is said to be \emph{simple}.
  \item[\rm (b)] {\bf (Semi-orthogonal linear transform of proximable function; see, e.g., \cite[Proposition 24.14]{bauschke11} and \cite[Table 10.1]{plc-pesq11})} \\
    For $g\in \Gamma_0(\mathcal{K})$ and $A\in \mathcal{B}(\mathcal{X},\mathcal{K})$ such that $A A^*=\nu {\rm I}$ with some $\nu >0$,
    \begin{align}
      \label{eq:proxcomp}
      (\forall x \in {\cal X}) \ \operatorname{prox}_{g \circ A}(x)=
      x + \nu^{-1} A^*(\operatorname{prox}_{\nu g}(Ax) -Ax).
    \end{align}
  \item[\rm (c)] {\bf (Hinge loss function; see, e.g., \cite{Argyriou} and \rev{\cite[Example 24.36]{bauschke11})}}
    For $\gamma>0$ and
    \begin{align}\label{def-h}
      & h:{\mathbb R}\rightarrow [0,\infty): t \mapsto \max\{0, 1-t\}, \\
      \label{eq:proxhinge}
      & (\forall t \in \mathbb{R}) \quad \operatorname{prox}_{\gamma h}(t)=\min\{ t+ \gamma, \max\{t,1\}\}.
    \end{align}
    \item[\rm (d)] {\bf ($\ell_1$ norm; see, e.g., \rev{\cite{plc-vrw05,bauschke11}})}
      For $\gamma \geq 0$ and the $\ell_1$ norm $\|\cdot\|_1 \in \Gamma_0(\mathbb{R}^N)$ 
      \begin{equation*}
        (\forall {\mathbf x}=(x_1,x_2,\ldots, x_N) \in \mathbb{R}^N) \quad \|{\mathbf x}\|_1:= \sum_{j=1}^N |x_i|,
      \end{equation*}
      the $i$-th component of the proximity operator of $\gamma \|\cdot \|_1$ is given as 
      \begin{equation*}
         \quad
        (\forall \mathbf{x}=(x_1,x_2,\dots, x_N) \in \mathbb{R}^N) \quad 
        [{\rm prox}_{\gamma \| \cdot \|_1}(\mathbf{x})]_i =\left\{
          \begin{array}{ll}
            x_i - \operatorname{sgn}(x_i)\gamma,  & {\rm if} \ |x_i|>\gamma; \\
            0, & {\rm otherwise,}
          \end{array}
        \right.
      \end{equation*}
      where $\operatorname{sgn}\colon \mathbb{R} \to \mathbb{R}$ is the signum function, i.e., $\operatorname{sgn}(x)=0$ if $x=0$ and $\operatorname{sgn}(x)=x/|x|$ otherwise.
      ${\rm prox}_{\gamma \| \cdot \|_1}$ is also known as soft-thresholding \cite{D-J94,donoho95}. 
    \item[\rm (e)] {\bf (Proximity operator of \emph{perspective} of $\| \cdot \|^q$; see, e.g., \cite{plc-muller18})}
      Let $\beta>0$ and $q>1$. The perspective $\widetilde{\varphi}_{q}$ of
      $\varphi_{q}(\cdot):=\|\cdot\|^{q}/\beta$ 
      (see also \eqref{def-perspective} in Example \ref{def-prop-perspective})
      is given by
\begin{equation}\label{def-perspective-phi2}
\widetilde{\varphi}_{q}:{\mathbb R}\times {\mathbb R}^{N}\rightarrow (-\infty, \infty]: (\eta,{\mathbf y})\mapsto 
\left\{
\begin{array}{ll}
\frac{\|{\mathbf y}\|^{q}}
{\beta \eta^{q-1}},
&\mbox{if }\eta>0;\\
0,&\mbox{if }\eta=0\mbox{ and }{\mathbf y}={\mathbf 0};\\
+\infty,&\mbox{otherwise,} 
\end{array}
\right. \\
\end{equation}
and its proximity operator can be expressed as
\begin{eqnarray}
{\rm prox}_{\widetilde{\varphi}_{q}}&:&{\mathbb R}\times {\mathbb R}^{N}\rightarrow {\mathbb R}\times {\mathbb R}^{N}\nonumber \\
&:& 
 (\eta, {\mathbf y})\mapsto 
 \left\{
 \begin{array}{ll}
 \left(\eta+\frac{\varrho}{q^{*}}\|{\mathbf p}\|^{q^{*}},{\mathbf y}-{\mathbf p}\right),&
\mbox{ if }q^{*}\eta+\varrho \|{\mathbf y}\|^{q^{*}}>0;\\
(0,{\mathbf 0}),&\mbox{ if }q^{*}\eta+\varrho \|{\mathbf y}\|^{q^{*}}\leq 0,
  \end{array}
                  \right.\olabel{cl-ex-prox-perspective} \\
  \text{where }&&
                  q^{*}:=\frac{q}{q-1}, \ \varrho:=\left(\beta(1-1/q^{*})\right)^{q^{*}-1}, \ 
                  {\mathbf p}:=\left\{
\begin{array}{ll}
\tau\frac{\mathbf y}{\|{\mathbf y}\|}, &\mbox{ if }{\mathbf y}\neq {\mathbf 0};\\
{\mathbf 0}, &\mbox{ if }{\mathbf y}={\mathbf 0}, 
\end{array}
\right. \nonumber 
\end{eqnarray}
and $\tau\in (0,\infty)$ is uniquely determined as the solution to the equation: 
\begin{equation*}
\tau^{2q^{*}-1}+\frac{q^{*}\eta}{\varrho}\tau^{q^{*}-1}+\frac{q^{*}\|{\mathbf y}\|}{\varrho^{2}}=0. 
\end{equation*}
The the proximity operator of
the translation, by $(a,{\mathbf b}) \in \mathbb{R} \times \mathbb{R}^N$, of $\widetilde{\varphi}_{q}$
\begin{align*}
  & \ \mbox{}_{\tau_{(a,{\mathbf b})}}  \widetilde{\varphi}_{q}\colon \mathbb{R} \times \mathbb{R}^N \to \mathbb{R} \times \mathbb{R}^N\colon (\eta,{\mathbf y}) \mapsto \widetilde{\varphi}_{q}\left(\eta-a, {\mathbf y}-{\mathbf b}\right),
\end{align*}
which
can be expressed as
\begin{align*}
& \ {\rm prox}_{\mbox{}_{\tau_{(a,{\mathbf b})}}  \widetilde{\varphi}_{q}}:{\mathbb R}\times {\mathbb R}^{N}\rightarrow {\mathbb R}\times {\mathbb R}^{N}:  (\eta, {\mathbf y})\mapsto  
      (a, {\mathbf b})+{\rm prox}_{\widetilde{\varphi}_{q}}\left({\eta}-a, {\mathbf y}-{\mathbf b}\right), 
\end{align*}
will play an important role in Section \ref{sec:5}.
  \end{enumerate}
\end{bold-example}

\begin{fact}\label{fa:MoreauEnvelope}{\bf (Moreau Envelope (See, e.g., \cite[Section 12.4]{bauschke11}, \cite{Moreau62,Moreau65}))} \\
  For $f \in \Gamma_0(\mathcal{X})$, 
  \begin{equation*}
    \mbox{}^{\gamma} f\colon \mathcal{X} \to \mathbb{R}: x \mapsto \min_{y \in \mathcal{X}}\left(f(y) +\frac{1}{2\gamma}\| x - y \|^2\right)
  \end{equation*}
  is called the \emph{Moreau envelope} (or \emph{Moreau-Yosida regularization)} \cite{Moreau62,Moreau65} of $f$ of the index $\gamma >0$.
  The function $\mbox{}^{\gamma}f$ is G\^{a}teaux differentiable convex with Lipschitzian gradient 
  \begin{equation*}
    \nabla \mbox{}^{\gamma}f \colon {\cal X} \to {\cal X}\colon x \mapsto \frac{1}{\gamma} ({\rm I} - {\rm prox}_{\gamma f})(x).
  \end{equation*}
  The Moreau envelope of $f$ converges pointwise to $f$ on $\operatorname{dom}(f)$ as $\gamma \downarrow 0$ (see, e.g., \cite[Proposition 12.33(ii)]{bauschke11}), i.e., $\lim_{\gamma \downarrow 0} \mbox{}^{\gamma}f(x)=f(x) \ (\forall x \in \operatorname{dom}(f))$.
\end{fact}  

\subsection{Proximal Splitting Algorithms and Their Fixed Point Characterizations}\label{subsec:2.2}

In this section, we introduce the Douglas-Rachford splitting method\footnote{
  See \cite{plc-pesq11,bauschke17} for the history of the Douglas-Rachford splitting method, originated from Douglas-Rachford's seminal paper \cite{DouglasRachford56} for solving matrix equations of the form $u=Ax+Bx$, where $A$ and $B$ are positive-definite matrices (see also \cite{Varga}).
  For recent applications, of the Douglas-Rachford splitting method, to image recovery, see, e.g., \cite{plc-jcp07,Chaux09,Dupe09,Durand10}, and to data sciences, see, e.g., \cite{Gandy10,Gandy11,plc-muller18}.
  Lastly, we remark that it was shown in \cite{E-B92} that \emph{the alternating direction method of multipliers (ADMM)} \cite{L-M,Gabay83,Boyd11,E-Y15,YY2017} can be seen as a dual variant of the Douglas-Rachford splitting method.
}
(see, e.g., \cite{L-M,E-B92,plc-jcp07,Comb2009,bauschke11,bauschke17}) and the linearized augmented Lagrangian method (see, e.g., \cite{yang2013linearized,YY2017}) as examples of \emph{the proximal splitting algorithms} built on computable nonexpansive operators with a great deal of potential in their applications to the hierarchical convex optimization problem. As explained briefly just after (\ref{generic-view1}--\ref{def-nexp}), these proximal splitting algorithms are essentially realized by applying Fact~\ref{fa:Mann} (see Section \ref{sec:Monotone}) to certain computable nonexpansive operators.

\footnotetext[10]{ 
  \rev{We should remark that Proposition \ref{pr:DRS} can also be reproduced from \cite[Proposition 26.1(iii) and Theorem 26.11(i)(iii)]{bauschke11} in the context of the monotone inclusion problems. For completeness, we present Proposition 1 and its proof in the scenario of convex optimization.}}
\begin{prop}[DRS Operator and Douglas-Rachford Splitting Method\rev{\footnotemark}]\label{pr:DRS}
  \ \\
  Let $(\mathcal{X}, \langle \cdot, \cdot \rangle_{\mathcal{X}}, \|\cdot\|_{\mathcal{X}})$ be a real Hilbert space and $f,g \in \Gamma_0(\mathcal{X})$. Suppose that \begin{align}
 &\operatorname{argmin}(f+g)(\mathcal{X}) \not = \varnothing, \label{eq:DRSq1}\\
 &\operatorname{argmin}(f^*+g^*\circ(-{\rm I}))(\mathcal{X}) \not = \varnothing, \label{eq:DRSq2}\\
 &\min(f+g)(\mathcal{X})=-\min(f^*+g^*\circ(-{\rm I}))(\mathcal{X}). \label{eq:DRSq3}
\end{align}
Then the DRS operator 
\begin{align}
  \label{eq:DRSoriginal}
  T_{\rm DRS}:=(2 \operatorname{prox}_f -{\rm I})\circ(2 \operatorname{prox}_g-{\rm I})
\end{align}
satisfies:
\begin{enumerate}
  \item[\rm (a)] $\operatorname{prox}_g(\operatorname{Fix}(T_{\rm DRS})) = \operatorname{argmin}(f+g)(\mathcal{X})$; 
    \item[\rm (b)] $T_{\rm DRS}$ is nonexpansive; 
    \item[\rm (c)] By using $(\alpha_n)_{n \in \mathbb{N}} \subset [0,1]$
      \revb{satisfying $\sum_{n\in \mathbb{N}} \alpha_n(1-\alpha_n) = \infty$}
      in Fact \ref{fa:Mann} (see Section~\ref{sec:Monotone}), the sequence $(y_n)_{n \in \mathbb{N}} \subset \mathcal{X}$ generated by
      \begin{align}
          \label{eq:DRSoriginalalg}
  y_{n+1} = (1-\alpha_n) y_n + \alpha_n T_{\rm DRS}(y_n)
\end{align}
converges weakly to a point in $\operatorname{Fix}(T_{\rm DRS})$. Moreover, $(\operatorname{prox}_g(y_n))_{n \in \mathbb{N}}$ converges weakly to a point in $\operatorname{argmin}(f+g)(\mathcal{X})$. 
\end{enumerate}
\end{prop}
\noindent The iterative algorithm to produce $(\operatorname{prox}_g(y_n))_{n \in \mathbb{N}}$ with \eqref{eq:DRSoriginalalg} can be seen as a simplest example of the so-called Douglas-Rachford splitting method.

The proof of Proposition \ref{pr:DRS}(a) is given in Appendix A because the conditions (\ref{eq:DRSq1}--\ref{eq:DRSq3}) are newly imposed for applications of $T_{\rm DRS}$ (in \eqref{eq:DRSoriginal})  to hierarchical convex optimizations in Theorem \ref{th:DRS} and in Theorem \ref{th:DRSave} (see Remark \ref{re:DRSbehind}(b) and Remark \ref{re:DRSavebehind}(b) in Section~\ref{sec:3.1})
and different from \cite[Condition (6)]{plc-jcp07} which is also in the context of convex optimization. Proposition \ref{pr:DRS}(b) is obvious from the properties of the proximity operator just after \eqref{eq:eqiv}.
For weak convergence of $(\operatorname{prox}_g(y_n))_{n \in \mathbb{N}}$ in Proposition \ref{pr:DRS}(c), see, e.g., \cite[Corollary 28.3(iii)]{bauschke11} while the weak convergence of $(y_n)_{n \in \mathbb{N}}$ is obvious from Fact \ref{fa:Mann}.

\emph{ The linearized augmented Lagrangian method (LALM)} seems to have been proposed originally as an algorithmic solution to the minimization of the nuclear norm of a matrix subject to a linear constraint \cite{yang2013linearized}. Inspired by the operator defined as the iterative update \cite[(3.7)]{yang2013linearized} in the method for this special convex optimization problem, we extended in \cite{YY2017} the operator to $T_{\rm LAL}$ in \eqref{eq:LALoriginal} to be applicable to the general convex optimization problem \eqref{typical-conv-opt} and showed the nonexpansiveness of $T_{\rm LAL}$ for solving efficiently the hierarchical convex optimization \eqref{viscosity-conv-prob} by plugging the extended operator $T_{\rm LAL}$ into the HSDM.

\begin{prop}[LAL Operator
  and Linearized Augmented Lagrangian Method]
  \label{pr:LAL}
  Let $(\mathcal{X}, \langle \cdot, \cdot \rangle_{\mathcal{X}}, \|\cdot\|_{\mathcal{X}})$ and $(\mathcal{K}, \langle \cdot, \cdot \rangle_{\mathcal{K}}, \|\cdot\|_{\mathcal{K}})$ be real Hilbert spaces.
  Suppose that $f \in \Gamma_0(\mathcal{X})$, $g=\iota_{\{0\}} \in \Gamma_0(\mathcal{K})$ and $A \in \mathcal{B}(\mathcal{X}, \mathcal{K})$ satisfy
\begin{align}
  & \ {\cal S}_{\rm pLAL}:=   \operatorname{argmin}(f+\iota_{\{0\}} \circ A)(\mathcal{X}) \not = \varnothing, \label{eq:LALq1}\\
  & \ {\cal S}_{\rm dLAL}:=  \operatorname{argmin}(f^*\circ A^*)(\mathcal{K}) \not = \varnothing, \label{eq:LALq2} \\
  & \ \min(f+\iota_{\{0\}} \circ A)(\mathcal{X})=-\min(f^*\circ A^*)(\mathcal{K}), \label{eq:LALq3}
\end{align}
where ${\cal S}_{\rm pLAL}$ is the solution set of the primal problem and 
${\cal S}_{\rm dLAL}$ is the solution set of the dual problem.
Define the LAL operator $T_{\rm LAL}\colon \mathcal{X} \times \mathcal{K} \to \mathcal{X} \times \mathcal{K}\colon (x,\nu) \mapsto (x_T,\nu_T)$ by
\begin{align}
  \label{eq:LALoriginal}
  \left[
  \begin{array}{l}
    x_{T} :=  \operatorname{prox}_{f}(x - A^*Ax+A^*\nu) \\
    \nu_{T} :=  \nu - Ax_{T}.
  \end{array}
  \right.
\end{align}
Then
\begin{enumerate}
\item[\rm (a)] $\operatorname{Fix}(T_{\rm LAL}) = {\cal S}_{\rm pLAL} \times {\cal S}_{\rm dLAL}$;
\item[\rm (b)] $T_{\rm LAL}$ is nonexpansive if $\|A\|_{\rm op} \leq 1$;
\item[\rm (c)]  By using $(\alpha_n)_{n \in \mathbb{N}} \subset [0,1]$ satisfying $\sum_{n\in \mathbb{N}} \alpha_n(1-\alpha_n) = \infty$ in Fact \ref{fa:Mann} (see Section~\ref{sec:Monotone}),
  the sequence $(x_n, \nu_n)_{n \in \mathbb{N}} \subset \mathcal{X} \times {\cal K}$ generated by 
    \begin{align}
      \label{eq:LALMMANN}
  (x_{n+1},\nu_{n+1}) = (1-\alpha_n) (x_n,\nu_n) + \alpha_n T_{\rm LAL}(x_n,\nu_n)
\end{align}
converges weakly to a point in ${\cal S}_{\rm pLAL} \times {\cal S}_{\rm dLAL}$ if $\|A\|_{\rm op}\leq 1$;
\item[\rm (d)] If $\|A\|_{\rm op}< 1$, the sequence $(x_n, \nu_n)_{n \in \mathbb{N}} \subset \mathcal{X} \times {\cal K}$ generated by \eqref{eq:LALMMANN} with $\alpha_n=1 \ (n \in \mathbb{N})$ converges weakly to a point in ${\cal S}_{\rm pLAL} \times {\cal S}_{\rm dLAL}$.
  \end{enumerate}
\end{prop}
\noindent The iterative algorithms, in Proposition~\ref{pr:LAL}(c) and (d), to produce $(x_n)_{n \in \mathbb{N}}$ with \eqref{eq:LALMMANN} can be seen as simplest examples of the so-called linearized augmented Lagrangian method.

The proof of Proposition \ref{pr:LAL}(a) is given in Appendix B for completeness because the conditions (\ref{eq:LALq1}--\ref{eq:LALq3}) are newly imposed for applications of $T_{\rm LAL}$ to hierarchical convex optimizations in Theorem \ref{th:HSDMLALI} and in Theorem \ref{th:HSDMLAL} (see Remark \ref{re:LALIIbehind}(b) and Remark~\ref{re:LALbehind}(a) in Section \ref{sec:3.2}) and different from \cite[(32)]{YY2017}.
For the proof of Proposition \ref{pr:LAL}(b), see \cite{YY2017}. Proposition \ref{pr:LAL}(c) is a straightforward application of Fact \ref{fa:Mann} to Proposition \ref{pr:LAL}(b). 
The proof of Proposition \ref{pr:LAL}(d) is given in Appendix B.

\begin{bold-remark}
  A primitive idea behind the update of the LAL operator $T_{\rm LAL}$ is found in minimization of the augmented Lagrangian function \cite{hestenes1969multiplier,powell69}: 
\begin{equation}
  \label{eq:augmentedLagrangianfunction}
  \mathcal{L}\colon \mathcal{X} \times \mathcal{K} \to (-\infty,\infty]: (x,\nu) \mapsto f(x) -\langle \nu, Ax \rangle_{\mathcal{K}} + \frac{1}{2}\|Ax \|_{\mathcal{K}}^2.
\end{equation}
Indeed, by introducing
\begin{align*}
  (\forall \hat{x} \in \mathcal{X})(\forall \hat{\nu} \in \mathcal{K}) \left[
  \begin{array}{l}
    {\cal L}_1^{(\hat{\nu})}\colon {\cal X} \to (-\infty,\infty]\colon x \mapsto {\cal L}(x,\hat{\nu});  \\
    {\cal L}_2^{(\hat{x})}\colon {\cal K} \to (-\infty,\infty]\colon \nu \mapsto {\cal L}(\hat{x},\nu),
  \end{array}
  \right.
\end{align*}
the zero $(x_{\star},\nu_{\star}) \in {\cal X} \times {\cal K}$ of the partial subdifferentials of \eqref{eq:augmentedLagrangianfunction} is characterized as
 \begin{align}
   \left[
   \begin{array}{l}
     0 \in \partial {\cal L}_1^{(\nu_{\star})}(x_{\star}) \\
     0 \in \partial {\cal L}_2^{(x_{\star})}(\nu_{\star})
    \end{array}
   \right]
    \ \Leftrightarrow & \ 
  \left[
    \begin{array}{l}
      0 \in \partial f(x_{\star}) -A^*\nu_{\star} +A^*(Ax_{\star}) \\
      0 = -Ax_{\star}
    \end{array} 
  \right] \nonumber \\
  \ \Leftrightarrow & \
    \left[
      \begin{array}{l}
        x_{\star}={\rm prox}_{f}(x_{\star}-A^*Ax_{\star} +A^*\nu_{\star}) \\       
        \nu_{\star} = \nu_{\star}-Ax_{\star}
    \end{array}
  \right]
   \nonumber \\
   \ \Leftrightarrow & \ (x_{\star},\nu_{\star}) \in \operatorname{Fix}(T_{\rm LAL}).                        \nonumber
 \end{align}
\end{bold-remark}

\subsection{Hybrid Steepest Descent Method}\label{subsec:2.3}
Consider the problem 
\begin{equation}
  \label{eq:problemtarget} 
  {\rm find } \ x_{\star} \in \underset{x \in \operatorname{Fix}(T)}{\operatorname{argmin}} \  \Theta(x) =: \Omega \not = \varnothing,
\end{equation}
where $\Theta \in \Gamma_0(\mathcal{H})$ is G\^{a}teaux differentiable \revb{over $T(\mathcal{H})$} and $T \colon \mathcal{H} \to \mathcal{H}$ is a nonexpansive operator with $\operatorname{Fix}(T) \not = \varnothing$. The hybrid steepest descent method (HSDM) 
\begin{eqnarray}
  \label{eq:updateHSDM}
  x_{n+1} = T(x_n) - \lambda_{n+1} \nabla \Theta (T(x_n))
\end{eqnarray}
generates a sequence $(x_n)_{n \in \mathbb{N}}$ to approximate successively a solution of Problem \eqref{eq:problemtarget}.
\begin{fact}[Hybrid Steepest Descent Method for Nonexpansive Operators] 
  \label{fa:HSDM} \ \\[-6mm]
  \begin{itemize}
  \item[I.] \ {\cite[special case of Theorems 3.2 and 3.3 for more general variational inequality problems]{yama2001a}}
      \revb{Let $T\colon \mathcal{H} \to \mathcal{H}$ be a nonexpansive mapping with $\operatorname{Fix}(T) \not = \varnothing$. Suppose that the gradient $\nabla \Theta$ is $\kappa$-Lipschitzian and $\eta$-strongly monotone over $T(\mathcal{H}):=\{T(x) \in \mathcal{H} \mid x \in \mathcal{H} \}$, which guarantees $|\Omega|=1$.
  Then, by using any sequence $(\lambda_{n+1})_{n \in \mathbb{N}}\subset [0,\infty)$ satisfying (W1) $\lim_{n \to +\infty} \lambda_n=0$, (W2) $ \sum_{n \in \mathbb{N}} \lambda_{n+1} = +\infty$, (W3) $\sum_{n \in \mathbb{N}} |\lambda_{n+1} - \lambda_{n+2}| <\infty$ [or $(\lambda_{n+1})_{n \in \mathbb{N}}\subset (0,\infty)$ satisfying (L1) $\lim_{n \to +\infty} \lambda_n=0$, (L2) $\sum_{n \in \mathbb{N}} \lambda_{n+1} = +\infty$, (L3) $\lim_{n \to +\infty} (\lambda_n - \lambda_{n+1}) \lambda_{n+1}^{-2} = 0$], the sequence $(x_n)_{n \in \mathbb{N}}\subset\mathcal{H}$ generated, for arbitrary $x_0 \in \mathcal{H}$, by (\ref{eq:updateHSDM}) converges strongly to the uniquely existing solution of Problem (\ref{eq:problemtarget}).}

\item[II.] \ \revb{(Nonstrictly convex case \cite{OY2001,OY2003a,yamada2011minimizing}) Assume that $\operatorname{dim}(\mathcal{H}) < \infty$. Suppose that (i) $T \colon \mathcal{H} \to \mathcal{H}$ is an attracting nonexpansive operator with bounded $\operatorname{Fix}(T) \not = \varnothing$, (ii) $\nabla \Theta$ is $\kappa$-Lipschitzian over $T(\mathcal{H})$, which guarantees $\Omega \not = \varnothing$. Then, by using\footnote{
      $\ell_{+}^1$ denotes the set of all summable nonnegative sequences.
      $\ell_{+}^2$ denotes the set of all square-summable nonnegative sequences.
} $(\lambda_{n+1})_{n \in \mathbb{N}} \in \ell_{+}^2 \setminus \ell_{+}^1$,  the sequence $(x_n)_{n \in \mathbb{N}}$ generated by (\ref{eq:updateHSDM}), for arbitrary $x_0 \in \mathcal{H}$,  satisfies $\lim_{n \to \infty} d_{\Omega}(x_n)=0$, where $d_{\Omega}(x_n):=\min_{y \in \Omega} \|x_n - y\|$.}

\end{itemize}
\end{fact}

\begin{bold-remark}
  \label{rem:HSDMabd}
\ 
  \begin{enumerate} 
  \item[\rm (a)] {\bf (Comparison between Fact \ref{fa:HSDM}(I) and Fact \ref{fa:Mann})} Fact \ref{fa:Mann} in Section \ref{sec:Monotone} is available for generation of a weak convergent sequence to a point in $\operatorname{Fix}(T)$, where the weak limit depends on the choices of $x_0$ and $(\alpha_n)_{n \in \mathbb{N}}$. Fact \ref{fa:HSDM}(I) guarantees \emph{the strong convergence} of $(x_n)_{n \in \mathbb{N}}$ to a point in $\operatorname{Fix}(T)$, where the strong limit is optimal in $\operatorname{Fix}(T)$ because it minimizes $\Theta$ \revb{able} to be designed strategically for many applications. Note that, thanks to Fact \ref{fa:HSDM}(I), we present that the LAL operator plugged into the HSDM yields an iterative approximation, of a solution of Problem \eqref{eq:problemtarget}, whose the strong convergence is guaranteed if $\Theta$ has the strongly monotone \revb{Lipschitzian} gradient \revb{over $\mathcal{H}$} (see Theorem \ref{th:HSDMLALI} below).
    \item[\rm (b)] 
{\bf (Boundedness assumption of $\operatorname{Fix}(T)$ in Fact \ref{fa:HSDM}(II))} 
For readers who get worried about the boundedness assumption in Fact \ref{fa:HSDM}(II), we present some sufficient conditions, in Section \ref{sec:3.3}, to guarantee the boundedness for $\operatorname{Fix}(T)$ in the context of DRS operators and LAL operators. These conditions hold automatically in the application to the hierarchical enhancement of the Lasso, in Section \ref{subsec:5.2}.
However, the boundedness assumption in Fact \ref{fa:HSDM}(II) may not be restrictive for most practitioners by just modifying our original target \eqref{eq:problemtarget} into
\begin{align}
  \label{eq:minHmodify}
  \operatorname{minimize}  \Theta(x) \ \text{\rm subject to } x\in \overline{B}(0,r)\cap \operatorname{Fix}(T) \not = \varnothing
\end{align}
with a sufficiently large closed ball $\overline{B}(0,r)$. Note that Fact \ref{fa:HSDM}(II) is applicable to \eqref{eq:minHmodify} because ${P}_{\overline{B}(0,r)}\circ T$ is nonexpansive and satisfies
$\operatorname{Fix}({P}_{\overline{B}(0,r)}\circ T) = \overline{B}(0,r) \cap \operatorname{Fix}(T)$ (see \cite[Proposition 1(d)]{yo2004a}).
\revb{Similar} strategy will be utilized in the application to the hierarchical enhancement of the SVM in Section \ref{subsec:4.2}.
\item[\rm (c)] {\bf (Conditions for $\Theta$)}
  The condition for $\Theta \in \Gamma_0(\mathcal{H})$ in \eqref{eq:problemtarget}, where it is required to have the Lipschitzian gradient $\nabla \Theta$, may not be restrictive as well for practitioners just by passing through the smooth regularizations, e.g., Moreau-Yosida regularization (see \eqref{viscosity-conv-prob-2} and Fact \ref{fa:MoreauEnvelope} in Section \ref{sec:Monotone}).
    \end{enumerate}
\end{bold-remark}

\begin{bold-remark}[On the Hybrid Steepest Descent Method]
  \ 
\begin{enumerate}
\item[\rm (a)] The HSDM was established originally as a generalization of the so-called Halpern-type iteration (or anchor method) \cite{hal,lions,Bac96b} for iteratively computing ${P}_{\operatorname{Fix}(T)}(x)$ for a nonexpansive operator $T\colon \mathcal{H} \to \mathcal{H}$ and $x \in \mathcal{H}$. Indeed, by choosing $\Psi(\cdot):=\frac{1}{2}\|\cdot - x\|^2$, the iteration \eqref{eq:updateHSDM} is reduced to the Halpern-type iteration.
\item[\rm (b)] One can relax (L3) to
$\lim_{n \to \infty}\frac{\lambda_n}{\lambda_{n+1}} =1$ in \cite{Xu}. Moreover, if $T$ is an averaged nonexpansive operator it was shown in \cite{Iemoto} that only (W1) and (W2) are required to guarantee the strong convergence.
\item[\rm (c)] The HSDM can be robustified against the numerical errors produced possibly in the computation of $T$ \cite{y-o-s2002}. 
\item[\rm (d)] Parallel versions of the HSDM were developed in \cite{Tak-Yam2008}. Specifically, convex optimization over the Cartesian product of the intersections of the fixed point sets of nonexpansive operators is considered, where strong convergence theorems are established under a certain contraction assumption with respect to the weighted maximum norm. 
\item[\rm (e)] The HSDM has been extended for the variational inequality problems over the fixed point set of certain class of quasi-nonexpansive operators including subgradient projection operators \cite{yo2004a,yamada2011minimizing} and has been applied to signal processing problems (see, e.g., \cite{yamada2011minimizing,OY2015}). 
\item[\rm (f)] The mathematical properties of the HSDM, e.g., in \cite{yama2001a,yo2004a} have been studied extensively in various directions by many mathematicians (see, e.g., \cite{mainge,Chidume} for extensions in Banach spaces).
\end{enumerate}
\end{bold-remark}

\section{Hierarchical Convex Optimization with Proximal Splitting Operators}\label{sec:3}
In this section, we present our central strategy for iterative approximation of the solution of the hierarchical convex optimization \eqref{viscosity-conv-prob} by plugging proximal splitting operators into the HSDM. For simplicity, we focus on the DRS and the LAL operators as such proximal splitting operators.\footnote{In \cite[Section 17.5]{yamada2011minimizing}, the authors introduced briefly the central strategy of plugging the Douglas-Rachford splitting operator into the HSDM for hierarchical convex optimization. For applications of the HSDM to other proximal splitting operators, e.g., the forward-backward splitting operator \cite{plc-vrw05}, the primal-dual splitting operator \cite{Condat13,Vu13} for the hierarchical convex optimization of different types from \eqref{viscosity-conv-prob}, see \cite{yamada2011minimizing,OYTSP}.} Assume that Problem \eqref{viscosity-conv-prob} has a solution, i.e., there exists at least one minimizer of $\Psi$ over ${\cal S}_p$, and that $(f,g,A)$ satisfies its qualification condition \eqref{eq:qualification0}
(Note: The condition \eqref{eq:qualification0} holds automatically for many instances of \eqref{typical-conv-opt}, see, e.g., Section~\ref{subsec:4.2} [\eqref{eq:qualSVMred}]
and Section~\ref{subsec:5.2} [Lemma \ref{pr:TREXqualcond} and \eqref{eq:qualTRXred}]).
As explained briefly just around (\ref{trans-to-optfix}) in Section \ref{sec:1}, for applications of the HSDM \eqref{eq:updateHSDM} to Problem \eqref{viscosity-conv-prob}, we need characterization of the constraint set as ${\cal S}_{p}=\Xi(\operatorname{Fix}(T))$ with a computable nonexpansive operator $T\colon {\cal H} \to {\cal H}$ and with a bounded linear operator $\Xi\in {\cal B}({\cal H},{\cal X})$ which ensures the G\^{a}teaux differentiability of $\Theta :=\Psi \circ \Xi \in \Gamma_0({\cal H})$ with Lipschitzian gradient $\nabla \Theta$.
In the following, we introduce three examples of such pair of computable nonexpansive operator $T\colon {\cal H} \to {\cal H}$ and $\Xi \in \mathcal{B}({\cal H},{\cal X})$.

\subsection{Plugging DRS Operators into Hybrid Steepest Descent Method}\label{sec:3.1}

We introduce a nonexpansive operator \revb{called ${\mathbf T}_{\rm DRS_{I}}$ of Type-I, as an instance of the DRS operator,} that can characterize ${\cal S}_{p}$ (see \eqref{eq:DRScharacterization}) and demonstrate how this nonexpansive operator can be plugged into the HSDM for \eqref{viscosity-conv-prob}.

\begin{theorem}[HSDM with the DRS Operator in Product Space of Type-I]
  \label{th:DRS}
  Let $f\in \Gamma_0( \mathcal{X})$, $g\in \Gamma_0( \mathcal{K})$, and $A \in \mathcal{B}(\mathcal{X}, \mathcal{K})$ in Problem \eqref{viscosity-conv-prob} satisfy ${\cal S}_p \not = \varnothing$ and the qualification condition \eqref{eq:qualification0}. 
  Suppose that $\Psi \in \Gamma_0(\mathcal{X})$ is \revb{G\^{a}teaux} differentiable with Lipschitzian gradient $\nabla \Psi$ \revb{over $\mathcal{X}$} and that $\Omega:= \underset{x^{\star} \in {\cal S}_p}{\operatorname{argmin}} \ \Psi(x^{\star}) \not = \varnothing$. Then the operator
  \begin{align}
    \label{eq:DRS31}
    {\mathbf T}_{\rm DRS_{I}}\colon \mathcal{X} \times \mathcal{K} \to \mathcal{X} \times \mathcal{K}\colon (x,y) \mapsto (x_T,y_T),
  \end{align}
  where
  \begin{align} \label{eq:DRSIo}
    \left[
    \begin{array}{l}
    p=  x-A^*({\rm I} + AA^*)^{-1}  (Ax -y)\\
    (x_{1/2},y_{1/2}) =  (2p-x, 2Ap-y) \\
      (x_{T},y_{T}) =  (2\operatorname{prox}_f(x_{1/2}) -x_{1/2}, 2\operatorname{prox}_g(y_{1/2}) -y_{1/2}),
    \end{array}
    \right.
  \end{align}
can be plugged into the HSDM \eqref{eq:updateHSDM}, with any $\alpha \in (0,1)$ and any $(\lambda_{n+1})_{n \in \mathbb{N}} \in \ell_{+}^2 \setminus \ell_{+}^1$, as 
  \begin{align}
    \label{eq:HSDM311}
    \left[
    \begin{array}{l}
    (x_{n+1/2},y_{n+1/2}) =  (1-\alpha)(x_{n}, y_{n})+ \alpha {\mathbf T}_{\rm DRS_{I}}(x_{n}, y_{n})  \\
    x_{n+1}^{\star}=  x_{n+1/2}-A^*({\rm I} + AA^*)^{-1}  (Ax_{n+1/2} -y_{n+1/2}) \\
    x_{n+1} =  x_{n+1/2} - \lambda_{n+1} ({\rm I}- A^*({\rm I}+AA^*)^{-1}A)\circ \nabla \Psi (x_{n+1}^{\star})\\
        y_{n+1} =  y_{n+1/2} - \lambda_{n+1} (({\rm I}+AA^*)^{-1}A)\circ \nabla \Psi (x_{n+1}^{\star}).
    \end{array}
    \right.
  \end{align}
  The algorithm \eqref{eq:HSDM311} generates,
  for any  $(x_0, y_0) \in \mathcal{X} \times \mathcal{K}$,
  a sequence $(x_{n+1}^{\star})_{n \in \mathbb{N}}\subset \mathcal{X}$ which satisfies
\begin{align}
  \label{eq:DRSconvergence}
  \lim_{n \to \infty} d_{\Omega}(x_{n}^{\star})=0
\end{align}
if \ $\operatorname{dim}(\mathcal{X} \times \mathcal{K}) < \infty$ and $\operatorname{Fix}({\mathbf T}_{\rm DRS_{I}})$ is bounded.
\end{theorem}
\begin{bold-remark}[Idea Behind the Derivation of Theorem \ref{th:DRS}]
  \ 
  \label{re:DRSbehind}
  \begin{enumerate}
  \item[\rm (a)] The operator ${\mathbf T}_{\rm DRS_{I}}$ in \eqref{eq:DRS31} can be expressed as\footnote{
      \rev{The use of the DRS operator in a product space as in \eqref{eq:DRSIimplicit} is found explicitly or implicitly in various applications, mainly for solving \eqref{multiple-typical-conv-opt} (see, e.g., \cite{C-P-invP-08,Gandy10,Gandy11,Pustelnik11,Dupe12,Chaari14,plc-pesq15}). 
      }      
}      
      \begin{align}
        \label{eq:DRSIimplicit}
        \hspace{-5mm} {\mathbf T}_{\rm DRS_{I}}=(2\operatorname{prox}_F -{\rm I}) \circ (2\operatorname{prox}_{\iota_{\mathcal{N}(\check{A})}} -{\rm I}) =(2\operatorname{prox}_F -{\rm I}) \circ (2{P}_{\mathcal{N}(\check{A})} -{\rm I})
  \end{align}
  which is nothing but the DRS operator in the sense of Proposition \ref{pr:DRS} (see Section \ref{subsec:2.2})
  specialized for  
  \begin{align}
    \label{eq:DRSreformulation}
    \operatorname{minimize} \ (F+\iota_{\mathcal{N}(\check{A})})( \mathcal{X} \times \mathcal{K}),
  \end{align}
  where
  \begin{align}
 \label{eq:F31}
  F\colon & \ \mathcal{X} \times \mathcal{K} \to (-\infty, \infty]\colon (x,y) \mapsto f(x)+g(y), \\
  \label{eq:checkA}
  \check{A}\colon & \ \mathcal{X} \times \mathcal{K} \to \mathcal{K}\colon (x,y) \mapsto Ax-y,
  \end{align}
  and $\mathcal{N}(\check{A})$ stands for the null space of $\check{A} \in {\cal B}(\mathcal{X} \times \mathcal{K}, \mathcal{K} )$.
  Note that exactly in the same way as in \eqref{eq:intro_proxdecompose}, $\operatorname{prox}_F\colon{\cal X} \times {\cal K} \to {\cal X} \times {\cal K}\colon(x,y) \mapsto(\operatorname{prox}_f(x),\operatorname{prox}_g(y))$ can be used as a computational tool if 
  $\operatorname{prox}_f$ and $\operatorname{prox}_g$ are available. Moreover, $\operatorname{prox}_{\iota_{\mathcal{N}(\check{A})}}={P}_{\mathcal{N}(\check{A})}\colon{\cal X} \times {\cal K} \to {\cal N}(\check{A})\colon (x,y) \mapsto (p,Ap)$ is also available if  $p$ in \eqref{eq:DRSIo} is computable, hence Problem \eqref{eq:DRSreformulation} is minimization of the sum of two proximable functions.
Obviously, Problem \eqref{eq:DRSreformulation} is a reformulation of Problem \eqref{viscosity-conv-prob} in a higher dimensional space in the sense of 
  \begin{align}
    & \ {\cal S}_{p}\mbox{[in \eqref{viscosity-conv-prob}]} =  \mathcal{Q}_{\mathcal{X}}\left[
      \underset{(x,y) \in {\cal X} \times {\cal K}}{\operatorname{argmin}}(F(x,y)+\iota_{\mathcal{N}(\check{A})}(x,y))
      \right], \olabel{eq:DRSI1}
  \end{align}
 where
  \begin{align}
  \label{eq:Q}
    \mathcal{Q}_{\mathcal{X}}\colon & \ \mathcal{X} \times \mathcal{K} \to \mathcal{X}\colon (x,y) \mapsto x, 
  \end{align}
which is verified by
  \begin{align}
  & \ \operatorname{argmin}_{x \in \mathcal{X}} f(x)+g(Ax)  \nonumber \\
  = & \ \mathcal{Q}_{\mathcal{X}}\left[\operatorname{argmin}_{(x,y) \in \mathcal{X} \times \mathcal{K}} f(x)+g(y)+\iota_{\{0\}}(Ax - y)\right] \nonumber \\
    = & \ \mathcal{Q}_{\mathcal{X}}\left[\operatorname{argmin}_{(x,y) \in \mathcal{X} \times \mathcal{K}} F(x,y)+\iota_{\mathcal{N}(\check{A})}(x,y)\right].
        \nonumber
\end{align}

\item[\rm (b)] For application of the HSDM (based on Fact \ref{fa:HSDM}(II) in Section \ref{subsec:2.3}), Theorem \ref{th:DRS} uses the convenient expression:
  \begin{align}
    & \ {\cal S}_{p}\mbox{[in \eqref{viscosity-conv-prob}]}  \overset{\mathrm{see \ below}}{=}  {\cal Q}_{\mathcal{X}}(\operatorname{prox}_{\iota_{\mathcal{N}(\check{A})}}(\operatorname{Fix}({\mathbf T}_{\rm DRS_{I}})) \label{eq:DRSI2}  \\
    & \ = \mathcal{Q}_{\mathcal{X}}({P}_{\mathcal{N}(\check{A})}(\operatorname{Fix}({\mathbf T}_{\rm DRS_{I}})) \label{eq:DRScharacterization0} \\
    & \ =\Xi_{\rm DRS_I}(\operatorname{Fix}({\mathbf T}_{\rm DRS_{I}})) =\Xi_{\rm DRS_I}(\operatorname{Fix}((1-\alpha){\rm I}+\alpha {\mathbf T}_{\rm DRS_{I}}))
      \label{eq:DRScharacterization}
  \end{align}
  in terms of attracting operator $(1-\alpha){\rm I}+\alpha {\mathbf T}_{\rm DRS_{I}}$ with $\alpha \in (0,1)$ (see \eqref{eq:averaged1}), where
  \begin{align}
    \label{eq:XIDRSI}
    \Xi_{\rm DRS_I}:= & \ \mathcal{Q}_{\mathcal{X}}\circ  {P}_{\mathcal{N}(\check{A})} \in \mathcal{B}({\cal X} \times {\cal K}, {\cal X}).
  \end{align}
Note that the characterization \eqref{eq:DRScharacterization} is illustrated in Figure \ref{fig:TREX} (see Section \ref{subsec:5.1})
and is utilized, in Section \ref{subsec:5.2}, in the context of the hierarchical enhancement of Lasso.
  To prove \eqref{eq:DRSI2} based on Proposition \ref{pr:DRS}(a) in Section \ref{subsec:2.2}, we need: \\[2mm]
  {\bf Claim \ref{th:DRS}:} {\rm
  If $f\in \Gamma_0( \mathcal{X})$, $g\in \Gamma_0( \mathcal{K})$, and $A \in \mathcal{B}(\mathcal{X}, \mathcal{K})$ in Problem \eqref{viscosity-conv-prob} satisfy ${\cal S}_p \not = \varnothing$ and the qualification condition \eqref{eq:qualification0}, we have
 \begin{align}
 &\operatorname{argmin}(F+\iota_{\mathcal{N}(\check{A})})(\mathcal{X} \times \mathcal{K}) \not = \varnothing, \label{eq:DRSq11}\\
 &\operatorname{argmin}(F^*+\iota_{\mathcal{N}(\check{A})}^*\circ(-{\rm I}))(\mathcal{X} \times \mathcal{K}) \not = \varnothing, \label{eq:DRSq12}\\
 &\min(F+\iota_{\mathcal{N}(\check{A})})(\mathcal{X} \times \mathcal{K})=-\min(F^*+\iota_{\mathcal{N}(\check{A})}^*\circ(-{\rm I}))(\mathcal{X} \times \mathcal{K}). \label{eq:DRSq13}
 \end{align}}
\ \\[-5mm]
Note that (\ref{eq:DRSq11}--\ref{eq:DRSq13}) correspond to (\ref{eq:DRSq1}--\ref{eq:DRSq3}) in Proposition \ref{pr:DRS} for minimization of $F+ \iota_{{\cal N}(\check{A})}$ and therefore Claim \ref{th:DRS} is the main step in the proof of Theorem~\ref{th:DRS}.
 \\
 \item[\rm (c)]  To plug the operator ${\mathbf T}_{\rm DRS_{I}}\colon \mathcal{H} \to \mathcal{H}$, with $\mathcal{H} := \mathcal{X} \times \mathcal{K}$, into the HSDM based on Fact \ref{fa:HSDM}(II) in Section \ref{subsec:2.3}, the characterization  
  ${\cal S}_{p}=\Xi_{\rm DRS_I}(\operatorname{Fix}((1-\alpha){\rm I}+\alpha {\mathbf T}_{\rm DRS_{I}}))$ in \eqref{eq:DRScharacterization}
  is utilized in the translation [exactly in the same way as in \eqref{trans-to-optfix}]:
\begin{align}
  & \ \Omega\mbox{[in Theorem \ref{th:DRS}]}=\Xi_{\rm DRS_I} (\Omega_{\rm DRS_{I}}),
                                        \label{eq:theorem31c} \\ 
& \ \text{where } 
               \Omega_{\rm DRS_{I}}:= \!
  \underset{{\mathbf z} \in \operatorname{Fix}({\mathbf T}_{\rm DRS_{I}})}{\operatorname{argmin}} \Theta_{\rm DRS_I}({\mathbf z})
  =\! \! \! \underset{{\mathbf z} \in \operatorname{Fix}((1-\alpha){\rm I} + \alpha {\mathbf T}_{\rm DRS_{I}})}{\operatorname{argmin}} \Theta_{\rm DRS_I}({\mathbf z}),
   \label{eq:problemHSDMDRS}
\end{align}
and $\Theta_{\rm DRS_I}=\Psi\circ \Xi_{\rm DRS_{I}} \in \Gamma_0(\mathcal{X} \times \mathcal{K})$.

\item[\rm (d)] 
  Application of the HSDM to \eqref{eq:problemHSDMDRS} yields
  \begin{align}
    \left[
    \begin{array}{l}
      {\mathbf z}_{n+1/2}=  [(1-\alpha){\rm I}+\alpha {\mathbf T}_{\rm DRS_{I}}]({\mathbf z}_n), \\
      {\mathbf z}_{n+1}=  {\mathbf z}_{n+1/2} - \lambda_{n+1}\nabla \Theta_{\rm DRS_I}({\mathbf z}_{n+1/2}) \\
      \phantom{{\mathbf z}_{n+1}}=  
          {\mathbf z}_{n+1/2} - \lambda_{n+1}\Xi_{\rm DRS_{I}}^* \nabla \Psi(\Xi_{\rm DRS_{I}}{\mathbf z}_{n+1/2}),
    \end{array}
          \right.
          \label{eq:DRSHSDMc}
  \end{align}
where $\Xi_{\rm DRS_{I}}^*$ is the conjugate of $\Xi_{\rm DRS_{I}}$ in \eqref{eq:XIDRSI}.
By letting ${\mathbf z}_n=:(x_n,y_n) \in \mathcal{X} \times \mathcal{K}$, ${\mathbf z}_{n+1/2}=:(x_{n+1/2},y_{n+1/2}) \in \mathcal{X} \times \mathcal{K}$, and
$x_{n+1}^{\star}:=\Xi_{\rm DRS_{I}}{\mathbf z}_{n+1/2} \in \mathcal{X}$, as well as, by noting
  \begin{align*}
   \hspace{-1cm} \Xi_{\rm DRS_{I}}^* = {P}_{\mathcal{N}(\check{A})} \circ \mathcal{Q}_{\mathcal{X}}^*  
     \colon  \mathcal{X} \to \mathcal{X} \times \mathcal{K} \colon 
     x \mapsto (({\rm I}- A^*({\rm I}+AA^*)^{-1}A)x, ({\rm I}+AA^*)^{-1}Ax),
    \end{align*}
    we can verify the equivalence between \eqref{eq:DRSHSDMc} and \eqref{eq:HSDM311}.
    
    \item[\rm (e)] Fact \ref{fa:HSDM}(II) in Section \ref{subsec:2.3} guarantees
      $\lim_{n \to \infty}d_{\Omega_{\rm DRS_{I}}}({\mathbf z}_n)=0$. Moreover, by noting that
      $\Xi_{\rm DRS_{I}} {P}_{\Omega_{\rm DRS_{I}}}({\mathbf z}_{n+1/2}) \in \Omega$ (see \eqref{eq:theorem31c}) and 
$\Omega_{\rm DRS_{I}} \subset \operatorname{Fix}((1-\alpha){\rm I}+\alpha {\mathbf T}_{\rm DRS_{I}})$ (see \eqref{eq:problemHSDMDRS}),
      \eqref{eq:DRSconvergence} is verified as
\begin{align*}
  & \ d_{\Omega}(x_{n+1}^{\star})=
  d_{\Omega}(\Xi_{\rm DRS_{I}}{\mathbf z}_{n+1/2}) \\  
    \leq & \ \|\Xi_{\rm DRS_{I}}{\mathbf z}_{n+1/2} - \Xi_{\rm DRS_{I}}{P}_{\Omega_{\rm DRS_{I}}}({\mathbf z}_{n+1/2}) \|_{\cal X} \\
    \leq & \ \|\Xi_{\rm DRS_{I}}\|_{\rm op}\|{\mathbf z}_{n+1/2} -{P}_{\Omega_{\rm DRS_{I}}}({\mathbf z}_{n+1/2}) \|_{\cal H} \\
    \leq & \ \|\Xi_{\rm DRS_{I}}\|_{\rm op} d_{\Omega_{\rm DRS_{I}}}({\mathbf z}_{n}) \to 0  \ (n \to \infty).
  \end{align*}
 \end{enumerate}
\end{bold-remark}
(The proof of Theorem \ref{th:DRS} is given in Appendix C). \\

Next, 
we introduce another nonexpansive operator \revb{called ${\mathbf T}_{\rm DRS_{II}}$ of Type-II, as an instance of the DRS operator,} that can characterize ${\cal S}_{p}$ (see \eqref{eq:DRSavecharacterization}) and demonstrate how this nonexpansive operator can be plugged into the HSDM for \eqref{viscosity-conv-prob}. The operator ${\mathbf T}_{\rm DRS_{II}}$ is designed based on Example \ref{ex:proximablefunctions}(b) in Section \ref{sec:Monotone}.

\begin{theorem}[HSDM with the DRS Operator in Product Space of Type-II]
  \label{th:DRSave}
  Let $\mathcal{K}=\mathbb{R}^m$.
  Let $f\in \Gamma_0( \mathcal{X})$, $g=\bigoplus_{i=1}^mg_i\in \Gamma_0(\mathcal{K})$, $A \colon\mathcal{X} \to \mathcal{K}\colon x \mapsto Ax=(A_1x,A_2x,\ldots, A_mx)$ with $A_i \in \mathcal{B}(\mathcal{X}, \mathbb{R}) \setminus\{0\} \ (i=1,2,\ldots,m)$ in Problem \eqref{viscosity-conv-prob} satisfy ${\cal S}_p \not = \varnothing$ and the qualification condition \eqref{eq:qualification0}.
  Suppose that $\Psi \in \Gamma_0(\mathcal{X})$ is \revb{G\^{a}teaux} differentiable with Lipschitzian gradient $\nabla \Psi$ \revb{over $\mathcal{X}$} and that $\Omega:= \underset{x^{\star} \in {\cal S}_p}{\operatorname{argmin}} \  \Psi(x^{\star}) \not = \varnothing$.
  Then the operator 
  \begin{align}
    \label{eq:DRSave}
    \hspace{-6mm} {\mathbf T}_{\rm DRS_{II}}\colon \mathcal{X}^{m+1}\to \mathcal{X}^{m+1}\colon (x^{(1)}, x^{(2)}, \ldots, x^{(m+1)}) \mapsto (x_T^{(1)}, x_T^{(2)}, \ldots, x_T^{(m+1)}),
  \end{align}
where
\begin{align*}
  \left[
  \begin{array}{l}
    \bar{x}=  \frac{1}{m+1}\sum_{j=1}^{m+1} x^{(j)} \\
    x_T^{(i)}= (2\bar{x}-x^{(i)}) + 2(A_iA_i^*)^{-1} A_i^*(\operatorname{prox}_{(A_iA_i^*) g_i}[A_i(2\bar{x}-x^{(i)})] -A_i(2\bar{x}-x^{(i)}))
    \quad (i=1,2,\ldots, m)  \\
    x_T^{(m+1)}=  2 \operatorname{prox}_f(2\bar{x}-x^{(m+1)})-(2\bar{x}-x^{(m+1)}),
  \end{array}
  \right.
\end{align*}
can be plugged into the HSDM \eqref{eq:updateHSDM}, with any $\alpha \in (0,1)$ and any $(\lambda_{n+1})_{n \in \mathbb{N}} \in \ell_{+}^2 \setminus \ell_{+}^1$, as 
  \begin{align}
    \label{eq:HSDMave1}
    \left[
    \begin{array}{l}
      \left(x_{n+1/2}^{(1)}, \ldots, x_{n+1/2}^{(m+1)}\right) =  (1-\alpha)\left(x_n^{(1)}, \ldots, x_n^{(m+1)}\right)
      + \alpha {\mathbf T}_{\rm DRS_{II}}\left(x_n^{(1)}, \ldots, x_n^{(m+1)}\right)  \\
    x_{n+1}^{\star}=  \frac{1}{m+1}\sum_{j=1}^{m+1} x_{n+1/2}^{(j)} \\
    x_{n+1}^{(i)} =  x_{n+1/2}^{(i)} - \frac{\lambda_{n+1}}{m+1} \nabla \Psi(x_{n+1}^{\star}) \   \quad (i=1,2,\ldots, m+1).
    \end{array}
    \right.
  \end{align}
  The algorithm \eqref{eq:HSDMave1} generates,
  for any $\left(x_0^{(1)}, \ldots, x_0^{(m+1)}\right) \in \mathcal{X}^{m+1}$, a sequence $(x_{n+1}^{\star})_{n \in \mathbb{N}} \subset \mathcal{X}$
  which satisfies
  \begin{align}
  \label{eq:DRSaveconvergence}
  \lim_{n \to \infty} d_{\Omega}(x_{n}^\star)=0
\end{align}
if  \ $\operatorname{dim}(\mathcal{X}) < \infty$ and $\operatorname{Fix}({\mathbf T}_{\rm DRS_{II}})$ is bounded.

\end{theorem}
\begin{bold-remark}[Idea Behind the Derivation of Theorem \ref{th:DRSave}]
  \ 
  \label{re:DRSavebehind}
  \begin{enumerate}
    \item[\rm (a)] The operator ${\mathbf T}_{\rm DRS_{II}}$ in \eqref{eq:DRSave} can be expressed as
      \begin{align}
        \label{eq:DRSIIexpress}
        {\mathbf T}_{\rm DRS_{II}}=(2\operatorname{prox}_H -{\rm I}) \circ (2\operatorname{prox}_{\iota_{D}} -{\rm I}) = (2\operatorname{prox}_H -{\rm I}) \circ (2{P}_D -{\rm I})
  \end{align}
  which is the DRS operator in the sense of Proposition \ref{pr:DRS} (see Section \ref{subsec:2.2}) 
  specialized for  
  \begin{align}
    \label{eq:DRSavereformulation}
    \operatorname{minimize} \ (H+\iota_{D})( \mathcal{X}^{m+1}),
  \end{align}
  where
  \begin{align}
  \label{eq:H33}
  & \ \hspace{-4mm} H\colon   \mathcal{X}^{m+1} \to (-\infty,\infty]\colon (x^{(1)},\ldots,x^{(m+1)}) \mapsto \sum_{i=1}^m g_i(A_ix^{(i)})+f(x^{(m+1)}), \\
  \label{eq:SDS}
  & \ \hspace{-4mm} D:=  \{(x^{(1)},\ldots,x^{(m+1)}) \in \mathcal{X}^{m+1} \mid x^{(i)}=x^{(j)} \ (i,j =1,2,\ldots, m+1)  \}.
  \end{align}
  Note that exactly in the same way as in \eqref{eq:intro_proxdecompose},
\begin{align*}
   & \ \operatorname{prox}_{H}(x^{(1)},x^{(2)}, \ldots, x^{(m+1)}) \\
  = & \ (\operatorname{prox}_{g_1 \circ A_1}(x^{(1)}),\ldots,\operatorname{prox}_{g_m \circ A_m}(x^{(m)}), \operatorname{prox}_f(x^{(m+1)}))
\end{align*}
 can be used with \eqref{eq:proxcomp}, in Example \ref{ex:proximablefunctions}(b) (see Section \ref{sec:Monotone}), as a computational tool if 
 $\operatorname{prox}_f$ and $\operatorname{prox}_{A_iA_i^* g} \ (i=1,2,\ldots, m)$ are available.
 Moreover, $\operatorname{prox}_{\iota_{D}}={P}_D\colon \mathcal{X}^{m+1} \to \mathcal{X}^{m+1}\colon (x^{(1)},x^{(2)}, \ldots, x^{(m+1)}) \mapsto (\bar{x},\ldots, \bar{x})$ with $\bar{x}:=\frac{1}{m+1}\sum_{i=1}^{m+1}x^{(i)}$ is also available. Hence Problem \eqref{eq:DRSavereformulation} is minimization of the sum of two proximable functions (Note: Thanks to
 $A_iA_i^* \in \mathbb{R}_{++}:=\{r \in \mathbb{R} \mid r>0 \}$, the computation of ${\mathbf T}_{\rm DRS_{II}}$ in \eqref{eq:DRSIIexpress} does not require any matrix inversion). Obviously, Problem \eqref{eq:DRSavereformulation} is a reformulation of Problem \eqref{viscosity-conv-prob} in a higher dimensional space in the sense of 
\begin{align}
    & \  {\cal S}_{p}\mbox{[in \eqref{viscosity-conv-prob}]} =  \mathcal{Q}_{\mathcal{X}^{(1)}}\left[
\underset{(x^{(1)},\ldots,x^{(m+1)}) \in \mathcal{X}^{m+1}}{\operatorname{argmin}}(H+\iota_{D})(x^{(1)},\ldots,x^{(m+1)}))\right], \olabel{eq:DRSIave1} 
\end{align}
 where
  \begin{align}
    \olabel{eq:Qnp}
  & \  \mathcal{Q}_{\mathcal{X}^{(1)}}\colon  \mathcal{X}^{m+1} \to \mathcal{X}\colon (x^{(1)},\ldots, x^{(m+1)}) \mapsto x^{(1)},
\end{align}
which is verified by
\begin{align}
  & \ \operatorname{argmin}_{x \in \mathcal{X}} g(Ax)+f(x)  \nonumber \\
  = & \ \operatorname{argmin}_{x \in \mathcal{X}} \sum_{i=1}^m g_i(A_ix)+f(x)  \nonumber \\
  = & \ \mathcal{Q}_{\mathcal{X}^{(1)}}\left[\operatorname{argmin}_{(x^{(1)},\ldots,x^{(m+1)}) \in \mathcal{X}^{m+1}} \left(H+ \iota_{D}\right)(x^{(1)},\ldots,x^{(m+1)}) \right]. \olabel{eq:uvviyicqehboqh}
\end{align}  
\item[\rm (b)] For application of the HSDM (based on Fact \ref{fa:HSDM}(II) in Section \ref{subsec:2.3}), Theorem \ref{th:DRSave} uses the convenient expression:
  \begin{align}
    & \ {\cal S}_{p}\mbox{[in \eqref{viscosity-conv-prob}]}  \nonumber \\
& \  \overset{\mathrm{see \ below}}{=}  \mathcal{Q}_{\mathcal{X}^{(1)}}(\operatorname{prox}_{\iota_{D}}(\operatorname{Fix}({\mathbf T}_{\rm DRS_{II}})) \label{eq:DRSIave2} \\
    & \ = \mathcal{Q}_{\mathcal{X}^{(1)}}({P}_{D}(\operatorname{Fix}({\mathbf T}_{\rm DRS_{II}})) \label{eq:DRSIave2b}  \\
    & \  =\Xi_{\rm DRS_{II}}(\operatorname{Fix}({\mathbf T}_{\rm DRS_{II}})) =\Xi_{\rm DRS_{II}}(\operatorname{Fix}((1-\alpha){\rm I}+\alpha {\mathbf T}_{\rm DRS_{II}}))
      \label{eq:DRSavecharacterization}
  \end{align}
  in terms of attracting operator $(1-\alpha){\rm I}+\alpha {\mathbf T}_{\rm DRS_{II}}$ with $\alpha \in (0,1)$ (see \eqref{eq:averaged1}), where
  \begin{align}
    \label{eq:XIDRSII}
    & \ \Xi_{\rm DRS_{II}}:=  \mathcal{Q}_{\mathcal{X}^{(1)}}\circ {P}_D \in \mathcal{B}({\cal X}^{m+1}, {\cal X}).
  \end{align}
  To prove \eqref{eq:DRSIave2} based on Proposition \ref{pr:DRS}(a) in Section \ref{subsec:2.2}, we need: \\[2mm]
  {\bf Claim \ref{th:DRSave}:} {\rm 
  If $\operatorname{dim}(\mathcal{K}) < \infty$, $f\in \Gamma_0( \mathcal{X})$, $g=\bigoplus_{i=1}^mg_i\in \Gamma_0(\mathcal{K})$, $A \colon\mathcal{X} \to \mathcal{K}\colon x \mapsto Ax=(A_1x,A_2x,\ldots, A_mx)$ with $A_i \in \mathcal{B}(\mathcal{X}, \mathbb{R}) \setminus\{0\} \ (i=1,2,\ldots,m)$ in Problem \eqref{viscosity-conv-prob} satisfy ${\cal S}_p \not = \varnothing$ and the qualification condition \eqref{eq:qualification0}, we have
 \begin{align}
 &\operatorname{argmin}(H+\iota_{D})(\mathcal{X}^{m+1}) \not = \varnothing, \label{eq:DRSaveq11}\\
 &\operatorname{argmin}(H^*+\iota_{D}^*\circ(-{\rm I}))(\mathcal{X}^{m+1}) \not = \varnothing, \label{eq:DRSaveq12}\\
 &\min(H+\iota_{D})(\mathcal{X}^{m+1})=-\min(H^*+\iota_{D}^*\circ(-{\rm I}))(\mathcal{X}^{m+1}). \label{eq:DRSaveq13}
 \end{align}}
\ \\[-5mm]
Note that (\ref{eq:DRSaveq11}--\ref{eq:DRSaveq13}) correspond to (\ref{eq:DRSq1}--\ref{eq:DRSq3}) in Proposition \ref{pr:DRS} for minimization of $H + \iota_{D}$ and therefore Claim \ref{th:DRSave} is the main step in the proof of Theorem~\ref{th:DRSave}.
 \item[\rm (c)]  To plug the operator ${\mathbf T}_{\rm DRS_{II}}\colon \mathcal{H} \to \mathcal{H}$, with $\mathcal{H} := \mathcal{X}^{m+1}$, into the HSDM based on Fact \ref{fa:HSDM}(II) in Section \ref{subsec:2.3}, the characterization  
  ${\cal S}_{p}=\Xi_{\rm DRS_{II}}(\operatorname{Fix}((1-\alpha){\rm I}+\alpha {\mathbf T}_{\rm DRS_{II}}))$ in \eqref{eq:DRSavecharacterization}
  is utilized in the translation [exactly in the same way as in \eqref{trans-to-optfix}]:
\begin{align}
  & \ \hspace{-7mm} \Omega\mbox{[in Theorem \ref{th:DRSave}]}=\Xi_{\rm DRS_{II}} (\Omega_{\rm DRS_{II}}),
                                       \nonumber \\
& \ \hspace{-7mm} \text{where } 
               \Omega_{\rm DRS_{II}}:=\!\!
  \underset{{\mathbf X} \in \operatorname{Fix}({\mathbf T}_{\rm DRS_{II}})}{\operatorname{argmin}} \Theta_{\rm DRS_{II}}({\mathbf X})
  =\!\!\!\!\!\!\underset{{\mathbf X} \in \operatorname{Fix}((1-\alpha){\rm I} + \alpha {\mathbf T}_{\rm DRS_{II}})}{\operatorname{argmin}} \Theta_{\rm DRS_{II}}({\mathbf X}),
   \label{eq:problemHSDMDRSave}
\end{align}
and $\Theta_{\rm DRS_{II}}=\Psi\circ \Xi_{\rm DRS_{II}} \in \Gamma_0(\mathcal{X}^{m+1})$.

\item[\rm (d)] 
  Application of the HSDM to \eqref{eq:problemHSDMDRSave} yields
  \begin{align}
    \left[
    \begin{array}{l}
      {\mathbf X}_{n+1/2}=  [(1-\alpha){\rm I}+\alpha {\mathbf T}_{\rm DRS_{II}}]({\mathbf X}_n), \\
      {\mathbf X}_{n+1}=  {\mathbf X}_{n+1/2} - \lambda_{n+1}\nabla \Theta_{\rm DRS_{II}}({\mathbf X}_{n+1/2}) \\
      \phantom{{\mathbf X}_{n+1}}=  
          {\mathbf X}_{n+1/2} - \lambda_{n+1}\Xi_{\rm DRS_{II}}^* \nabla \Psi(\Xi_{\rm DRS_{II}}{\mathbf X}_{n+1/2}),
    \end{array}
          \right.
          \label{eq:DRSHSDMavec}
  \end{align}
where $\Xi_{\rm DRS_{II}}^*$ is the conjugate of $\Xi_{\rm DRS_{II}}$ in \eqref{eq:XIDRSII}.
  By letting ${\mathbf X}_n=:(x_n^{(1)},\ldots, x_n^{(m+1)}) \in \mathcal{X}^{m+1}$, ${\mathbf X}_{n+1/2}=:(x_{n+1/2}^{(1)},\ldots, x_{n+1/2}^{(m+1)}) \in \mathcal{X}^{m+1}$, and $x_{n+1}^{\star}:=\Xi_{\rm DRS_{II}}{\mathbf X}_{n+1/2} \in \mathcal{X}$, as well as, by noting
  \begin{align*}
   \Xi_{\rm DRS_{II}}^* = {P}_D \circ \mathcal{Q}_{\mathcal{X}^{(1)}}^*
     \colon  \mathcal{X} \to \mathcal{X}^{m+1} \colon 
    x \mapsto \frac{1}{m+1}(x, x, \ldots, x),
    \end{align*}
    we can verify the equivalence between \eqref{eq:DRSHSDMavec} and \eqref{eq:HSDMave1}.
    
    \item[\rm (e)] 
In the same way as in Remark \ref{re:DRSbehind}(e), Fact \ref{fa:HSDM}(II) in Section \ref{subsec:2.3} guarantees $\lim_{n \to \infty}d_{\Omega_{\rm DRS_{II}}}({\mathbf X}_n)=0$ and \eqref{eq:DRSaveconvergence}.
 \end{enumerate}
\end{bold-remark}
(The proof of Theorem \ref{th:DRSave} is given in Appendix D).

\subsection{Plugging LAL Operator into Hybrid Steepest Descent Method}\label{sec:3.2}
\newcommand{\sca}{\mathfrak{u}}

We introduce a nonexpansive operator \revb{called ${\mathbf T}_{\rm LAL}$, as an instance of the LAL operator,} that can characterize ${\cal S}_p$ (see \eqref{eq:LALcharacterization}) and demonstrate how this nonexpansive operator can be plugged into the HSDM for \eqref{viscosity-conv-prob}.
In particular, if $\nabla \Psi$ is strongly monotone \revb{over $\mathcal{X}$}, ${\mathbf T}_{\rm LAL}$ can be plugged into the HSDM based on Fact \ref{fa:HSDM}(I) in Section \ref{subsec:2.3}, which results in a strongly convergent iterative algorithm for \eqref{viscosity-conv-prob} (see Theorem \ref{th:HSDMLALI}).
Of course, ${\mathbf T}_{\rm LAL}$ can also be plugged into the HSDM based on Fact \ref{fa:HSDM}(II) 
(see Theorem \ref{th:HSDMLAL}).

\begin{theorem}[Strong Convergence Achieved by HSDM with LAL Operator]
  \label{th:HSDMLALI}
  Let $f\in \Gamma_0( \mathcal{X})$, $g\in \Gamma_0( \mathcal{K})$ and $A \in \mathcal{B}(\mathcal{X}, \mathcal{K})$ in Problem \eqref{viscosity-conv-prob} satisfy not only ${\cal S}_p \not = \varnothing$ and the qualification condition \eqref{eq:qualification0} but also $\|\check{A}\|_{\rm op} \leq \frac{1}{\mathfrak{u}}$ $(\exists \mathfrak{u}>0)$ with $\check{A}$ in \eqref{eq:checkA}.
Suppose also that $\Psi \in \Gamma_0(\mathcal{X})$ is G\^{a}teaux differentiable with Lipschitzian as well as strongly monotone gradient $\nabla \Psi $ \revb{over $\mathcal{X}$}. 
Then the operator
  \begin{align}
    \label{eq:LAL32}
    {\mathbf T}_{\rm LAL}\colon & \ \mathcal{X} \times \mathcal{K} \times \mathcal{K} \to \mathcal{X} \times \mathcal{K} \times \mathcal{K} \\
    \colon & \ 
             \begin{pmatrix}
               x\\
               y\\
               \nu
               \end{pmatrix} \mapsto
\begin{pmatrix}
               x_T\\
               y_T\\
               \nu_T
             \end{pmatrix}
    =
  \left(
  \begin{array}{l}
    \operatorname{prox}_{f}(x -\sca^{2}(A^*Ax- A^* y) + \sca A^*\nu) \\
    \operatorname{prox}_{g}(y -\sca^{2}(-Ax+y) -\sca\nu) \\
    \nu - \sca(Ax_{T}-y_{T})
  \end{array}
\right)
\nonumber
\end{align}
can be plugged into HSDM \eqref{eq:updateHSDM}, with any $\alpha \in (0,1]$ and any $\eta_{xy},\eta_{\nu}>0$, as
  \begin{align}
    \label{eq:HSDM322I}
    \hspace{-3mm}
    \left[
    \begin{array}{l}
       (x_{n+1/2},y_{n+1/2},\nu_{n+1/2}) = (1-\alpha)(x_{n}, y_{n},\nu_n)+\alpha{\mathbf T}_{\rm LAL}(x_{n}, y_{n},\nu_n)\\
     x_{n+1} = x_{n+1/2} - \lambda_{n+1} (
                \nabla \Psi (x_{n+1/2}) + \eta_{xy} A^*(Ax_{n+1/2} -y_{n+1/2})) \\
     y_{n+1} = y_{n+1/2} + \lambda_{n+1} \eta_{xy} (Ax_{n+1/2} -y_{n+1/2})   \\
     \nu_{n+1} = \nu_{n+1/2} -\lambda_{n+1} \eta_{\nu} v_{n+1/2}.
      \end{array}
      \right.
  \end{align}
  The algorithm \eqref{eq:HSDM322I} generates,
  for any $(x_0, y_0, \nu_0) \in \mathcal{X} \times \mathcal{K} \times \mathcal{K}$, a sequence $(x_n)_{n \in \mathbb{N}} \subset \mathcal{X}$ which converges strongly to the uniquely existing solution of Problem (\ref{viscosity-conv-prob}) \revb{if $(\lambda_{n+1})_{n \in \mathbb{N}}\subset [0,\infty)$ satisfies conditions ($W1$--$W3$) [or  $(\lambda_{n+1})_{n \in \mathbb{N}}\subset (0,\infty)$ satisfies ($L1$--$L3$)]} in Fact \ref{fa:HSDM}(I) in Section \ref{subsec:2.3}.
\end{theorem}
\begin{bold-remark}[Idea Behind the Derivation of Theorem \ref{th:HSDMLALI}]
  \ 
  \label{re:LALIIbehind}
  \begin{enumerate}
    \item[\rm (a)] The operator ${\mathbf T}_{\rm LAL}$ in \eqref{eq:LAL32} can be expressed as
      \begin{align}
        \label{eq:LALexpress}
        ({\mathbf z},\nu) \mapsto ({\mathbf z}_T,\nu_T)  \text{ with }
        \left\{
        \begin{array}{l}
          {\mathbf z}_T=\operatorname{prox}_F({\mathbf z}-(\sca\check{A})^*(\sca\check{A}){\mathbf z}+(\sca\check{A})^*\nu) \\
          \nu_T=\nu - \sca\check{A}{\mathbf z}_T
        \end{array}
        \right.
  \end{align}
  by introducing ${\mathbf z}:=(x,y)$ and ${\mathbf z}_T:=(x_T,y_T)$, which is the LAL operator of Proposition \ref{pr:LAL} (see Section \ref{subsec:2.2})
  specialized for  
  \begin{align}
    \label{eq:LALreformulation}
    \operatorname{minimize} \ (F+\iota_{\{0\}}\circ (\sca\check{A}))(\mathcal{X} \times \mathcal{K}),
  \end{align}
  where
  $F$ and $\check{A}$ are defined, respectively, in \eqref{eq:F31} and in \eqref{eq:checkA}.
  Note that exactly in the same way as in \eqref{eq:intro_proxdecompose}, $\operatorname{prox}_F\colon{\cal X} \times {\cal K} \to {\cal X} \times {\cal K}\colon(x,y) \mapsto(\operatorname{prox}_f(x),\operatorname{prox}_g(y))$ can be used as a computational tool if 
  $\operatorname{prox}_f$ and $\operatorname{prox}_g$ are available. 
Obviously, Problem \eqref{eq:LALreformulation} is a reformulation of Problem \eqref{viscosity-conv-prob} in a higher dimensional space in the sense of 
\begin{align}
  & \ \hspace{-6mm} {\cal S}_{p}\mbox{[in \eqref{viscosity-conv-prob}]} =  \mathcal{Q}_{\mathcal{X}}\left[\underset{(x,y) \in {\cal X} \times {\cal K}}{\operatorname{argmin}}F(x,y)+\iota_{\{0\}}(\sca\check{A}(x,y))\right], \olabel{eq:LALI1} 
\end{align}
where $\mathcal{Q}_{\mathcal{X}}$ is defined as in \eqref{eq:Q}, which is verified by
\begin{align}
  {\cal S}_p = & \ \operatorname{argmin}_{x \in \mathcal{X}} f(x)+g(Ax)  \nonumber \\
  = & \ \mathcal{Q}_{\mathcal{X}}\left[\operatorname{argmin}_{(x,y) \in \mathcal{X} \times \mathcal{K}} f(x)+g(y)+\iota_{\{0\}}(Ax - y)\right] \nonumber \\
  = & \ \mathcal{Q}_{\mathcal{X}}\left[\operatorname{argmin}_{(x,y) \in \mathcal{X} \times \mathcal{K}} F(x,y)+\iota_{\{0\}}(\sca\check{A}(x,y))\right]. \olabel{eq:obwrhr2LAL} 
\end{align}
\item[\rm (b)] For application of the HSDM (based on Fact \ref{fa:HSDM}(I) in Section \ref{subsec:2.3}), Theorem \ref{th:HSDMLALI} uses the convenient expression:
  \begin{align}
    & \ \hspace{-6mm} {\cal S}_{p}\mbox{[in \eqref{viscosity-conv-prob}]} \nonumber\\
    & \ \hspace{-6mm} = \mathcal{Q}_{\mathcal{X}} \circ \mathcal{Q}_{\mathcal{X} \times \mathcal{K}}\left[\operatorname{argmin}(F+\iota_{\{0\}}\circ (\sca\check{A}))(\mathcal{X} \times \mathcal{K}) \times \operatorname{argmin}(F^*\circ (\sca\check{A})^*)(\mathcal{K})\right] \olabel{eq:LALI15} \\
    & \ \hspace{-8mm} \overset{\mathrm{see \ below}}{=}
      \mathcal{Q}_{\mathcal{X}}\circ \mathcal{Q}_{\mathcal{X} \times \mathcal{K}}
     (\operatorname{Fix}({\mathbf T}_{\rm LAL})) \label{eq:LALI2} \\
    & \ \hspace{-6mm} =\Xi_{\rm LAL}(\operatorname{Fix}({\mathbf T}_{\rm LAL})) 
    = \Xi_{\rm LAL}(\operatorname{Fix}((1-\alpha){\rm I}+\alpha {\mathbf T}_{\rm LAL}))
      \label{eq:LALcharacterization}
  \end{align}
  in terms of nonexpansive operator $(1-\alpha){\rm I}+\alpha {\mathbf T}_{\rm LAL}$ with $\alpha \in (0,1]$ (Note: The nonexpansiveness of ${\mathbf T}_{\rm LAL}$ is ensured by Proposition \ref{pr:LAL}(b) in Section \ref{subsec:2.2} with $\|\sca \check{A}\|_{\rm op}\leq 1$)
\begin{align}
  \olabel{eq:Qxy}
  & \ \mathcal{Q}_{\mathcal{X} \times \mathcal{K}}\colon  \mathcal{X} \times \mathcal{K} \times \mathcal{K} \to \mathcal{X} \times \mathcal{K}\colon (x,y,\nu) \mapsto (x,y) \\
    \label{eq:XILALI}
  & \
    \Xi_{\rm LAL}:= \mathcal{Q}_{\mathcal{X}}\circ \mathcal{Q}_{\mathcal{X} \times \mathcal{K}} \in \mathcal{B}({\cal X} \times {\cal K} \times {\cal K}, {\cal X}),
\end{align}
where $\mathcal{Q}_{\mathcal{X}}$ is defined as in \eqref{eq:Q}.
To prove \eqref{eq:LALI2} based on Proposition \ref{pr:LAL}(a) in Section \ref{subsec:2.2}, we need:  \\[2mm]
{\bf Claim \ref{th:HSDMLALI}:} {\rm 
If $f\in \Gamma_0( \mathcal{X})$, $g\in \Gamma_0( \mathcal{K})$ and $A \in \mathcal{B}(\mathcal{X}, \mathcal{K})$ in Problem \eqref{viscosity-conv-prob} satisfy not only ${\cal S}_p \not = \varnothing$ and the qualification condition \eqref{eq:qualification0} but also $\|\check{A}\|_{\rm op} \leq \frac{1}{\mathfrak{u}}$ $(\exists \mathfrak{u}>0)$ with $\check{A}$ in \eqref{eq:checkA}, we have
 \begin{align}
 &\operatorname{argmin}(F+\iota_{\{0\}}\circ (\sca\check{A}))(\mathcal{X} \times \mathcal{K}) \not = \varnothing, \label{eq:LALq11}\\
 &\operatorname{argmin}(F^*\circ (\sca\check{A})^*)(\mathcal{K}) \not = \varnothing, \label{eq:LALq12}\\
 &\min(F+\iota_{\{0\}}\circ (\sca\check{A}))(\mathcal{X} \times \mathcal{K})=-\min(F^*\circ (\sca\check{A})^*)(\mathcal{K}). \label{eq:LALq13}
 \end{align}}
\ \\[-5mm]
 Note that (\ref{eq:LALq11}--\ref{eq:LALq13}) correspond to (\ref{eq:LALq1}--\ref{eq:LALq3}) in Proposition \ref{pr:LAL} for minimization of $F+ \iota_{\{0\}} \circ (\sca\check{A})$ and therefore Claim \ref{th:HSDMLALI} is the main step in the proof of Theorem \ref{th:HSDMLALI}.
 In Claim \ref{th:HSDMLALI}, we also remark that $(\mathfrak{u} \check{A})^*$ in \eqref{eq:LALq12} is the conjugate of $\mathfrak{u} \check{A}$ and given by
 \begin{align}\label{eq:conjugateuAcheck}
   (\mathfrak{u} \check{A})^*\colon{\cal K} \to  {\cal X} \times {\cal K}\colon \nu \mapsto (\mathfrak{u}A^* \nu, -\mathfrak{u} \nu).
   \end{align}
 \item[\rm (c)]  
To plug the operator ${\mathbf T}_{\rm LAL}\colon \mathcal{H} \to \mathcal{H}$, with $\mathcal{H} := \mathcal{X} \times \mathcal{K} \times \mathcal{K}$, into the HSDM based on Fact \ref{fa:HSDM}(I) in Section \ref{subsec:2.3}, the characterization
 ${\cal S}_{p}=\Xi_{\rm LAL}(\operatorname{Fix}((1-\alpha){\rm I}+\alpha {\mathbf T}_{\rm LAL}))$ in \eqref{eq:LALcharacterization} 
 is utilized in the translation:
\begin{align}
  & \ \hspace{-5mm} \Omega\mbox{[in Theorem \ref{th:HSDMLALI}]} =\Xi_{\rm LAL} (\Omega_{\rm LAL}^{\rm reg}), \label{eq:omegatheorem3} \\
  & \ \hspace{-5mm} \text{where } 
               \Omega_{\rm LAL}^{\rm reg}:=\underset{{\mathbf w} \in \operatorname{Fix}({\mathbf T}_{\rm LAL})}{\operatorname{argmin}} \Theta_{\rm LAL}^{\rm reg}({\mathbf w}) = \underset{{\mathbf w} \in \operatorname{Fix}((1-\alpha){\rm I}+\alpha {\mathbf T}_{\rm LAL})}{\operatorname{argmin}} \Theta_{\rm LAL}^{\rm reg}({\mathbf w}),
   \label{eq:problemHSDMLAL2II} 
\end{align}
\ \\[-11mm]
\begin{align}
  \Theta_{\rm LAL}^{\rm reg} \colon & \mathcal{X} \times \mathcal{K} \times \mathcal{K} \to \mathbb{R} \nonumber \\
  \colon &  {\mathbf w}_{\star} \mapsto \Psi(\Xi_{\rm LAL} {\mathbf w}_{\star})+\frac{\eta_{xy}}{2}\|\check{A}\circ \mathcal{Q}_{\mathcal{X} \times \mathcal{K}}{\mathbf w}_{\star}\|^2_{\cal K}
+\frac{\eta_{\nu}}{2}\| \mathcal{Q}_{\mathcal{K}}{\mathbf w}_{\star}\|^2_{\cal K}, \nonumber 
\end{align}
for $\eta_{xy},\eta_{\nu}>0$
with $\mathcal{Q}_{\mathcal{K}}\colon \mathcal{X} \times \mathcal{K} \times \mathcal{K}\colon (x,y,\nu) \mapsto \nu$. \revb{Note that, since $\nabla \Psi$ is strongly monotone \revb{over $\mathcal{X}$}, the gradient $\nabla \Theta_{\rm LAL}^{\rm reg}$ is strongly monotone \revb{over $\mathcal{X} \times \mathcal{K} \times \mathcal{K}$}} (for the proof, see \cite[Theorem 2(d)]{YY2017}).

\item[\rm (d)] 
  Application of the HSDM to \eqref{eq:problemHSDMLAL2II} yields
  \begin{align}
    \left[
    \begin{array}{l}
       {\mathbf w}_{n+1/2} = [(1-\alpha){\rm I}+\alpha {\mathbf T}_{\rm LAL}]({\mathbf w}_{n}) \\
     {\mathbf w}_{n+1} = {\mathbf w}_{n+1/2} - \lambda_{n+1} \nabla \Theta_{\rm LAL}^{\rm reg}({\mathbf w}_{n+1/2}).
    \end{array}
    \right.
    \label{eq:HSDMLALI1}
  \end{align}
  By letting ${\mathbf w}_n:=(x_n,y_n,\nu_n) \in \mathcal{X} \times \mathcal{K} \times \mathcal{K}$ and $
{\mathbf w}_{n+1/2}:=(x_{n+1/2},y_{n+1/2},\nu_{n+1/2}) \in \mathcal{X} \times \mathcal{K} \times \mathcal{K}$, as well as, by noting 
  \begin{align}
    \label{eq:conjXILAL}
   \Xi_{\rm LAL}^* =  \mathcal{Q}_{\mathcal{X} \times \mathcal{K}}^*\circ \mathcal{Q}_{\mathcal{X}}^*
     \colon & \ \mathcal{X} \to \mathcal{X} \times \mathcal{K} \times \mathcal{K} \colon 
    x \mapsto (x, 0 ,0), \\
   (\check{A}\circ \mathcal{Q}_{\mathcal{X} \times \mathcal{K}})^*\colon & \ \mathcal{K} \to \mathcal{X} \times \mathcal{K} \times \mathcal{K}\colon y \mapsto (A^*y, -y,0), \nonumber \\
    \mathcal{Q}_{\mathcal{K}}^*\colon & \ \mathcal{K} \to \mathcal{X} \times \mathcal{K} \times \mathcal{K}\colon \nu \mapsto (0, 0, \nu), \nonumber
\end{align}
    we can verify the equivalence between \eqref{eq:HSDMLALI1} and \eqref{eq:HSDM322I}.
    
  \item[\rm (e)] Fact \ref{fa:HSDM}(I) in Section \ref{subsec:2.3} guarantees that
    $({\mathbf w}_n)_{n \in \mathbb{N}}$ converges strongly to a point in $\Omega_{\rm LAL}^{\rm reg}$. Hence,
    $(\Xi_{\rm LAL}{\mathbf w}_n(=x_n))_{n \in \mathbb{N}}$ also converges strongly to a point in $\Omega$ (see \eqref{eq:omegatheorem3}).
 \end{enumerate}
\end{bold-remark}
(The proof of Theorem \ref{th:HSDMLALI} is given in Appendix E).

\begin{theorem}[HSDM with the LAL Operator Based on Fact \ref{fa:HSDM}(II)]
  \label{th:HSDMLAL}
  Let $f\in \Gamma_0(\mathcal{X})$, $g\in \Gamma_0(\mathcal{K})$ and $A \in \mathcal{B}(\mathcal{X}, \mathcal{K})$ in Problem \eqref{viscosity-conv-prob} satisfy not only ${\cal S}_p \not =\varnothing$ and the qualification condition \eqref{eq:qualification0}  but also $\|\check{A}\|_{\rm op} \leq \frac{1}{\mathfrak{u}}$ $(\exists \mathfrak{u}>0)$ with $\check{A}$ in \eqref{eq:checkA}.
  Suppose also that $\Psi \in \Gamma_0(\mathcal{X})$ is \revb{G\^{a}teaux} differentiable with Lipschitzian gradient $\nabla \Psi$ \revb{over $\mathcal{X}$} and that $\Omega:= \underset{x^{\star} \in {\cal S}_p}{\operatorname{argmin}} \ \Psi(x^{\star}) \not = \varnothing$.
  Then the operator ${\mathbf T}_{\rm LAL}$ in \eqref{eq:LAL32} can be plugged into HSDM \eqref{eq:updateHSDM}, with any $\alpha \in (0,1)$ and any $(\lambda_{n+1})_{n \in \mathbb{N}} \in \ell_{+}^2 \setminus \ell_{+}^1$, as
  \begin{align}\label{eq:HSDMLAL31}
    \left[
    \begin{array}{l}
    (x_{n+1/2},y_{n+1},\nu_{n+1}) = (1-\alpha)(x_{n}^{\star}, y_{n},\nu_n)+ \alpha {\mathbf T}_{\rm LAL}(x_{n}^{\star}, y_{n}, \nu_n)  \\
    x_{n+1}^{\star} =  x_{n+1/2} - \lambda_{n+1} \nabla \Psi (x_{n+1/2}).
    \end{array}
\right.
  \end{align}
  The algorithm \eqref{eq:HSDMLAL31} generates,
  for any $(x_0^{\star}, y_0, \nu_0) \in \mathcal{X} \times \mathcal{K} \times \mathcal{K}$,
  a sequence $(x_n^{\star})_{n \in \mathbb{N}} \subset \mathcal{X}$
  which satisfies
  \begin{align}
  \label{eq:LALconvergence}
  \lim_{n \to \infty} d_{\Omega}(x_{n}^{\star})=0
  \end{align}
if $\operatorname{dim}(\mathcal{X} \times \mathcal{K} \times \mathcal{K}) < \infty$ and $\operatorname{Fix}({\mathbf T}_{\rm LAL})$ is bounded.

\end{theorem}
\begin{bold-remark}[Idea Behind the Derivation of Theorem \ref{th:HSDMLAL}]
  \ 
  \label{re:LALbehind}
  \begin{enumerate}
 \item[\rm (a)]  
 Following Remark \ref{re:LALIIbehind}(a)(b), we obtain the characterization 
 \begin{align*}
   {\cal S}_{p}\mbox{[in \eqref{viscosity-conv-prob}]}=\Xi_{\rm LAL}(\operatorname{Fix}((1-\alpha){\rm I}+\alpha {\mathbf T}_{\rm LAL})),
 \end{align*}
 in \eqref{eq:LALcharacterization} (see also $\Xi_{\rm LAL}$ in \eqref{eq:XILALI}), with the attracting operator $(1-\alpha){\rm I}+\alpha {\mathbf T}_{\rm LAL}$ for $\alpha \in (0,1)$ (see \eqref{eq:averaged1}).
This characterization is utilized, to plug ${\mathbf T}_{\rm LAL}\colon \mathcal{H} \to \mathcal{H}$ ($\mathcal{H} := \mathcal{X} \times \mathcal{K} \times \mathcal{K}$) into the HSDM based on Fact \ref{fa:HSDM}(II) in Section \ref{subsec:2.3}, in the translation [see also \eqref{trans-to-optfix}]:
\begin{align}
  & \ \hspace{-7mm} \Omega\mbox{[in Theorem \ref{th:HSDMLAL}]}=  \Xi_{\rm LAL} (\Omega_{\rm LAL}), \nonumber \\
& \ \hspace{-7mm} \text{where } 
               \Omega_{\rm LAL}:=
  \underset{{\mathbf w} \in \operatorname{Fix}({\mathbf T}_{\rm LAL})}{\operatorname{argmin}} \Theta_{\rm LAL}({\mathbf w})
  =\underset{{\mathbf w} \in \operatorname{Fix}((1-\alpha){\rm I} + \alpha {\mathbf T}_{\rm LAL})}{\operatorname{argmin}} \Theta_{\rm LAL}({\mathbf w})
   \label{eq:problemHSDMLAL}
\end{align}
and $\Theta_{\rm LAL}:=\Psi\circ \Xi_{\rm LAL} \in \Gamma_0(\mathcal{X} \times {\cal K} \times {\cal K})$. 

\item[\rm (b)] 
  Application of the HSDM to \eqref{eq:problemHSDMLAL} yields
  \begin{align}
    \left[
    \begin{array}{l}
      {\mathbf w}_{n+1/2}=  [(1-\alpha){\rm I}+\alpha {\mathbf T}_{\rm LAL}]({\mathbf w}_n), \\
      {\mathbf w}_{n+1}=  {\mathbf w}_{n+1/2} - \lambda_{n+1}\nabla \Theta_{\rm LAL}({\mathbf w}_{n+1/2}) \\
      \phantom{X_{n+1}}=  
          {\mathbf w}_{n+1/2} - \lambda_{n+1}\Xi_{\rm LAL}^* \nabla \Psi(\Xi_{\rm LAL}{\mathbf w}_{n+1/2}),
    \end{array}
          \right.
          \label{eq:LALHSDMavec}
  \end{align}
where $\Xi_{\rm LAL}^*$ is the conjugate of $\Xi_{\rm LAL}$.
  By letting ${\mathbf w}_n=:(x_n^{\star},y_n,\nu_n) \in \mathcal{X} \times \mathcal{K} \times \mathcal{K}$ and $
{\mathbf w}_{n+1/2}=:(x_{n+1/2},y_{n+1/2},\nu_{n+1/2}) \in \mathcal{X} \times \mathcal{K} \times \mathcal{K}$, as well as, by noting \eqref{eq:conjXILAL},
    we can verify the equivalence between \eqref{eq:LALHSDMavec} and \eqref{eq:HSDMLAL31}.
    
    \item[\rm (c)] In the same way as in Remark \ref{re:DRSbehind}(e),
Fact \ref{fa:HSDM}(II) in Section \ref{subsec:2.3} guarantees
      $\lim_{n \to \infty}d_{\Omega_{\rm LAL}}({\mathbf w}_n)=0$ and 
      \eqref{eq:LALconvergence}.
 \end{enumerate}
\end{bold-remark}
(The proof of Theorem \ref{th:HSDMLAL} is omitted, see Remark \ref{re:LALbehind}).

\subsection{Conditions for Boundedness of Fixed Point Sets of DRS and LAL Operators}\label{sec:3.3}
In Theorems \ref{th:DRS}, \ref{th:DRSave}, and \ref{th:HSDMLAL}, the boundednesses of
$\operatorname{Fix}({\mathbf T}_{\rm DRS_{I}})$, $\operatorname{Fix}({\mathbf T}_{\rm DRS_{II}})$, and $\operatorname{Fix}({\mathbf T}_{\rm LAL})$ are required for 
the algorithms \eqref{eq:HSDM311}, \eqref{eq:HSDMave1}, and \eqref{eq:HSDMLAL31} to produce $(x_{n+1}^{\star})_{n \in \mathbb{N}}$ satisfying $\lim_{n \to \infty} d_{\Omega}(x_n^{\star})=0$. Theorem \ref{th:wrhwrh} \revb{below} presents sufficient conditions for the boundednesses of these fixed point sets.
Corollary \ref{pr:uibqehrwh4} \revb{below} presents a stronger condition which will be used in Section \ref{subsec:5.2} to guarantee the boundedness of $\operatorname{Fix}({\mathbf T}_{\rm DRS_{I}})$ in the context of the hierarchical enhancement of Lasso. 
\begin{theorem}
  \label{th:wrhwrh}
  Let $f \in \Gamma_0(\mathcal{X})$, $g \in \Gamma_0(\mathcal{K})$ and $A \in \mathcal{B}(\mathcal{X}, \mathcal{K})$ in Problem \eqref{viscosity-conv-prob}
  satisfy ${\cal S}_p \not = \varnothing$ and the qualification condition \eqref{eq:qualification0}.
  Let $(A^*)^{-1}\colon {\cal X} \to 2^{\cal K}\colon x \mapsto \{y \in {\cal K} \mid  x = A^*y\}$.
  Then we have \\
{\rm (a)} \revb{\text{${\cal S}_p$, $\partial f({\cal S}_p)$, and \ 
$\bigcup_{x \in {\cal S}_p}\left([-(A^*)^{-1}(\partial f(x))] \cap \partial g(Ax)\right) \subset \mathcal{K}$
              are bounded}}
            \vspace{-2mm}
\begin{align}
  \Rightarrow \ \text{$\operatorname{Fix}\left({\mathbf T}_{\rm DRS_{I}} \right)\subset  \mathcal{X} \times \mathcal{K}$ in Theorem \ref{th:DRS} is bounded.}
    \nonumber
\end{align}
\noindent {\rm (b)} If $\check{A}$ in \eqref{eq:checkA} satisfies $\|\check{A}\|_{\rm op} \leq \frac{1}{\mathfrak{u}}$ $(\exists \mathfrak{u}>0)$, then \\[4mm]
\revb{\hspace{5mm}${\cal S}_p$ and
  $\displaystyle \bigcup_{x \in {\cal S}_p}\left([-(A^*)^{-1}(\partial f(x))] \cap \partial g(Ax)\right) \subset \mathcal{K}$ are bounded}
\vspace{-4mm}
\begin{align*}
    \Rightarrow \ \text{
  $\operatorname{Fix}\left( 
    {\mathbf T}_{\rm LAL}
  \right) \subset  \mathcal{X} \times \mathcal{K} \times \mathcal{K}$ in Theorem \ref{th:HSDMLAL} is bounded.
    }
    \nonumber
\end{align*}

\noindent {\rm (c)} If, in particular, $\mathcal{K}=\mathbb{R}^m$, $g=\bigoplus_{i=1}^mg_i\in \Gamma_0(\mathbb{R}^m)$, $A \colon \mathcal{X} \to \mathbb{R}^m\colon x \mapsto Ax=(A_1x,A_2x,\ldots, A_mx)$ with $A_i \in \mathcal{B}(\mathcal{X}, \mathbb{R})\setminus\{0\} \ ( i=1,2,\ldots,m)$ in Problem \eqref{viscosity-conv-prob}, then \\
 ${\cal S}_p$ and \
\revb{$\bigcup_{x \in {\cal S}_p }
  \left(
    \left[\varprod_{j=1}^{m} A_j^* \partial g_j (A_jx)\right] \times 
    \left[\partial f(x) \cap \left(-\sum_{i=1}^{m}  A_i^*\partial g_i( A_ix)\right)\right]
\right)$} are bounded
\vspace{-5mm}
\begin{align*}
    \Rightarrow \ \text{ 
    $\operatorname{Fix}\left( 
    {\mathbf T}_{\rm DRS_{II}}
  \right) \subset \mathcal{X}^{m+1}$
  in Theorem \ref{th:DRSave} is bounded.}
    \nonumber
\end{align*}
\end{theorem}
(The proof of Theorem \ref{th:wrhwrh} is given in Appendix F.)

The following simple relations 
\begin{align}
  & \ \begin{array}{l}
        \text{$\bigcup_{x \in {\cal S}_p}\left([-(A^*)^{-1}(\partial f(x))] \cap \partial g(Ax)\right)\text{[in Theorem \ref{th:wrhwrh}(a)(b)]}$} \nonumber \\
        \hspace{5mm}\text{$\subset  \left[ (-(A^{*})^{-1}(\partial f({\cal S}_p)))
            \cap \partial g({\cal K})\right] \text{ and }$}
        \end{array} \\
    & \ \begin{array}{l}
          \text{$\revb{\bigcup_{x \in {\cal S}_p }
          \left(
          \left[\varprod_{j=1}^{m} A_j^* \partial g_j (A_jx)\right] \times 
          \left[\partial f(x) \cap \left(-\sum_{i=1}^{m}  A_i^*\partial g_i( A_ix)\right)\right]      
          \right)\text{[in Theorem \ref{th:wrhwrh}(c)]}}$} \\
    \hspace{5mm} \revb{\subset
      \left[\varprod_{j=1}^{m} A_j^*\partial g_j(\mathbb{R})\right] \times 
          \left[-\sum_{i=1}^{m} A_i^*\partial g_i(\mathbb{R}) \right]}
        \end{array}
      \nonumber
\end{align}
lead to the corollary below.
\begin{cor}
  \label{pr:uibqehrwh4}
  Let $f \in \Gamma_0(\mathcal{X})$, $g \in \Gamma_0(\mathcal{K})$, $A \in \mathcal{B}(\mathcal{X}, \mathcal{K})$ in Problem \eqref{viscosity-conv-prob}, and $\check{A}$ in \eqref{eq:checkA} satisfy ${\cal S}_p \not = \varnothing$, the qualification condition \eqref{eq:qualification0}, and $\|\check{A}\|_{\rm op} \leq \frac{1}{\mathfrak{u}}$ $(\exists \mathfrak{u}>0)$.
  Then we have \\
  (a) ${\cal S}_p, \partial f ({\cal S}_p)$ and $(-(A^{*})^{-1}(\partial f ({\cal S}_p))) \cap \partial g(\mathcal{K}) $ are bounded
  \vspace{-7mm}
  \begin{align}
    \olabel{eq:prop3cond}
    \\
    \Rightarrow \ \left[
    \begin{array}{l}
      \text{$\operatorname{Fix}\left({\mathbf T}_{\rm DRS_{I}} \right)\subset  \mathcal{X} \times \mathcal{K}$ in Theorem \ref{th:DRS} is bounded;} \\
\text{
  $\operatorname{Fix}\left( 
    {\mathbf T}_{\rm LAL}
  \right) \subset  \mathcal{X} \times \mathcal{K} \times \mathcal{K}$ in Theorem \ref{th:HSDMLAL} is bounded.
      }
    \end{array}
    \right.
                            \nonumber 
  \end{align}
  \noindent \revb{{\rm (b)} If, in particular, $\mathcal{K}=\mathbb{R}^m$, $g=\bigoplus_{i=1}^mg_i\in \Gamma_0(\mathbb{R}^m)$, $A \colon \mathcal{X} \to \mathbb{R}^m\colon x \mapsto Ax=(A_1x,A_2x,\ldots, A_mx)$ with $A_i \in \mathcal{B}(\mathcal{X}, \mathbb{R})\setminus\{0\} \ ( i=1,2,\ldots,m)$ in Problem \eqref{viscosity-conv-prob}, then} \\
  \vspace{-5mm}
  \begin{align*}
    & \ \text{ ${\cal S}_p$ and $\partial g_i(\mathbb{R}) \ (i=1,2,\ldots,m)$ are bounded} \\
    \Rightarrow &\ \text{ 
    $\operatorname{Fix}\left( 
    {\mathbf T}_{\rm DRS_{II}}
  \right) \subset \mathcal{X}^{m+1}$
  in Theorem \ref{th:DRSave} is bounded.}
    \nonumber
\end{align*}
\end{cor}

\section{Application to Hierarchical Enhancement of Support Vector Machine}
\label{sec:4}
\subsection{Support Vector Machine}\label{subsec:4.1}

Consider a supervised learning problem for estimating a binary function 
\begin{equation}\olabel{target-function}
{\mathfrak L}:{\mathbb R}^{p}\rightarrow \{-1,1\} 
\end{equation}
with a given training dataset ${\cal D}:=\{({\mathbf x}_{i}, y_{i})\in {\mathbb R}^{p}\times \{-1,1\}
\mid i=1,2,\ldots,N \}$, where $y_{i}$ is a possibly corrupted version of the label ${\mathfrak L}({\mathbf x}_{i})$ of 
the point ${\mathbf x}_{i}$. The {\sl support vector machine} (SVM) has been recognized as one of the most successful supervised machine learning algorithms for such a learning problem. 
For simplicity, we focus on the linear SVM because 
the nonlinear SVM exploiting the so-called {\sl Kernel trick} can be viewed as an instance of the linear classifiers in the \emph{Reproducing Kernel Hilbert Spaces} (RKHS).  

The dataset ${\cal D}$ is said to be {\sl linearly separable} if there exists $({\mathbf w},b)\in \left({\mathbb R}^{p}\setminus \{{\mathbf 0}\}\right)\times {\mathbb R}$ defining a $(p-1)$-dimensional hyperplane 
\begin{equation}\label{separating-hp}
\Pi_{({\mathbf w},b)}:=\{{\mathbf x}\in {\mathbb R}^{p}\mid {\mathbf w}^{\top}{\mathbf x}-b=0\}=\Pi_{t({\mathbf w},b)}\quad (\forall t>0)
\end{equation}
which  satisfies 
\begin{equation}
\left.
\begin{array}{l}
\{{\mathbf x}_{i}\in {\mathbb R}^{p}\mid ({\mathbf x}_{i}, 1)\in {\cal D}\}
\subset \Pi_{({\mathbf w},b)}^{+}:=\{{\mathbf x}\in {\mathbb R}^{p}\mid {\mathbf w}^{\top}{\mathbf x}-b>0\}\\
\{{\mathbf x}_{i}\in {\mathbb R}^{p}\mid ({\mathbf x}_{i}, -1)\in {\cal D}\}
\subset \Pi_{({\mathbf w},b)}^{-}:=\{{\mathbf x}\in {\mathbb R}^{p}\mid {\mathbf w}^{\top}{\mathbf x}-b<0\}
\end{array}
\right\}.
  \label{halfspaces+class0}
\end{equation}

In such a case, the so-called linear classifier is defined as a mapping 
\begin{equation}\label{lin-classifier}
{\cal L}_{({\mathbf w},b)}:{\mathbb R}^{p}\rightarrow \{-1,1\}:{\mathbf x}\mapsto \left\{
\begin{array}{ll}
+1&\mbox{ if }{\mathbf x}\in \Pi_{({\mathbf w},b)}^{+},\\
-1&\mbox{ if }{\mathbf x}\in \Pi_{({\mathbf w},b)}^{-},
\end{array}
\right.
\end{equation}
which is hopefully a good approximation of the function ${\mathfrak L}$ observed partially through the training dataset ${\cal D}$. 
If ${\cal D}$ is linearly separable, there also exists infinitely many $({\mathbf w},b)\in {\mathbb R}^{p}\times {\mathbb R}$ satisfying 
\begin{equation}
\left.
\begin{array}{l}
{\cal D}_{+}:=\{{\mathbf x}_{i}\in {\mathbb R}^{p}\mid ({\mathbf x}_{i}, 1)\in {\cal D}\}
\subset \Pi_{({\mathbf w},b)}^{\geq 1}:=\{{\mathbf x}\in {\mathbb R}^{p}\mid {\mathbf w}^{\top}{\mathbf x}-b\geq 1\}\\
{\cal D}_{-}:=\{{\mathbf x}_{i}\in {\mathbb R}^{p}\mid ({\mathbf x}_{i}, -1)\in {\cal D}\}
\subset \Pi_{({\mathbf w},b)}^{\leq -1}:=\{{\mathbf x}\in {\mathbb R}^{p}\mid {\mathbf w}^{\top}{\mathbf x}-b\leq -1\}\\
\end{array}
\right\}, \label{halfspaces+class+margin}
\end{equation}
which is confirmed by rescaling $({\mathbf w},b)\in {\mathbb R}^{p}\times {\mathbb R}$  in (\ref{separating-hp}) with a constant $t\geq 1/\min\{|{\mathbf w}^{\top}{\mathbf x}_{i}-b|\}_{i=1}^{N}>0$.
 
The half-spaces $\Pi_{({\mathbf w},b)}^{\geq 1}$ and $\Pi_{({\mathbf w},b)}^{\leq -1}$ defined in (\ref{halfspaces+class+margin}) are main players in the following consideration on the linear classifier ${\cal L}_{({\mathbf w},b)}$ even for linearly non-separable data ${\cal D}$. 
In this paper, the \emph{margin} of the linear classifier ${\cal L}_{({\mathbf w},b)}$ in (\ref{lin-classifier})
is defined by   
\begin{equation}\olabel{width-gap}
\frac{1}{2}{\rm dist}\left(\Pi_{({\mathbf w},b)}^{\geq 1}, \Pi_{({\mathbf w},b)}^{\leq -1}\right)
=\frac{1}{2}\min_{\substack{
{\mathbf x}_{+}\in \Pi_{({\mathbf w},b)}^{\geq 1}, {\mathbf x}_{-}\in \Pi_{({\mathbf w},b)}^{\leq -1}}}
\|{\mathbf x}_{+}-{\mathbf x}_{-}\|=\frac{1}{\|{\mathbf w}\|}. 
\end{equation}
By using the function $h$ in \eqref{def-h} and 
\[
(\forall {\mathbf z}\in {\mathbb R}^{p})\quad 
\left[
\begin{array}{l}
{\rm d}\left({\mathbf z}, \Pi_{({\mathbf w},b)}^{\geq 1}\right)
=
\left\{
\begin{array}{ll}
\frac{\left|
{\mathbf w}^{\top}{\mathbf z}-b-1
\right|}{\|{\mathbf w}\|}=
\frac{
1-\left({\mathbf w}^{\top}{\mathbf z}-b\right)}
{\|{\mathbf w}\|}
&\mbox{ if }{\mathbf z}\not\in \Pi_{({\mathbf w},b)}^{\geq 1},\\
0&\mbox{ otherwise, }
\end{array}
\right.\\
{\rm d}\left({\mathbf z}, \Pi_{({\mathbf w},b)}^{\leq -1}\right)
=
\left\{
\begin{array}{ll}
\frac{
\left|
{\mathbf w}^{\top}{\mathbf z}-b+1
\right|}
{\|{\mathbf w}\|}
=\frac{1+\left({\mathbf w}^{\top}{\mathbf z}-b\right)}{\|{\mathbf w}\|}&\mbox{ if }{\mathbf z}\not\in 
\Pi_{({\mathbf w},b)}^{\leq -1},\\
0&\mbox{ otherwise, }
\end{array}
\right.
\end{array}
\right.
\]
we deduce 
\begin{equation}\label{hinge-dist}
\|{\mathbf w}\|
\left[
\sum_{{\mathbf z}\in {\cal D}_{+}}{\rm d}\left({\mathbf z}, \Pi_{({\mathbf w},b)}^{\geq 1}\right)
+\sum_{{\mathbf z}\in {\cal D}_{-}}
{\rm d}\left({\mathbf z}, \Pi_{({\mathbf w},b)}^{\leq -1}\right)
\right]
=\sum_{i=1}^{N}h\left(y_{i}\left({\mathbf w}^{\top}{\mathbf x}_{i }-b\right)\right)
\end{equation}
which clarifies the geometric interpretation of ``the \emph{empirical hinge loss} of  ${\cal L}_{({\mathbf w},b)}$''  defined in the right hand side of (\ref{hinge-dist}) and ensures    
\begin{equation}
\mbox{\rm Condition (\ref{halfspaces+class+margin})}\Leftrightarrow 
\sum_{i=1}^{N}h\left(y_{i}\left({\mathbf w}^{\top}{\mathbf x}_{i }-b\right)\right)=0. 
\label{halfspaces+class+margin2}
\end{equation}

For linearly separable data ${\cal D}$, among all ${\cal L}_{({\mathbf w},b)}$ satisfying (\ref{halfspaces+class+margin}), the \emph{Support Vector Machine (SVM)} ${\cal L}_{({\mathbf w^{\star}},b^{\star})}$ was proposed in 1960s by Vapnik (see, e.g.,\cite{VL63,V98}) as a special linear classifier  which achieves maximal margin, i.e.,  
\begin{equation}
\frac{1}{2}{\rm dist}\left(\Pi_{({\mathbf w}^{\star},b^{\star})}^{\geq 1}, \Pi_{({\mathbf w}^{\star},b^{\star})}^{\leq -1}\right)
=\max_{({\mathbf w},b)\mbox{ satisfies (\ref{halfspaces+class+margin2})} } 
\frac{1}{2}{\rm dist}\left(\Pi_{({\mathbf w},b)}^{\geq 1}, \Pi_{({\mathbf w},b)}^{\leq -1}\right).\label{def-orig-SVM}
\end{equation}

Therefore the SVM  ${\cal L}_{({\mathbf w^{\star}},b^{\star})}$ for linearly separable ${\cal D}$ is given as the solution of the following convex optimization problem:
\begin{eqnarray}
&&\mbox{\rm minimize }\|{\mathbf w}\|^{2}\mbox{\rm subject to }\sum_{i=1}^{N}h\left(y_{i}\left({\mathbf w}^{\top}{\mathbf x}_{i }-b\right)\right)=0 \label{SVM-separable-opt2} \\
&&\Updownarrow\nonumber \\ 
&&\hspace*{-1cm} \mbox{\rm minimize }\|{\mathbf w}\|^{2}
\mbox{\rm subject to }({\mathbf w},b)\in 
   \underset{(\hat{\mathbf w},\hat{b})\in {\mathbb R}^{p}\times {\mathbb R}}{\operatorname{argmin}}    \ 
\sum_{i=1}^{N}h\left(y_{i}\left({\hat{\mathbf w}}^{\top}{\mathbf x}_{i }-{\hat b}\right)\right), \label{SVM-separable-opt3} 
\end{eqnarray}
where the last equivalence holds true under the linear separability of ${\cal D}$ because of the nonnegativity of $h$ 
in (\ref{def-h}).  

The SVM defined equivalently in (\ref{def-orig-SVM}) or (\ref{SVM-separable-opt2}) or (\ref{SVM-separable-opt3}) for linearly separable training data has been extended for applications to even possibly linearly nonseparable training data ${\cal D}$  where the existence of $({\mathbf w},b)\in {\mathbb R}^{p}\times {\mathbb R}$ 
satisfying (\ref{halfspaces+class+margin}) is no longer guaranteed. 
One of the most widely accepted extensions of (\ref{SVM-separable-opt3}) is known as \emph{the soft margin hyperplane} 
\cite{Cortes-Vapnik95,Chang-Lin11,Bishop06,HastieTibshiraniFreidman09} which is characterized as a solution to the optimization problem:
\begin{equation}\label{SVM-nonseparable-opt1}
\mbox{\rm minimize, w.r.t. $({\mathbf w},b)$, }
\quad \frac{1}{2}\|{\mathbf w}\|^{2}+\reg
\sum_{i=1}^{N}h\left(y_{i}({\mathbf w}^{\top}{\mathbf x}_{i }-b)\right) 
\end{equation}
or equivalently
\begin{eqnarray}
&&\mbox{\rm minimize, w.r.t. $({\mathbf w},b, {\mathbf \xi})$, }\qquad \frac{1}{2}\|{\mathbf w}\|^{2}+\reg
\sum_{i=1}^{N}\xi_{i}\nonumber \\
&&\mbox{subject to }y_{i}\left({\mathbf w}^{\top}{\mathbf x}-b\right)\geq 1-\xi_{i}\mbox{ {\rm and }}\ \xi_{i}\geq 0\ (i=1,2,\ldots,N), 
\label{SVM-nonseparable-opt2}
\end{eqnarray} 
where $\reg>0$ is a tuning parameter, and $\xi_{i}$ ($i=1,2,\ldots,N$) are slack variables.

Along the Cover's theorem (on the capacity of a space in linear dichotomies) \cite{Cover95}, saying that
the probability of any grouping of the points ${\mathbf x}_1, {\mathbf x}_2, \ldots, {\mathbf x}_N, \in \mathbb{R}^l$, in general position, into two classes to be linearly separable tends to unity as $l \to \infty$,
another extension of the strategy ${\cal L}_{({\mathbf w}^{\star},b^{\star})}$ of (\ref{def-orig-SVM}) into higher dimensional spaces was made in \cite{BoserGuyonVapnik92,Cover95}, for application to possibly linearly nonseparable training data ${\cal D}$, by passing through a certain nonlinear transform $\mathfrak{N}\colon \mathbb{R}^p \to \mathbb{R}^l$ $(l \gg p)$ of the original training data ${\cal D}:=\{({\mathbf x}_{i}, y_{i})\in {\mathbb R}^{p}\times \{-1,1\} \mid i=1,2,\ldots,N \}$ to $\mathfrak{D}:=\{(\mathfrak{N}({\mathbf x}_{i}), y_{i})\in {\mathbb R}^{l}\times \{-1,1\} \mid i=1,2,\ldots,N \}$, where the nonlinear transform $\mathfrak{N}$ is defined usually in terms of kernel built in the theory of \emph{the Reproducing Kernel Hilbert Space (RKHS)} \cite{Aronszajn,Saitoh,Scholkopf} for exploiting the so-called kernel trick.

\subsection{Optimal Margin Classifier with Least Empirical Hinge Loss}\label{subsec:4.2}

As suggested in \cite[Section 3]{Cortes-Vapnik95}, the original goal behind the soft margin hyperplane in 
(\ref{SVM-nonseparable-opt1}) or (\ref{SVM-nonseparable-opt2}) seems to determine $({\mathbf w}^{\star\star},b^{\star\star}) \in (\mathbb{R}^{p} \setminus \{\mathbf 0\}) \times \mathbb{R}$ as the solution of the following nonconvex hierarchical optimization:
\begin{align}\label{min-subset1}
  \operatorname{minimize} & \ \frac{1}{2}\|{\mathbf w}^{\star} \|^2 \\
  \text{ subject to} & \ 
  ({\mathbf w}^{\star},b^{\star})\in \underset{({\mathbf w},b)}{\operatorname{argmin}} \ \left|{\cal E}({\mathbf w},b)\right|,
  \label{min-subset2}
\end{align}
where $| \cdot |$ stands for the cardinality of a set and 
${\cal E}({\mathbf w},b)\subset {\cal D}_{+}\cup {\cal D}_{-}$ is the training error set defined as  
\begin{equation}\label{min-subset}
{\cal E}({\mathbf w},b):=\left\{{\mathbf z}\in 
{\cal D}_{+}\mid {\rm d}\left({\mathbf z}, \Pi_{({\mathbf w},b)}^{\geq 1}\right)>0\right\}
\cup
\left\{{\mathbf z}\in {\cal D}_{-}\mid {\rm d}\left({\mathbf z}, \Pi_{({\mathbf w},b)}^{\leq -1}\right)>0\right\}, 
\end{equation}
i.e., determining a special hyperplane $({\mathbf w}^{\star\star},b^{\star\star})$, which achieves maximal margin in the set
$\operatorname{argmin}_{({\mathbf w},b)}\left|{\cal E}({\mathbf w},b)\right|$, is desired.

Unfortunately, since the problem to determine $({\mathbf w}^{\star},b^{\star})$ in (\ref{min-subset2}) is in general NP-\revb{hard} \cite{Judd87,BlumRivest92,Cortes-Vapnik95} and since (\ref{hinge-dist}) implies that 
\begin{equation}
\mbox{(\ref{min-subset2})}\Leftrightarrow 
({\mathbf w}^{\star},b^{\star})\in \underset{({\mathbf w},b)}{\operatorname{argmin}}
\sum_{i=1}^{N} \left[\lim_{\sigma \downarrow 0} \ h^{\sigma}\left(y_{i}({\mathbf w}^{\top}{\mathbf x}_{i }-b)\right)\right],
\nonumber
\end{equation}
the original goal set in (\ref{min-subset1}--\ref{min-subset}) was replaced in \cite[Section3]{Cortes-Vapnik95} by a realistic goal (\ref{SVM-nonseparable-opt1}) [or  (\ref{SVM-nonseparable-opt2})] for a sufficiently large constant $\reg>0$.
However, unlike the desired solution of (\ref{min-subset1}--\ref{min-subset}),
the soft margin hyperplane in (\ref{SVM-nonseparable-opt1}) 
applied to linearly separable data has no guarantee to reproduce the original SVM in (\ref{def-orig-SVM}).    

The above observations induce a natural question: \\
\emph{Is the solution of (\ref{SVM-nonseparable-opt1}) for general training data  really a mathematically sound extension of the original SVM defined equivalently in (\ref{def-orig-SVM}) or (\ref{SVM-separable-opt2}) or (\ref{SVM-separable-opt3}) specialized for  linearly separable training data ?}\footnote{
  This question is common even for the soft margin SVM applied to the transformed data $\mathfrak{D}$ employed in \cite{BoserGuyonVapnik92} because the linear separability of $\mathfrak{D}$ is not always guaranteed. 
}

Clearly, this question comes from essentially common concern as seen in Scenario \ref{sce1}, therefore, 
an alternative natural extension of  the original SVM in (\ref{SVM-separable-opt3}) would be the solution of 
the optimization problem:
\begin{equation}
\mbox{\rm minimize }\frac{1}{2}\|{\mathbf w}^{\star}\|^{2}
\mbox{ {\rm subject to }}({\mathbf w}^{\star},b^{\star})\in \Gamma:=
\underset{({\mathbf w},{b})\in {\mathbb R}^{p}\times {\mathbb R}}{\operatorname{argmin}} \ 
\sum_{i=1}^{N}h\left(y_{i}({\mathbf w}^{\top}{\mathbf x}_{i }-b)\right)
\label{SVM-nonseparable-opt3} 
\end{equation}
which does not seem different from (\ref{SVM-separable-opt3}) at a glance but is defined even possibly for linearly nonseparable training data ${\cal D}$.
Remark that the hierarchical convex optimization problem (\ref{SVM-nonseparable-opt3}) is a  more faithful convex relaxation of (\ref{min-subset1}--\ref{min-subset}) than the convex optimization (\ref{SVM-nonseparable-opt1}) [or its equivalent formulation \eqref{SVM-nonseparable-opt2} with slack variables.\footnote{
 In terms of slack variables, Problem (\ref{SVM-nonseparable-opt3}) can also be restated as~\begin{eqnarray}
&&\mbox{\rm minimize }\frac{1}{2}\|{\mathbf w}^{\star}\|^{2} \ \mbox{\rm subject to }({\mathbf w}^{\star},b^{\star},{\mathbf \xi}^{\star})\in 
\underset{({\mathbf w},b, {\mathbf \xi})\in {\mathbb R}^{p}\times {\mathbb R}\times {\mathbb R}^{N}}{\operatorname{argmin}}\ 
\sum_{i=1}^{N}\left[
\xi_{i}+\iota_{S_{i}}({\mathbf w},b, {\mathbf \xi})+\iota_{\mathfrak{\Xi}_{i}}({\mathbf w},b, {\mathbf \xi})
\right],\nonumber \\
&&\mbox{ where }S_{i}:=\left\{({\mathbf w},b, {\mathbf \xi})\in {\mathbb R}^{p}\times {\mathbb R}\times {\mathbb R}^{N}
\mid y_{i}\left({\mathbf w}^{\top}{\mathbf x}_{i}-{b}\right)\geq 1-\xi_{i}
\right\}\mbox{ and }\mathfrak{\Xi}_{i}:=\left\{({\mathbf w},b, {\mathbf \xi})\in {\mathbb R}^{p}\times {\mathbb R}\times {\mathbb R}^{N}\mid \xi_{i}\geq 0\right\}\quad (i=1,2,\ldots, N).
   \nonumber
\end{eqnarray}}] for the soft margin hyperplane. This is because the solution of (\ref{SVM-nonseparable-opt3}) for linearly separable data certainly reproduces the original SVM in (\ref{def-orig-SVM}). As remarked in Example \ref{ex:into}(b) in Section \ref{sec:1}, in general, the soft margin hyperplane via (\ref{SVM-nonseparable-opt1}) for a fixed constant $\reg>0$ does not achieve the hierarchical optimality in the sense of \eqref{SVM-nonseparable-opt3}. Fortunately, Problem (\ref{SVM-nonseparable-opt3}) \revb{falls} in the class of the hierarchical convex optimization problems of type (\ref{viscosity-conv-prob}).

In the following, we demonstrate how Problem \eqref{SVM-nonseparable-opt3} can be solved by a proposed strategy in Section \ref{sec:3}. Let ${\cal X} := \mathbb{R}^p \times \mathbb{R}$,
$A_i \colon {\cal X} \to \mathbb{R} \colon ({\mathbf w},b) \mapsto y_i({\mathbf x}_i^{\top}{\mathbf w}-b)$ \ $(i=1,2,\ldots, N)$,
$f\colon {\cal X} \to \mathbb{R}\colon ({\mathbf w}, b)  \mapsto h(A_N({\mathbf w},b))$, $g:=\bigoplus_{i=1}^{N-1} h$, and
$A \colon \mathbb{R}^{p+1} \to \mathbb{R}^{N-1}\colon ({\mathbf w},b) \mapsto \left(A_i ({\mathbf w},b)\right)_{i=1}^{N-1}$.
By using these translations, we can express the hinge loss function, in the form of the first stage cost function of \eqref{viscosity-conv-prob}, as
\begin{align*}
    \sum_{i=1}^{N}h\left(y_{i}({\mathbf w}^{\top}{\mathbf x}_{i }-b)\right)
    = g\circ A({\mathbf w},{b})+f({\mathbf w},{b})
\end{align*}
and its associated qualification condition (see \eqref{eq:qualification0}) is verified by
\rev{
  \begin{align}
    \label{eq:qualSVMred}
  \operatorname{ri}\left(\operatorname{dom}(g) - A\operatorname{dom}(f) \right)
  =\operatorname{ri}\left( \mathbb{R}^{N-1} - A\operatorname{dom}(f) \right)
  =\operatorname{ri}\left( \mathbb{R}^{N-1} \right)
  = \mathbb{R}^{N-1} \ni 0.
\end{align}
}
Note that, 
for any $\gamma \in \mathbb{R}_{++}$, the proximity operator $\operatorname{prox}_{\gamma h}$  can be computed as \eqref{eq:proxhinge} in Example \ref{ex:proximablefunctions}(c) (see Section \ref{sec:Monotone})
and therefore $\operatorname{prox}_{f}=\operatorname{prox}_{h \circ A_N}$ can also be computed by applying  \eqref{eq:proxcomp} and \eqref{eq:proxhinge} in Example \ref{ex:proximablefunctions}(b)(c).
  Moreover, by introducing $\Psi\colon \mathbb{R}^p \times \mathbb{R} \to \mathbb{R}\colon ({\mathbf w}, b) \mapsto \frac{1}{2}\| {\mathbf w}\|^2$, we 
  can regard Problem \eqref{SVM-nonseparable-opt3} as an instance of Problem \eqref{viscosity-conv-prob} under the assumption of ${\cal S}_p:= \Gamma \not = \varnothing$.\footnote{\label{fn:slightmodification}If we need to guarantee ${\cal S}_p\text{[in \eqref{viscosity-conv-prob}]} \not = \varnothing$, we recommend the following slight modification of \eqref{SVM-nonseparable-opt3}:
    \begin{align*}
      \hspace{-5mm}
      \underset{{\mathbf w}^{\star} \in \widetilde{\Gamma}}{\operatorname{minimize}} & \ \frac{1}{2}\|{\mathbf w}^{\star}\|^2 
      \text{ subject to }  \widetilde{\Gamma}:= \underset{({\mathbf w}, b) \in \mathbb{R}^p \times \mathbb{R}}{\operatorname{argmin}}\left[
      \Phi({\mathbf w}, b):=
      \iota_{\overline{B}(0,r)}({\mathbf w}, b)+ \sum_{i=1}^{N}h\left(y_{i}({\mathbf w}^{\top}{\mathbf x}_{i }-b)\right)
      \right]
    \end{align*}
    \revb{with a sufficiently large closed ball $\overline{B}(0,r)$},
    where ${\cal S}_p:=\widetilde{\Gamma} \not = \varnothing$ is guaranteed due to the coercivity of $\Phi$. Fortunately, our strategies in Section \ref{sec:3} are still applicable to this modified problem because it is also
    an instance of 
    \eqref{viscosity-conv-prob-1} which can be translated into \eqref{viscosity-conv-prob} as explained in Section \ref{sec:1}. \revb{In the application of Theorem \ref{th:DRSave} in Section~\ref{sec:3.1} to this modification, the boundedness of $\operatorname{Fix}({\mathbf T}_{\rm DRS_{II}})$ is automatically guaranteed because of \revb{Corollary \ref{pr:uibqehrwh4}(b)} (see Section~\ref{sec:3.3}) and the boundedness of both $\widetilde{\Gamma} \subset \overline{B}(0,r)$ and $\partial h(\mathbb{R})=\partial h(\mathbb{R}\setminus\{1\}) \cup \partial h(\{1\})
      =[-1,0]$.}
  }
  In fact, we can apply Theorem~\ref{th:DRS}, Theorem~\ref{th:DRSave}, and Theorem~\ref{th:HSDMLAL} to \eqref{SVM-nonseparable-opt3} because $\Psi$ is not strictly convex.
  In the following numerical experiment, we applied Theorem \ref{th:DRSave} to \eqref{SVM-nonseparable-opt3}
  with slight modification$^{16}$.

  \subsection{Numerical Experiment:  Margin Maximization with Least Empirical Hinge Loss}
\begin{figure}[t]
  \centering
  \footnotesize
  \centerline{\includegraphics[width=14cm]{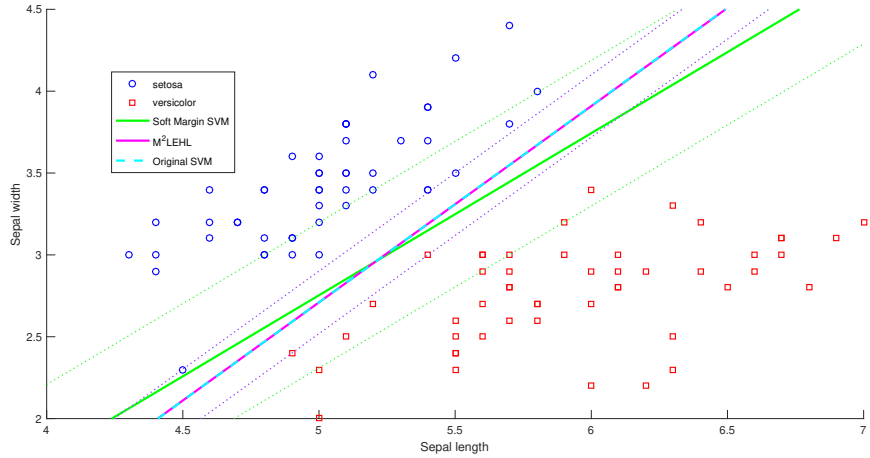}}
  \ \\[-4mm]
  \caption{Comparison between M$^2$LEHL, Original SVM, and Soft Margin SVM (Case of a separable training dataset ${\cal D}_{\rm sep}$): M$^2$LEHL reproduces Original SVM.}
  \label{fig:sep}
\end{figure}

\begin{figure}[t]
  \centering
  \footnotesize
  \centerline{\includegraphics[width=14cm]{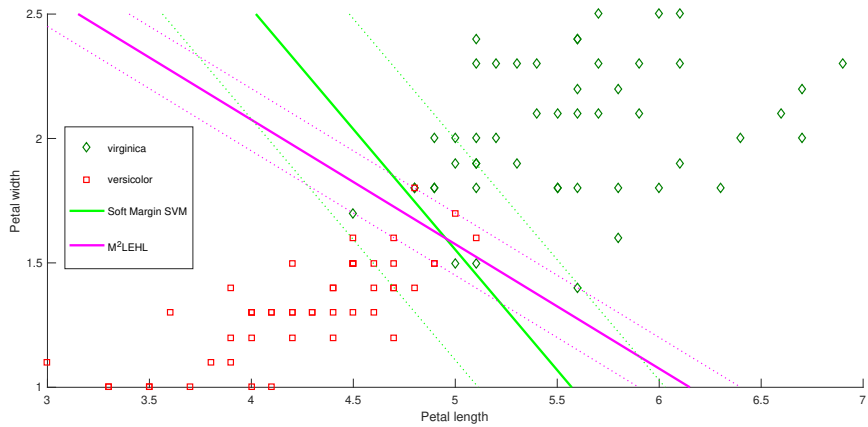}}
  \ \\[-4mm]
  \caption{Comparison between M$^2$LEHL and Soft Margin SVM (Case of a nonseparable training dataset ${\cal D}_{\rm nsep}$).
    }
  \label{fig:nosep}
\end{figure}

  We demonstrate that, as an extension of the original SVM in \eqref{def-orig-SVM}, the hierarchical enhancement of the SVM in \eqref{SVM-nonseparable-opt3} is more faithful to the original SVM than the soft margin SVM \eqref{SVM-nonseparable-opt1}.
  In our experiment, we applied the original SVM in \eqref{def-orig-SVM}, the soft margin SVM in \eqref{SVM-nonseparable-opt1},
  and the proposed hierarchical enhancement of the SVM in \eqref{SVM-nonseparable-opt3}
  to the Iris dataset which is famous dataset used firstly in Fisher's paper \cite{Fisher1936}. This data set has $150$ sample points, which are divided into three classes (I(setosa), II(versicolor), III(virginica)), and each sample point has four features (sepal length, sepal width, petal length, and petal width). 
  From Iris dataset, we construct two datasets: separable ${\cal D}_{\rm sep} \subset \mathbb{R}^2 \times \{-1,1\}$ with $|{\cal D}_{\rm sep}|=100$ comprising all the samples of Class I and Class II having only sepal length and sepal width; and non-separable ${\cal D}_{\rm nsep} \subset \mathbb{R}^2 \times \{-1,1\}$  with $|{\cal D}_{\rm nsep}|=100$ comprising  all the samples of Class II and Class III having only petal length and petal width.
  For each linear classifier ${\cal L}_{({\mathbf w},b)}$ of our interest,
  the three hyperplanes $\Pi_{({\mathbf w},b)}$, $\Pi_{({\mathbf w},b+1)}$, $\Pi_{({\mathbf w},b-1)}$ (see \eqref{halfspaces+class0}) are drawn in Figure \ref{fig:sep} and Figure \ref{fig:nosep}, in cyan for ``Original SVM'', in green for ``Soft Margin SVM'', and in magenta for the proposed ``M$^2$LEHL'' (which stands for \emph{the Margin Maximization with Least Empirical Hinge Loss}), respectively, where $({\mathbf w}, b)_{\rm org}$ is obtained by applying a quadratic programming solver \verb|quadprog| in Matlab to \eqref{def-orig-SVM}, $({\mathbf w}, b)_{\rm soft}$ is obtained by applying a soft margin SVM solver \verb|fitcsvm| (with the default setting, i.e., $\mathfrak{C}=1$) in Matlab to \eqref{SVM-nonseparable-opt1}, and $({\mathbf w}, b)_{\rm M^2LEHL}$ is obtained by applying the proposed algorithmic solution in Section \ref{subsec:4.2}, to \eqref{SVM-nonseparable-opt3}, designed based on Theorem \ref{th:DRSave} with slight modification$^{16}$.

Figure \ref{fig:sep} illustrates the resulting separating hyperplanes for the separable dataset ${\cal D}_{\rm sep}$.
Since the magenta lines are completely overlapped with cyan lines, ``M$^2$LEHL'' reproduces ``Original SVM'', as explained in just after \eqref{SVM-nonseparable-opt3}.
``Soft Margin SVM'' does not succeed in maximizing the margin, i.e., \eqref{SVM-nonseparable-opt3} is a more faithful extension of the original SVM in \eqref{def-orig-SVM} than the soft margin SVM \eqref{SVM-nonseparable-opt1}.

Figure \ref{fig:nosep} illustrates the resulting separating hyperplanes for the nonseparable dataset ${\cal D}_{\rm nsep}$. 
Since the original SVM \eqref{def-orig-SVM} has no solution, ``Original SVM'' is not depicted.
As the performance measure, we employ the number of errors $|{\cal E}(\cdot)|$ defined in \eqref{min-subset} along the original goal \eqref{min-subset1} (as suggested in  \cite[Sec.~3]{Cortes-Vapnik95}). 
 Though ``Soft Margin SVM'' has  $21$ errors, ``M$^2$LEHL'' achieves only $6$ errors, which demonstrates that \eqref{SVM-nonseparable-opt3} is more effective formulation for approaching to the original goal \eqref{min-subset1} than the soft margin SVM \eqref{SVM-nonseparable-opt1}.

\section{Application to Hierarchical Enhancement of Lasso}
\label{sec:5}  
\subsection{TREX : A Nonconvex Automatic Sparsity Control of Lasso}
\label{subsec:5.1}
Consider the estimation of a sparse vector ${\mathbf b}^{\rm tru}\in {\mathbb R}^{p}$ in the standard linear regression model:
\begin{equation}\label{linear-model}
{\mathbf z}={\mathbf X}{\mathbf b}^{\rm tru}+\sigma {\mathbf e},
\end{equation}
where ${\mathbf z}=(z_1,\ldots, z_N)^{\top}\in {\mathbb R}^{N}$ is a response vector, ${\mathbf X}\in {\mathbb R}^{N\times p}$ 
a design matrix, $\sigma>0$ a constant, ${\mathbf e}=(\varepsilon_1,\ldots, \varepsilon_N)^{\top}$ the noise vector, 
each $\varepsilon_i$ is the realization of a random variable with mean zero and variance 1.  

The Lasso (Least Absolute Shrinkage and Selection Operator) \cite{Ro-Tibshirani96} has been used widely as one of the most 
well-known sparsity aware statistical estimation methods \cite{HastieTibshiraniFreidman09,Hastie15}. The Lasso for (\ref{linear-model}) is defined as 
a  minimizer of the least squares criterion with $\ell_{1}$ penalty, i.e., 
\begin{equation}\olabel{lasso}
\mathbf b_{\rm Lasso}(\lambda)\in \underset{{{\mathbf b}\in {\mathbb R}^{p}}}{\operatorname{argmin}} \
\frac{1}{2N}\left\|{\mathbf z}-{\mathbf X}{\mathbf b}\right\|_{2}^{2}+\lambda \|{\mathbf b}\|_{1}, 
\end{equation}
where the tuning parameter $\lambda>0$ aims at controlling the sparsity of ${\mathbf b}_{\rm Lasso}(\lambda)$. 
However selection of $\lambda>0$ is highly influential to ${\mathbf b}_{\rm Lasso}(\lambda)$ and therefore its reliable way of selection has been strongly desired.
Among many efforts toward automatic sparsity control of Lasso, the following prediction bound offers 
a firm basis and has been applied widely in recent strategies including \cite{van-de-Geer13,habiri13,Dalalyan17,Lederer15}. 

\begin{fact} \label{pred-bound}(A Prediction Bound of Lasso \cite{Kolt11,Rigollet11})
For $\lambda\geq \frac{
                       2\left\|{\mathbf X}^{\top}\left({\mathbf z}-{\mathbf X}{\mathbf b}^{\rm tru}\right)\right\|_{\infty}
                     }
                     {N}$, it holds                      
                     $\frac{
                           \left\|              
                           {\mathbf X}{\mathbf b}_{\rm Lasso}(\lambda)-
                           {\mathbf X}{\mathbf b}^{\rm tru}
                           \right\|^{2}
                           }{N}
                      \leq 2\lambda \|{\mathbf b}^{\rm tru}\|_{1}.                      
                     $
\end{fact}

The TREX (Tuning-free Regression that adapts to the Entire noise $\sigma {\mathbf e}$ and the design matrix ${\mathbf X}$) \cite{Lederer15} is one of
the state-of-the-art strategies based on Fact \ref{pred-bound}. The TREX 
is defined as a solution of a nonconvex optimization problem:
\begin{equation}\label{trex}
\text{find } \ {\mathbf b}_{\rm TREX}\in 
\underset{{{\mathbf b}\in {\mathbb R}^{p}}}{\operatorname{argmin}} \
\frac{\left\|{\mathbf X}{\mathbf b}-{\mathbf z}\right\|^{2}}
{\left\|{\mathbf X}^{\top}\left({\mathbf X}{\mathbf b}-{\mathbf z}\right)\right\|_{\infty}}
+\beta \|{\mathbf b}\|_{1},
\end{equation}  
where $\left\|{\mathbf X}^{\top}\left({\mathbf X}{\mathbf b}-{\mathbf z}\right)\right\|_{\infty}=
\max_{1\leq j\leq p}\left|{\mathbf X}_{:j}^{\top}\left({\mathbf X}{\mathbf b}-{\mathbf z}\right)\right|$, ${\mathbf X}_{:j}$ denotes the $j$th column of ${\mathbf X}$, and the parameter $\beta$ can be set to a constant value ($\beta=1/2$ being the default choice).

The authors in \cite{Bien17} cleverly decomposed the nonconvex optimization (\ref{trex}) into $2p$ subproblems: 
\begin{equation}\label{trex-subp}
\text{find } \ {\mathbf b}^{(j)}_{\rm TREX}\in \underset{\substack{{\mathbf b}\in {\mathbb R}^{p}\\
{\mathbf x}_{j}^{\top}\left({\mathbf X}{\mathbf b}-{\mathbf z}\right)>0}}{\operatorname{argmin}}
\left[
\frac{\left\|{\mathbf X}{\mathbf b}-{\mathbf z}\right\|^{2}}
{
\beta{\mathbf x}_{j}^{\top}\left({\mathbf X}{\mathbf b}-{\mathbf z}\right)
}
+ \|{\mathbf b}\|_{1}
\right],
\end{equation}   
where 
\begin{align}
\label{eq:def-xj}
  {\mathbf x}_{j}=   \begin{cases}
    {\mathbf X}_{:j} & (j=1,2,\ldots, p); \\
    -{\mathbf X}_{:j-p} & (j=p+1,p+2,\ldots, 2p).
\end{cases}
\end{align}
More precisely, ${\mathbf b}_{\rm TREX}$ in (\ref{trex}) is characterized as 
\begin{align}\label{decomp-TREX}
{\mathbf b}_{\rm TREX}\in  \widehat{\cal Q}_{\mathbb{R}^p}\left[\underset{\substack{({\mathbf b},j)\in {\mathbb R}^{p} \times \{1,2,\ldots, 2p\}\\
{\mathbf x}_{j}^{\top}\left({\mathbf X}{\mathbf b}-{\mathbf z}\right)>0}}{\operatorname{argmin}}
\left(
\frac{\left\|{\mathbf X}{\mathbf b}-{\mathbf z}\right\|^{2}}
{
\beta{\mathbf x}_{j}^{\top}\left({\mathbf X}{\mathbf b}-{\mathbf z}\right)
}
+ \|{\mathbf b}\|_{1}
\right)\right], 
\end{align}
where 
\begin{align}
\label{eq:QTREX}
  \widehat{\cal Q}_{\mathbb{R}^p}\colon {\mathbb R}^{p}\times \{1,2,\ldots, 2p\} \to \mathbb{R}^p\colon ({\mathbf b},j) \mapsto {\mathbf b}.
\end{align}

Remarkably, each subproblem (\ref{trex-subp}) was shown to be a convex optimization and solved in \cite{Bien17} with a second-order cone program (SOCP) \cite{lebo98}. 

\begin{figure}[t]
  \centering
  \footnotesize
  \centerline{\includegraphics[width=11cm]{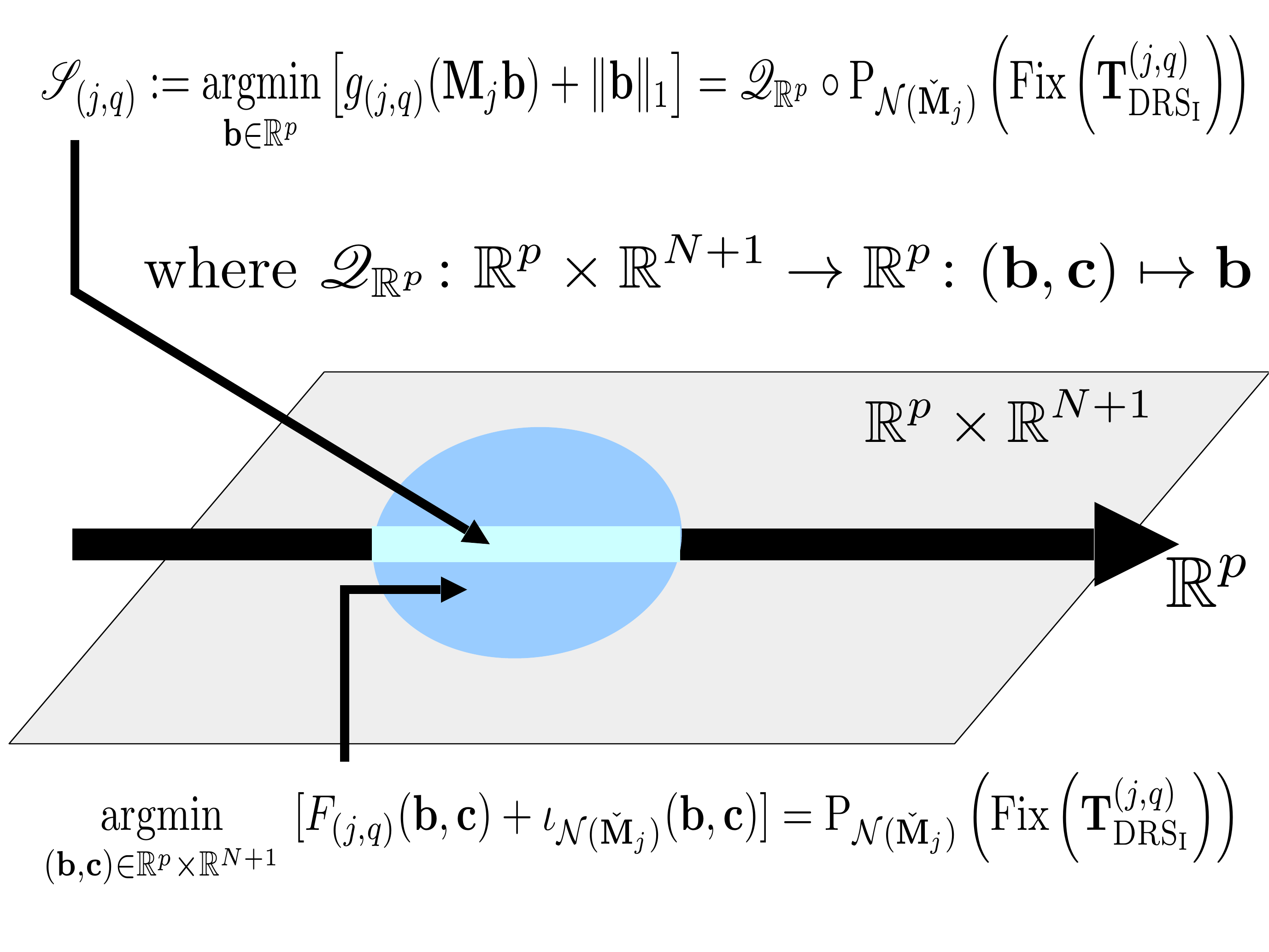}}
  \ \\[-4mm]
  \caption{Illustration of the fixed point characterization of ${\cal S}_{(j,q)}$ in \eqref{p-gtrex-subp} via the Douglas-Rachford splitting operator ${\mathbf T}_{\rm DRS_{I}}^{(j,q)}$ in \eqref{eq:DRS31svm}. }
  \label{fig:TREX}
\end{figure}

Recently, for sound extensions of the subproblem (\ref{trex-subp}) as well as for sound applications of proximal splitting, 
a successful reformulation of \eqref{decomp-TREX} was made for general $q>1$ in \cite{plc-muller18} as 
\begin{align}\label{decomp-TREXref}
{\mathbf b}_{\rm TREX_q}\in {\cal S}_{\rm TREX_q}:=\widehat{\cal Q}_{\mathbb{R}^p}\left[
\underset{({\mathbf b},j)\in  {\mathbb R}^{p} \times \{1,2,\ldots, 2p\}}{\operatorname{argmin}}
g_{(j,q)}({\mathbf M}_{j}{\mathbf b})
+ \|{\mathbf b}\|_{1}
\right]
\end{align}
whose solution ${\mathbf b}_{\rm TREX_q}$ is given, by passing through $2p$ convex subproblems, as ${\mathbf b}^{(j^{\star})}_{\rm TREX_q}$, where
\begin{equation}\label{p-gtrex-subp}
  \hspace{-6mm}\left[
  \begin{array}{l}
{\mathbf b}^{(j)}_{\rm TREX_q}\in {\mathcal S}_{(j,q)}:=
    \underset{{\mathbf b}\in {\mathbb R}^{p}}{\operatorname{argmin}}
\left[
g_{(j,q)}({\mathbf M}_{j}{\mathbf b})
+ \|{\mathbf b}\|_{1}
\right]  \ (j=1,2,\ldots, 2p); \\
    j^{\star} \in    
\underset{j \in \{1,2,\ldots, 2p\}}{\operatorname{argmin}}
\left[
g_{(j,q)}({\mathbf M}_{j}{\mathbf b}^{(j)}_{\rm TREX_q})
+ \|{\mathbf b}^{(j)}_{\rm TREX_q}\|_{1}
    \right],
  \end{array}
  \right.
\end{equation}
\vspace{-2mm}
\begin{equation}\label{def-gqj}
g_{(j,q)}:{\mathbb R}\times {\mathbb R}^{N}\rightarrow (-\infty, \infty]: (\eta,{\mathbf y})\mapsto 
\left\{
\begin{array}{ll}
\frac{\|{\mathbf y}-{\mathbf z}\|^{q}}
{\beta (\eta-{\mathbf x}_{j}^{\top}{\mathbf z})^{q-1}},
&\mbox{if }\eta>{\mathbf x}_{j}^{\top}{\mathbf z};\\
0,&\mbox{if }{\mathbf y}={\mathbf z}\mbox{ and }\eta={\mathbf x}_{j}^{\top}
{\mathbf z};\\
+\infty,&\mbox{otherwise} 
\end{array}
\right. 
\end{equation}
is a proper lower semicontinuous convex function, and 
\begin{equation}\label{def-Mj}
{\mathbf M}_{j}:{\mathbb R}^{p}\rightarrow {\mathbb R}\times {\mathbb R}^{N}: 
{\mathbf b}\mapsto 
\left({\mathbf x}_{j}^{\top}{\mathbf X}{\mathbf b}, {\mathbf X}{\mathbf b}\right) 
\end{equation}
is a bounded linear operator. The estimator ${\mathbf b}_{\rm TREX_q}$ in \eqref{decomp-TREXref} is called the generalized TREX in \cite{plc-muller18} where, as its specialization, ${\mathbf b}_{\rm TREX_2}$ is also called TREX. \revb{
  Note that, in view of Example \ref{ex:proximablefunctions}(d)(e) in Section~\ref{sec:Monotone} and a relation between $g_{(j,q)}$ and $\widetilde{\varphi}_{q}$ in \eqref{def-perspective-phi2} \cite[in Section 4.3.2]{plc-muller18}:
  \begin{equation}\label{rel-perspective}
    (\forall (\eta,{\mathbf y}) \in \mathbb{R} \times \mathbb{R}^N ) \quad g_{(j,q)}(\eta,{\mathbf y})=
    \mbox{}_{\tau_{\left({\mathbf x}_{j}^{\top}{\mathbf z},{\mathbf z}\right)}}  \widetilde{\varphi}_{q}\left(\eta, {\mathbf y}\right)=
    \widetilde{\varphi}_{q}\left(\eta-{\mathbf x}_{j}^{\top}{\mathbf z}, {\mathbf y}-{\mathbf z}\right),
  \end{equation}
  each convex subproblem in \eqref{p-gtrex-subp} is an instance of Problem \eqref{typical-conv-opt}.
}
  For the subproblem \eqref{p-gtrex-subp}, the Douglas-Rachford splitting method (see Proposition \ref{pr:DRS} in Section~\ref{subsec:2.2}) was successfully applied in \cite{plc-muller18}. For completeness, we reproduce this result in the style of (\ref{generic-view1}--\ref{def-nexp})
   followed by application of Fact \ref{fa:Mann} (in Section \ref{sec:Monotone}) to the characterization \eqref{eq:DRScharacterization0}. Suppose that for \eqref{p-gtrex-subp} the qualification condition (see \eqref{eq:qualification0})
\begin{align}
  \label{eq:qualcondTREX}
     \rev{0 \in \operatorname{ri}( \operatorname{dom}(g_{(j,q)}) - {\mathbf M}_j\operatorname{dom}(\|\cdot\|_1))}
\end{align}
holds\footnote{
  In \cite{plc-muller18}, the qualification condition \eqref{eq:qualcondTREX} seems to be assumed implicitly. If we assume additionally that ${\mathbf X}\in \mathbb{R}^{N \times p}$ has no zero column, it is automatically guaranteed as will be shown in Lemma \ref{pr:TREXqualcond} in Section \ref{subsec:5.2}.
   }.
Then, by using 
\begin{align}
  \left[
  \begin{array}{l}
  \check{\mathbf M}_j\colon  
  \mathbb{R}^p \times \mathbb{R}^{N+1} \to \mathbb{R}^{N+1}\colon ({\mathbf b},{\mathbf c}) \mapsto {\mathbf M}_j{\mathbf b} - {\mathbf c},
                               \\
    \mathcal{Q}_{\mathbb{R}^p}\colon  \mathbb{R}^p \times \mathbb{R}^{N+1} \to  \mathbb{R}^p\colon ({\mathbf b},{\mathbf c})  \mapsto {\mathbf b},
  \end{array}
                                       \right.
                                       \olabel{eq:defcMj}
\end{align}
we obtain
\begin{align}\label{eq:theorem31bTREX}
  \mathcal{S}_{(j,q)}= \mathcal{Q}_{\mathbb{R}^p} \circ  {P}_{\mathcal{N}(\check{\mathbf M}_j)}\left(\operatorname{Fix}\left({\mathbf T}_{\rm DRS_{I}}^{(j,q)}\right) \right)
\end{align}
(which is a specialization of \eqref{eq:DRScharacterization0} for \eqref{p-gtrex-subp}, see Figure \ref{fig:TREX}), where 
  \begin{align}
    \label{eq:DRS31svm}
    {\mathbf T}_{\rm DRS_{I}}^{(j,q)}\colon \mathbb{R}^p \times \mathbb{R}^{N+1} \to  \mathbb{R}^p \times \mathbb{R}^{N+1}\colon ({\mathbf b},{\mathbf c}) \mapsto ({\mathbf b}_T,{\mathbf c}_T)
  \end{align} 
is the DRS operator of Type-I (c.f., \eqref{eq:DRSIo} and \eqref{eq:DRSIimplicit}) specialized for \eqref{p-gtrex-subp} and is defined by 
  \begin{align}
    \olabel{eq:DRSTREXupdate}
    \left[
    \begin{array}{l}
      {\mathbf p}=  {\mathbf b}-{\mathbf M}_j^*( {\mathbf I} + {\mathbf M}_j{\mathbf M}_j^*)^{-1}  ({\mathbf M}_j{\mathbf b} -{\mathbf c}) \\
      ({\mathbf b}_{1/2},{\mathbf c}_{1/2}) =  (2{\mathbf p}-{\mathbf b}, 2{\mathbf M}_j{\mathbf p}-{\mathbf c})) \\
    ({\mathbf b}_{T},{\mathbf c}_{T}) =  (2\operatorname{prox}_{\|\cdot\|_1}({\mathbf b}_{1/2}) -{\mathbf b}_{1/2}, 2\operatorname{prox}_{g_{(j,q)}}({\mathbf c}_{1/2}) -{\mathbf c}_{1/2}),
      \end{array}
    \right.
  \end{align}
  or equivalently by
  \begin{align}
    \olabel{eq:DRSTREXexpression}
    {\mathbf T}_{\rm DRS_{I}}^{(j,q)}:= (2 \operatorname{prox}_{F_{(j,q)}} -{\mathbf I}) \circ(2 {P}_{{\cal N}(\check{\mathbf M}_j)} - {\mathbf I})
  \end{align}
   with $F_{(j,q)}\colon\mathbb{R}^p \times \mathbb{R}^{N+1} \to (-\infty,\infty]\colon ({\mathbf b}, {\mathbf c}) \mapsto g_{(j,q)}({\mathbf c})+\|{\mathbf b}\|_1$. Note that ${\mathbf T}_{\rm DRS_{I}}^{(j,q)}$ can be computed efficiently if $( {\mathbf I} + {\mathbf M}_j{\mathbf M}_j^*)^{-1}$ is available as a computational tool. 

 The above characterization \eqref{eq:theorem31bTREX} and Fact \ref{fa:Mann} (see Section \ref{sec:Monotone}) lead to the following algorithmic solution of \eqref{p-gtrex-subp}.
\begin{fact}[Douglas-Rachford Splitting Method for Subproblems of Generalized TREX]
\label{fa:DRSTREX}
Under the qualification condition \eqref{eq:qualcondTREX} for ${\cal S}_{(j,q)}$ in \eqref{p-gtrex-subp}, the sequence $({\mathbf b}_{n}, {\mathbf c}_{n})_{n \in \mathbb{N}} \subset \mathbb{R}^p \times \mathbb{R}^{N+1} $ generated, with $(\alpha_n)_{n \in \mathbb{N}} \subset [0,1]$ \revb{satisfying $\sum_{n\in \mathbb{N}} \alpha_n(1-\alpha_n) = \infty$} in Fact~\ref{fa:Mann} (see Section \ref{sec:Monotone}) and $({\mathbf b}_{0}, {\mathbf c}_{0}) \in \mathbb{R}^p \times \mathbb{R}^{N+1}$, by 
\begin{align}
  \label{eq:DRSgTREX}
  ({\mathbf b}_{n+1}, {\mathbf c}_{n+1}) = (1-\alpha_n)  ({\mathbf b}_{n}, {\mathbf c}_{n})+ \alpha_n {\mathbf T}_{\rm DRS_{I}}^{(j,q)}({\mathbf b}_{n}, {\mathbf c}_{n})
\end{align}
converges to a point $({\mathbf b}_{\star}, {\mathbf c}_{\star})$ in $\operatorname{Fix}\left({\mathbf T}_{\rm DRS_{I}}^{(j,q)}\right)$ as well as the sequence 
$(\mathcal{Q}_{\mathbb{R}^p} \circ  {P}_{\mathcal{N}(\check{\mathbf M}_j)} ({\mathbf b}_{n}, {\mathbf c}_{n}))_{n \in \mathbb{N}}$ converges to $\mathcal{Q}_{\mathbb{R}^p} \circ  {P}_{\mathcal{N}(\check{\mathbf M}_j)} ({\mathbf b}_{\star}, {\mathbf c}_{\star}) \in \mathcal{S}_{(j,q)}$, where
\begin{align*}
  \mathcal{Q}_{\mathbb{R}^p} \circ  {P}_{\mathcal{N}(\check{\mathbf M}_j)}\colon \mathbb{R}^p \times \mathbb{R}^{N+1} \to  \mathbb{R}^p\colon ({\mathbf b}, {\mathbf c}) \mapsto {\mathbf b}-{\mathbf M}_j^*( {\mathbf I} + {\mathbf M}_j{\mathbf M}_j^*)^{-1}  ({\mathbf M}_j{\mathbf b} -{\mathbf c}).
\end{align*}
\end{fact}
Note that the sequence $(\mathcal{Q}_{\mathbb{R}^p} \circ  {P}_{\mathcal{N}(\check{\mathbf M}_j)} ({\mathbf b}_{n}, {\mathbf c}_{n}))_{n \in \mathbb{N}}$ can be generated efficiently by \eqref{eq:DRSgTREX} if $( {\mathbf I} + {\mathbf M}_j{\mathbf M}_j^*)^{-1}$ is available as a computational tool. %

\subsection{Enhancement of Generalized TREX Solutions with Hierarchical Optimization}
\label{subsec:5.2}

Along Scenario 2 in Section \ref{sec:1}, suppose that we found newly an effective criterion $\Psi \in \Gamma_0(\mathbb{R}^p)$ whose gradient is Lipschitzian \revb{over $\mathbb{R}^p$} and 
we hope to select a most desirable vector, in the sense of $\Psi$, from the solution set ${\cal S}_{\rm TREX_q}$ in \eqref{decomp-TREXref}.
This task is formulated as a \emph{hierarchical nonconvex optimization problem} (see \eqref{eq:QTREX}, \eqref{def-gqj}, and \eqref{def-Mj} for $\widehat{\cal Q}_{\mathbb{R}^p}$,  $g_{(j,q)}$, and ${\mathbf M}_j$):
\begin{align}
  \label{eq:HNGTREX}
\hspace{-5mm}   \text{minimize } & \ \Psi({\mathbf b}^{\star}) \\
  \text{ subject to } & \                       
{\mathbf b}^{\star} \in {\cal S}_{\rm TREX_q}=
 \widehat{\cal Q}_{\mathbb{R}^p}
\left[
                      \underset{({\mathbf b},j) \in \mathbb{R}^p \times \{1,2,\ldots, 2p\}}{\operatorname{argmin}}\ g_{(j,q)}({\mathbf M}_j{\mathbf b})+\|{\mathbf b}\|_1\right]
                      \nonumber
\end{align}
whose solution ${\mathbf b}_{\rm HTREX_q}$ is given, by passing through $2p$ (\emph{hierarchical convex} optimization) subproblems, as ${\mathbf b}^{(j^{\star\star})}_{\rm HTREX_q}$, where
\begin{equation}\label{eq:HHE}
  \left[
  \begin{array}{l}
    {\mathbf b}^{(j)}_{\rm HTREX_q}\in \Omega_{\rm DRS_{I}}^{(j,q)}:=\underset{{\mathbf b}^{\star} \in {\cal S}_{(j,q)}}{\operatorname{argmin}}
\Psi({\mathbf b}^{\star}), \\
    {\mathcal S}_{(j,q)}=
\underset{{\mathbf b}\in {\mathbb R}^{p}}{\operatorname{argmin}}
\left[
g_{(j,q)}({\mathbf M}_{j}{\mathbf b})
+ \|{\mathbf b}\|_{1}
\right]  \quad (j=1,2,\ldots, 2p) \  \text{[in \eqref{p-gtrex-subp}]}; \\
    \mathfrak{J}^{\star} :=
\underset{j \in \{1,2,\ldots, 2p\}}{\operatorname{argmin}}
\left[
g_{(j,q)}({\mathbf M}_{j}{\mathbf b}^{(j)}_{\rm HTREX_q})
+ \|{\mathbf b}^{(j)}_{\rm HTREX_q}\|_{1}
    \right]; \\
    j^{\star\star} \in
    \underset{j^{\star} \in  \mathfrak{J}^{\star}}{\operatorname{argmin}} \
    \Psi({\mathbf b}^{(j^{\star})}_{\rm HTREX_q}).
  \end{array}
  \right.
\end{equation}

\noindent \revb{Note that the coercivity of $\|\cdot \|_1$ and the nonnegativity of $g_{(j,q)}$ ensure that ${\cal S}_{(j,q)}$ is nonempty and bounded (see Fact \ref{fa:uobqegbuqeh}(c) in Section \ref{subsec:2.1}), which also guarantees $\Omega_{\rm DRS_{I}}^{(j,q)}=\operatorname{argmin}(\iota_{{\cal S}_{(j,q)}}+\Psi)(\mathbb{R}^p) \not = \varnothing$ $(j=1,2,\ldots, 2p)$ by the classical Weierstrass theorem. }

In the following, we focus on how to compute the solution ${\mathbf b}^{(j)}_{\rm HTREX_q}$ ($j=1,2,\ldots, 2p$) %
in \eqref{eq:HHE} by a proposed strategy in Section \ref{sec:3}. We assume that the design matrix ${\mathbf X}\in \mathbb{R}^{N \times p}$ in \eqref{linear-model} has no zero column, to guarantee the qualification condition \eqref{eq:qualcondTREX} for ${\cal S}_{(j,q)}$ in \eqref{p-gtrex-subp} for each $j=1,2,\ldots, 2p$.
\begin{lem}
  \label{pr:TREXqualcond}
  Suppose that the design matrix ${\mathbf X}\in \mathbb{R}^{N \times p}$ has no zero column. Then
  the qualification condition \eqref{eq:qualcondTREX} for ${\cal S}_{(j,q)}$ in \eqref{p-gtrex-subp} is guaranteed automatically for each $j=1,2,\ldots, 2p$.
\end{lem}
(The proof of Lemma \ref{pr:TREXqualcond} is given in Appendix G).
\begin{theorem}[Algorithmic Solution to Hierarchical TREX$_q$] %
  \label{th:DRSTREX}
  \revb{Suppose that ${\mathbf X}$ has no zero column and $\Psi \in \Gamma_0(\mathbb{R}^p)$ is \revb{G\^{a}teaux} differentiable with Lipschitzian gradient $\nabla \Psi $ \revb{over $\mathbb{R}^p$}.}
  Then, for ${\mathbf T}_{\rm DRS_{I}}^{(j,q)}$ in \eqref{eq:DRS31svm} ($j=1,2,\ldots, 2p$),
  \begin{enumerate}
  \item[\rm (a)] $\operatorname{Fix}({\mathbf T}_{\rm DRS_{I}}^{(j,q)})$ is bounded; \\
  \item[\rm (b)]
    ${\mathbf T}_{\rm DRS_{I}}^{(j,q)}$ can be plugged into the HSDM \eqref{eq:updateHSDM}, with any $\alpha \in (0,1)$ and $(\lambda_{n+1})_{n \in \mathbb{N}} \in \ell_{+}^2 \setminus \ell_{+}^1$, as
  \begin{align}
    \label{eq:HSDM311TREX} 
    \hspace{-5mm} 
    \left[
\begin{array}{l}
({\mathbf b}_{n+1/2},{\mathbf c}_{n+1/2}) =  \ (1-\alpha)({\mathbf b}_{n}, {\mathbf c}_{n})+ \alpha {\mathbf T}_{\rm DRS_{I}}^{(j,q)}({\mathbf b}_{n}, {\mathbf c}_{n}) \\
    {\mathbf b}_{n+1}^{\star}=  {\mathbf b}_{n+1/2}-{\mathbf M}_j^*({\mathbf I} + {\mathbf M}_j{\mathbf M}_j^*)^{-1}  ({\mathbf M}_j{\mathbf b}_{n+1/2} -{\mathbf c}_{n+1/2})  \\
  {\mathbf b}_{n+1} =  {\mathbf b}_{n+1/2} - \lambda_{n+1} ({\mathbf I}- {\mathbf M}_j^*({\mathbf I}+{\mathbf M}_j{\mathbf M}_j^*)^{-1}{\mathbf M}_j)\circ \nabla \Psi ({\mathbf b}_{n+1}^\star) \\
        {\mathbf c}_{n+1} =  {\mathbf c}_{n+1/2} - \lambda_{n+1} (({\mathbf I}+{\mathbf M}_j{\mathbf M}_j^*)^{-1}{\mathbf M}_j)\circ \nabla \Psi ({\mathbf b}_{n+1}^\star).
\end{array}
\right.
  \end{align}
  The algorithm \eqref{eq:HSDM311TREX} generates, for any $({\mathbf b}_{0}, {\mathbf c}_{0})  \in \mathbb{R}^p \times \mathbb{R}^{N+1}$, a sequence $({\mathbf b}_{n+1}^{\star})_{n \in \mathbb{N}} \subset \mathbb{R}^p$ which satisfies 
\begin{align*}
  \lim_{n \to \infty} d_{\Omega_{\rm DRS_{I}}^{(j,q)}}({\mathbf b}_{n}^{\star})=0,
\end{align*}
where $\Omega_{\rm DRS_{I}}^{(j,q)} \not = \varnothing$ is defined in \eqref{eq:HHE}.
\end{enumerate}
\end{theorem}
\begin{bold-remark}[Idea Behind Derivation of Theorem \ref{th:DRSTREX}]
  \label{re:HSDMTREX}
  \
  \begin{enumerate}
  \item[\rm (a)] 
    Recall that  ${\mathbf T}_{\rm DRS_{I}}^{(j,q)}$ is a DRS operator of Type-I (see \eqref{eq:DRS31svm} and Theorem \ref{th:DRS} in Section~\ref{sec:3.1}).
    By applying Corollary \ref{pr:uibqehrwh4}\revb{(a)} in Section \ref{sec:3.3} to Lemma \ref{pr:TREXqualcond}, the boundedness of ${\cal S}_{(j,q)} \not = \varnothing$, and the boundedness of the image of 
$\partial \|\cdot \|_1\colon \mathbb{R}^p \to [-1,1]^p \colon
{\mathbf b}=(b_1,b_2,\ldots, b_p) \mapsto \varprod_{i=1}^p\partial |~\cdot~|(b_i)$, we deduce the relation:
    \begin{align}
      & \text{$\left[-({\mathbf M}_j^{\top})^{-1}(\partial \|\cdot\|_1 ({\cal S}_{(j,q)}))\right] \cap \partial g_{(j,q)}(\mathbb{R}^{N+1}) $ is bounded} \label{eq:claimbounded}
      \\
      & \Rightarrow \operatorname{Fix}({\mathbf T}_{\rm DRS_{I}}^{(j,q)}) \text{ is bounded}, \nonumber 
    \end{align}
    where %
    $({\mathbf M}_j^{\top})^{-1}\colon \mathbb{R}^{p} \to 2^{\mathbb{R}^{N+1}}\colon {\mathbf b} \mapsto \{{\mathbf c} \in \mathbb{R}^{N+1} \mid  {\mathbf b} = {\mathbf M}_j^{\top}{\mathbf c}\}$ \revc{(see \eqref{def-Mj} for ${\mathbf M}_{j}$)}. Now, by  
$\partial g_{(j,q)}(\mathbb{R}^{N+1}) = \partial \widetilde{\varphi}_q(\mathbb{R}^{N+1})$ (due to \eqref{rel-perspective})
and the supercoercivity of ${\varphi}_{q}$ and ${\varphi}_{q}^*$ (due to \cite[Example 13.2 and Example 13.8]{bauschke11}),
for proving the boundedness of $\operatorname{Fix}({\mathbf T}_{\rm DRS_{I}}^{(j,q)})$ from \eqref{eq:claimbounded},
it is sufficient to show the following claim: \\[2mm]
{\bf Claim~\ref{th:DRSTREX}:} {\rm
  Suppose that ${\mathbf X}$ has no zero column.
  Let $S \subset \mathbb{R}^{p}$ be bounded, and
  $\varphi \in \Gamma_0(\mathbb{R}^{N})$ a supercoercive function having supercoercive $\varphi^* \in \Gamma_0(\mathbb{R}^{N})$.
 Then $({\mathbf M}_j^{\top})^{-1}(S) \cap \partial \widetilde{\varphi}(\mathbb{R}^{N+1})$ is bounded.} \\[2mm]
Note that Claim~\ref{th:DRSTREX} is the main step in the proof of Theorem \ref{th:DRSTREX}.
\item[\rm (b)] \revb{We have already confirmed the qualification condition \eqref{eq:qualcondTREX} in Lemma \ref{pr:TREXqualcond}, ${\cal S}_{(j,q)} \not = \varnothing$, and $\Omega_{\rm DRS_{I}}^{(j,q)}\not = \varnothing$ ($j=1,2,\ldots, 2p$) (see the short remark just after \eqref{eq:HHE}).}
  Therefore,
  application of Theorem \ref{th:DRS} (in Section \ref{sec:3.1}) to the subproblems to compute ${\mathbf b}^{(j)}_{\rm HTREX_q}$ ($j=1,2,\ldots, 2p$) in \eqref{eq:HHE} guarantees the statement of Theorem \ref{th:DRSTREX}(b).
  \end{enumerate}
\end{bold-remark}
(The proof of Theorem \ref{th:DRSTREX} is given in Appendix H).

\subsection{Numerical Experiment: Hierarchical TREX$_2$}
\label{subsec:numTREX}

We demonstrate that the proposed estimator ${\mathbf b}_{\rm HTREX_2}$, i.e., Hierarchical TREX$_2$ in \eqref{eq:HNGTREX} (see Section \ref{subsec:5.2}) can enhance further the estimation accuracy achieved by ${\mathbf b}_{\rm TREX_2}$ in \eqref{decomp-TREXref} if we can exploit another new criterion $\Psi\colon \mathbb{R}^p \to \mathbb{R}$ for promoting characteristics, of ${\mathbf b}^{\rm tru}$, which is not utilized in TREX$_2$.
Consider the situation where we like to estimate unknown vector
\begin{align}
  \olabel{eq:btru}
  {\mathbf b}^{\rm tru}=\frac{1}{\sqrt{p}}(0,0,0,1,1,1,0,\ldots,0)^{\top} \in \mathbb{R}^p
\end{align}
from the noisy observation ${\mathbf z} \in \mathbb{R}^N$ in \eqref{linear-model}. We suppose to know that ${\mathbf b}^{\rm tru}$ is not only sparse but also \emph{fairly flat}. Here, the fairly flatness of ${\mathbf b}^{\rm tru}$ means that the energy of oscillations (i.e., the sum of the squared gaps between the adjacent components) of ${\mathbf b}^{\rm tru}$ is small, which
is supposed to be our additional knowledge not utilized in the TREX$_2$ and Lasso estimators. If we have such prior knowledge, suppression of 
\begin{align}
  \label{eq:fdf}
  \Psi\colon \mathbb{R}^p \to \mathbb{R}\colon {\mathbf b} \mapsto \frac{1}{2}\|{\mathbf D} {\mathbf b} \|^2, \ \text{  with }
  {\mathbf D}:=
  \begin{pmatrix}
    -1 & 1 & 0 & 0 & 0 & 0 \\
    0 & -1 & 1 & 0 & 0 & 0 \\
     &  & \ddots & \ddots &  &  \\
    0 & 0 & 0 & -1 & 1 & 0 \\
    0 & 0 & 0 & 0 & -1 & 1
  \end{pmatrix} \in \mathbb{R}^{(p-1) \times p},
\end{align}
is expected to be effective for estimation of ${\mathbf b}^{\rm tru}$ because $\Psi$ can distinguish ${\mathbf b}^{\rm tru}$ from $\widetilde{\mathbf b}:=\frac{1}{\sqrt{p}}(0,1,0,0,1,1,0,\ldots,0)^{\top}\in \mathbb{R}^p$ of the same sparsity as ${\mathbf b}^{\rm tru}$ (i.e., $\|{\mathbf b}^{\rm tru}\|_0=\|\widetilde{\mathbf b}\|_0$ and $\|{\mathbf b}^{\rm tru}\|_1=\|\widetilde{\mathbf b}\|_1$)
by $\Psi({\mathbf b}^{\rm tru})<\Psi(\widetilde{\mathbf b})$.
Now, our new goal for enhancement of TREX$_2$ is to minimize $\Psi$ while keeping the optimality of the TREX$_2$ in the sense of \eqref{decomp-TREXref} for $q=2$ (see Scenario \ref{sce2} in Section \ref{sec:1}). This goal is achieved by solving the hierarchical nonconvex optimization \eqref{eq:HNGTREX} for $q=2$.

In our experiments, the design matrix ${\mathbf X} \in \mathbb{R}^{N\times p}$ in \eqref{linear-model} is given to satisfy ${\mathbf X}_{:2} = {\mathbf X}_{:3} = {\mathbf X}_{:4}$
with a sample of zero-mean Gaussian random variable
 followed by normalization to satisfy $\|{\mathbf X}_{:j}\| = \sqrt{N}$ ($j=1,\ldots, p)$.
The additive noise ${\mathbf e} \in \mathbb{R}^N$ in \eqref{linear-model} is drawn from the unit white Gaussian distribution. 
We tested the performances of the estimators under $({\rm SNR})=10 \log\left(\frac{\|{\mathbf X}{\mathbf b}^{\rm tru}\|^2}{\| \sigma {\mathbf e}\|^2}\right) \in [10,1000] \cup \{+\infty\}$, where $\sigma \in \mathbb{R}$ is adjusted to obtain a specific SNR.
Note that, in this setting, 
$\left\{ {\mathbf b} \in \mathbb{R}^p \mid \mathbf{Xb}=\mathbf{Xb}^{\rm tru} \text{ and } \|\mathbf{b} \|_1=\|\mathbf{b}^{\rm tru} \|_1\right\}$ is apparently an infinite set containing both ${\mathbf b}^{\rm tru}$ and $\widetilde{\mathbf b}$.

\begin{figure}[t]
  \centering
  \footnotesize
  \centerline{\includegraphics[width=10cm]{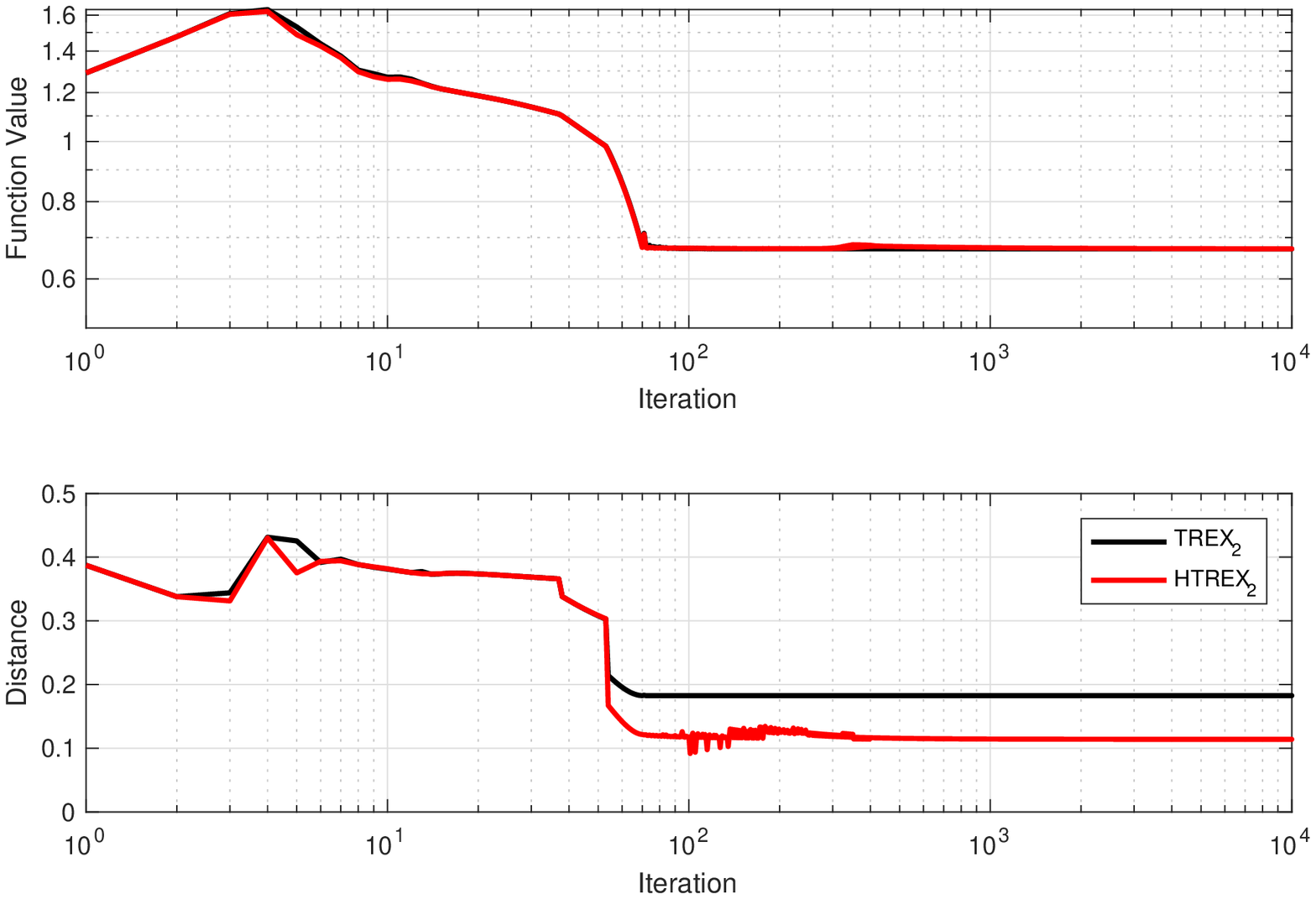}}
  \vspace{-8mm}
  (a) Comparison of TREX$_2$ and HTREX$_2$ in the process of convergences under the noise $\mathbf{e}=0$.
  \centerline{\includegraphics[width=10cm]{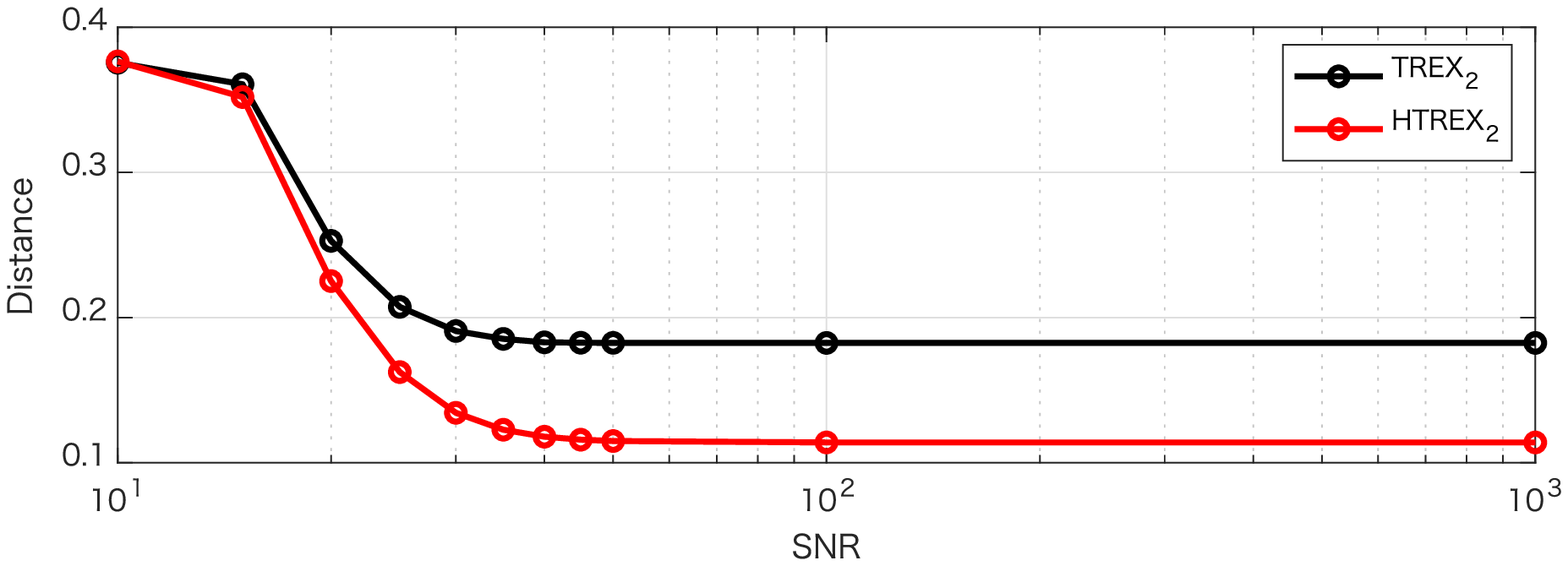}}
  (b)  Estimation accuracy achieved by TREX$_2$ and HTREX$_2$ for various SNR.
  \caption{Transient performances in an over-determined case; Criteria ``Function Value'', ``Distance'' are given in \eqref{eq:criteriaTREX} and $({\rm SNR})=10 \log\left(\frac{\|{\mathbf X}{\mathbf b}^{\rm tru}\|^2}{\| \sigma {\mathbf e}\|^2}\right) \ [{\rm dB}]$.}
  \label{fig:noisefree}
\end{figure}

\begin{figure}[t]
  \centering
  \footnotesize
  \centerline{\includegraphics[width=10cm]{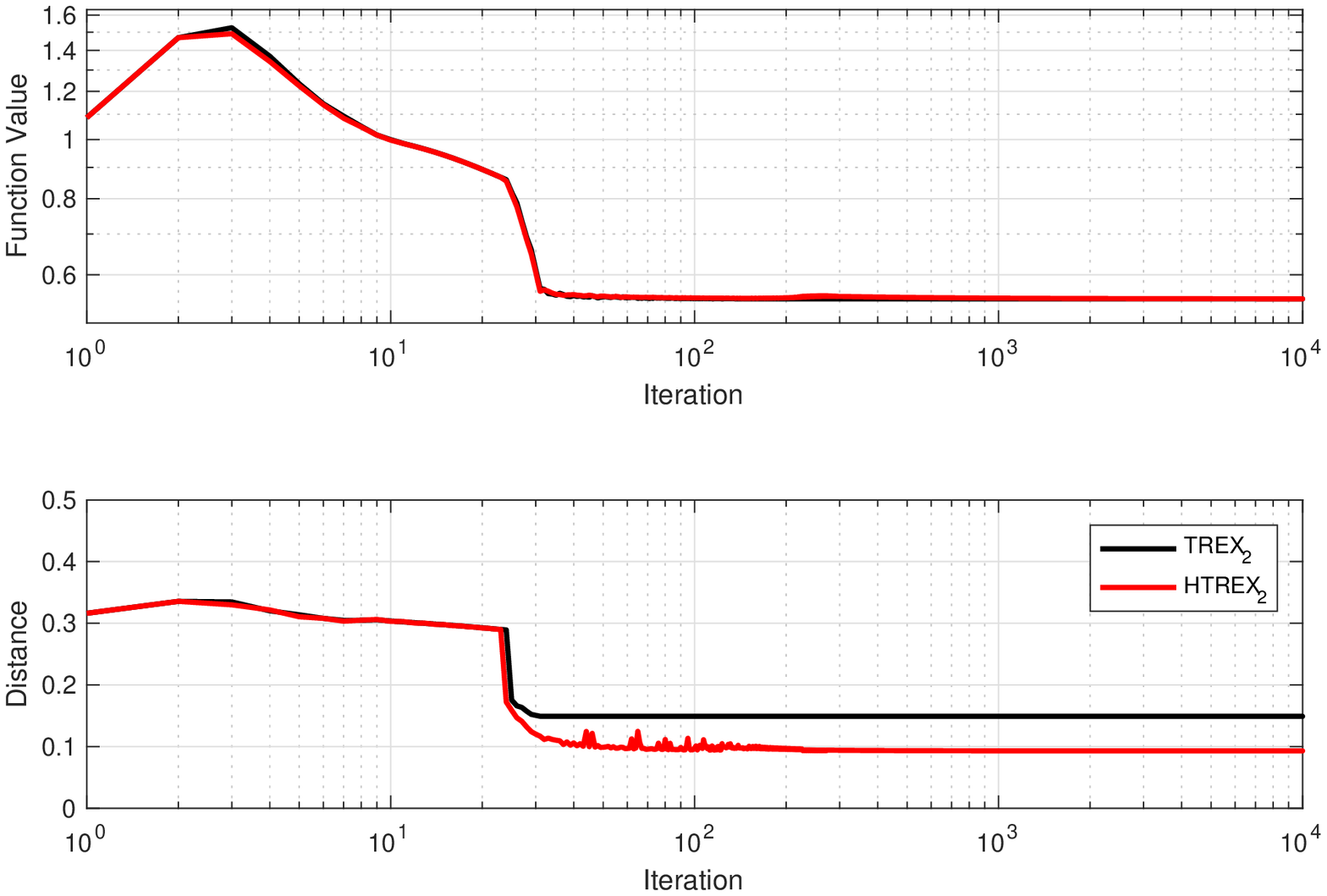}}
  \vspace{-8mm}
  (a)  Comparison of TREX$_2$ and HTREX$_2$ in the process of convergences under the noise $\mathbf{e}=0$.
  \centerline{\includegraphics[width=10cm]{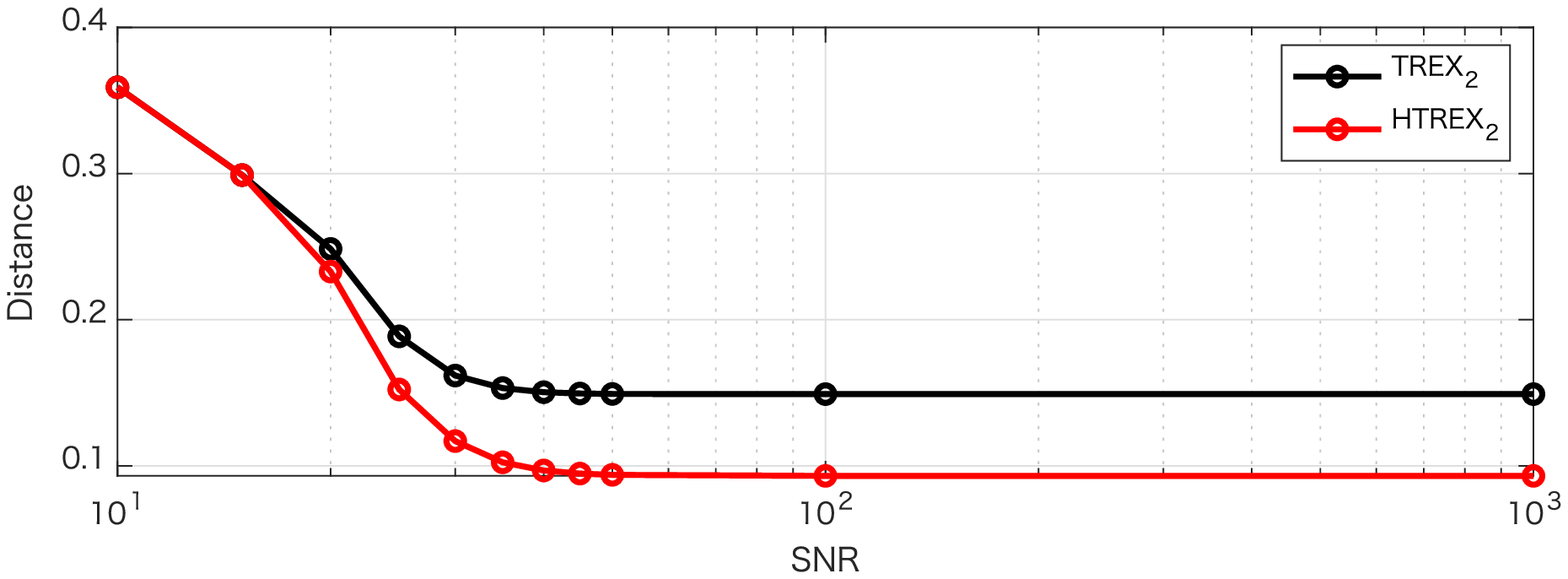}}
  (b) Estimation accuracy achieved by TREX$_2$ and HTREX$_2$ for various SNR.
  \caption{Transient performances in an under-determined case; Criteria ``Function Value'', ``Distance'' are given in \eqref{eq:criteriaTREX} and $({\rm SNR})=10 \log\left(\frac{\|{\mathbf X}{\mathbf b}^{\rm tru}\|^2}{\| \sigma {\mathbf e}\|^2}\right) \ [{\rm dB}]$.}
  \label{fig:noisefree5}
\end{figure}

 We compared the performances of ${\mathbf b}_{\rm TREX_2}$ in \eqref{decomp-TREXref} ($\beta =1/2$) and
 ${\mathbf b}_{\rm HTREX_2}$ in \eqref{eq:HNGTREX} employing $\Psi$ in \eqref{eq:fdf}.
 To approximate iteratively ${\mathbf b}_{\rm TREX_2}^{(j)}$ ($j=1,2,\ldots, 2p$) for \eqref{p-gtrex-subp} and ${\mathbf b}_{\rm HTREX_2}^{(j)}$ ($j=1,2,\ldots, 2p$) for \eqref{eq:HHE}, we used respectively TREX$_2$ \eqref{eq:DRSgTREX} (Fact \ref{fa:DRSTREX} with $\alpha_n=1.95$ ($n \in \mathbb{N}$)) and
the proposed algorithm \eqref{eq:HSDM311TREX} (HTREX$_2$ with $\alpha=1.95$ and $\lambda_{n}=\frac{1}{n}$ for $n \in \mathbb{N}$).
As performance measures, we used, in Figure \ref{fig:noisefree} and Figure \ref{fig:noisefree5},
\begin{align}
\left[
\begin{array}{ll}
\text{ \bf Function Value} \text{ (see \eqref{eq:HNGTREX})} & \min_{j=1,\ldots, 2p}(g_{(j,2)}({\mathbf M}_j{\mathbf b}_n)+\|{\mathbf b}_n\|_1),   \\
  \text{ \bf Distance} & \| {\mathbf b}_n - {\mathbf b}^{\rm tru}\|.
\end{array}
\right.
\label{eq:criteriaTREX}
\end{align}
The experiments were performed both in an over-determined case ($N=30$ and $p=20$) in Figure \ref{fig:noisefree} and an under-determined case ($N=20$ and $p=30$) in Figure \ref{fig:noisefree5}.

Figure \ref{fig:noisefree}(a) and \ref{fig:noisefree5}(a) illustrate the process of convergences of TREX$_2$ and HTREX$_2$ in the absence of noise, i.e., ${\mathbf e}={\mathbf 0} \in \mathbb{R}^N$. From these figures, we observe that (i) Function Values of TREX$_2$ and HTREX$_2$ converge to the same level, and that (ii) Distance (to ${\mathbf b}^{\rm tru}$) of HTREX$_2$ converges to a lower level than that of TREX$_2$. Figure \ref{fig:noisefree}(b) and \ref{fig:noisefree5}(b) summarize the behavior of Distance (to ${\mathbf b}^{\rm tru}$), against various SNR, by TREX$_2$ and HTREX$_2$ after $10000$ iterations.
For all the SNR, HTREX$_2$ seems to succeed in improving the performance of TREX$_2$.

\section{Concluding Remarks}
\label{sec:ConcludingRemarks}

In this paper, we have demonstrated how the modern proximal splitting operators can be plugged nicely into 
the hybrid steepest descent method (HSDM) for their applications to the hierarchical convex optimization problems 
which require further strategic selection of a most desirable vector from the set of all solutions of the standard convex optimization. For simplicity as well as for broad applicability, we have chosen to cast our target in the iterative approximation of a viscosity solution of the standard convex optimization problem, where the 1st stage cost function is given as a superposition of multiple nonsmooth convex functions, involving linear operators, while its viscosity solution is 
a minimizer of the 2nd stage cost function which is {\rm G\^{a}teaux} differentiable convex function with Lipschitzian gradient. The key ideas for the successful collaboration between 
the proximal splitting operators and the HSDM are not only in (i) the previously known expressions of the solution set of the standard convex optimization problem as the fixed point set of computable 
nonexpansive operators but also in (ii) linear relations build strategically between the solution set and the fixed point set. Fortunately, we have shown that such key ideas can be achieved by extending carefully 
the strategies behind the Douglas-Rachford splitting operators as well as the LAL operators defined in certain  product Hilbert spaces. We have also presented applications of the proposed algorithmic strategies to certain unexplored hierarchical enhancements of the support vector machine and the Lasso estimator.

\begin{acknowledgement}
Isao Yamada would like to thank Heinz H. Bauschke, Russell Luke, and Regina S. Burachik  
for their kind encouragement and invitation of the first author to the dream meeting: 
\emph{Splitting Algorithms, Modern Operator Theory, and Applications}  
(September 17-22, 2017) in  Oaxaca, Mexico where he had great opportunity to receive insightful deep comments by H\revb{\'{e}}dy At\revb{t}ouch. He would also like to thank Patrick Louis Combettes and Christian L. M$\ddot{\rm u}$ller for their invitation of the first author to a special mini-symposium \emph{Proximal Techniques for High-Dimensional Statistics} in the SIAM conference on Optimization 2017 (May 22-25, 2017) in Vancouver. Their kind invitations and their excellent approach to the TREX problem motivated very much the authors to study the application of the proposed strategies to the hierarchical enhancement of Lasso in this paper.
Isao Yamada would also like to thank Raymond Honfu Chan for his kind encouragement and invitation to \emph{the Workshop on Optimization in Image Processing} (June 27-30, 2016) at the Harvard University.
Lastly, the authors thank to Yunosuke Nakayama for his help in the numerical experiment related to the proposed hierarchical enhancement of the SVM.  
\end{acknowledgement}

 \section*{Appendices}
 \addcontentsline{toc}{section}{Appendix}
\renewcommand{\theequation}{A.\arabic{equation}}
\setcounter{equation}{0}

\subsection*{A: Proof of Proposition \ref{pr:DRS}(a)}
  Fact \ref{fa:pdsc}(i)$\Leftrightarrow$(ii) in Section \ref{subsec:2.1} yields
  \begin{align}
     \text{(\ref{eq:DRSq1}--\ref{eq:DRSq3})} \    \Leftrightarrow & \ (\exists \nu_{\star} \in {\cal K}) \ \nu_{\star} \in \partial f(x_{\star}) \text{ and } -\nu_{\star} \in \partial g(x_{\star}) \nonumber \\
    \ \Leftrightarrow & \ 0 \in \partial f(x_{\star})+\partial g(x_{\star}). \nonumber
\end{align}
The remaining follows from the proof in \cite[Proposition 18]{plc-jcp07}.
\hfill \hfill \qed
  
\subsection*{B: Proof of Proposition \ref{pr:LAL}(a)(d)}\label{asubsec:pr:proofs}
\newcommand{\Z}{{\cal X}}
\newcommand{\Hy}{{\cal K}}
  {\rm (a)} From \eqref{eq:LALq1} and \eqref{eq:LALq2}, there exists $(x_{\star},\nu_{\star}) \in  {\cal S}_{\rm pLAL} \times {\cal S}_{\rm dLAL}$.
Fact \ref{fa:pdsc}(i)$\Leftrightarrow$(ii) in Section \ref{subsec:2.1} yields the equivalence 
\begin{align}
  & \ 
      (x_{\star},\nu_{\star}) \in {\cal S}_{\rm pLAL} \times {\cal S}_{\rm dLAL} \text{ and } \eqref{eq:LALq3}
 \nonumber \\
\Leftrightarrow & \ 
   A^* \nu_{\star} \in \partial f(x_{\star}) \text{ and }
-\nu_{\star} \in \partial \iota_{\{0\}}(Ax_{\star})
 \olabel{eq:u9bqevbh2} \\ 
\Leftrightarrow & \
   A^* \nu_{\star} \in \partial f(x_{\star}) \text{ and }
Ax_{\star}=0
\label{eq:u9bqevbh25} \\ 
  \Leftrightarrow & \
           x_{\star}=  \operatorname{prox}_f(x_{\star} - A^*Ax_{\star} + A^*\nu_{\star})
                    \text{ and }
                    \nu_{\star}=  \nu_{\star} - A x_{\star}
\nonumber \\
\Leftrightarrow & \ (x_{\star}, \nu_{\star}) \in \operatorname{Fix}(T_{\rm LAL}). \label{eq:u9bqevbh3}
\end{align}
\hfill \hfill \qed \\
\noindent {\rm (d)}
Choose arbitrarily $(\bar{x},\bar{\nu}) \in \operatorname{Fix}(T_{\rm LAL})$, i.e., 
\begin{align*}
  (\bar{x},\bar{\nu})=T_{\rm LAL}(\bar{x},\bar{\nu})=
  \begin{pmatrix}
    \operatorname{prox}_{f}(\bar{x} - A^*A\bar{x}+A^*\bar{\nu}),
    \bar{\nu} - A\bar{x}
    \end{pmatrix}.
\end{align*}
Let $(x_n,\nu_n)_{n \in \mathbb{N}} \subset \mathcal{X} \times \mathcal{K}$ be generated, with any $(x_0,\nu_0) \in \mathcal{X} \times \mathcal{K}$, by
\begin{align}
  \label{eq:Aupdate}
   (x_{n+1}, \nu_{n+1})= T_{\rm LAL}(x_n, \nu_n)
  =  \begin{pmatrix}
    \operatorname{prox}_{f}(x_n - A^*Ax_n+A^*\nu_n),
    \nu_n - Ax_{n+1}
    \end{pmatrix}.
\end{align}
Then \cite[(B.3)]{YY2017} yields
\begin{align}
0 \leq  & 
  \|x_n-\bar{x} \|_{\Z}^2
  -\|x_{n+1} - \bar{x} \|_{\Z}^2
  -\| (x_{n+1}-x_n)-(\bar{x}-\bar{x}) \|_{\Z}^2
\nonumber \\ &\nonumber  \
  + \|Ax_{n+1} -  A\bar{x} - (Ax_n-A\bar{x})\|_{\Hy}^2 
               - \|Ax_n-A\bar{x}\|_{\Hy}^2
\nonumber \\ &\nonumber  \               
  -\|Ax_{n+1} -  A\bar{x}+\bar{\nu}-\nu_n\|_{\Hy}^2
               +\|\bar{\nu}- \nu_n \|_{\Hy}^2 \\
  = &
      \|x_n-\bar{x} \|_{\Z}^2
  -\|x_{n+1} - \bar{x} \|_{\Z}^2
  -\| x_{n+1}-x_n \|_{\Z}^2
\nonumber \\ &\nonumber  \
  + \|Ax_{n+1}  - Ax_n\|_{\Hy}^2 
               - \|Ax_n\|_{\Hy}^2
  -\|\bar{\nu}-\nu_{n+1}\|_{\Hy}^2
               +\|\bar{\nu}- \nu_n \|_{\Hy}^2  \nonumber \\ \nonumber
  \leq &
      (\|x_n-\bar{x} \|_{\Z}^2+\|\nu_n-\bar{\nu} \|_{\Hy}^2)
         -(\|x_{n+1} - \bar{x} \|_{\Z}^2
                        +\|\nu_{n+1}-\bar{\nu}\|_{\Hy}^2)
\nonumber \\ &  \
                        + (\|A\|_{\rm op}^2-1)\|x_{n+1}  - x_n\|_{\cal X}^2 - \|Ax_n\|_{\Hy}^2.
               \label{eq:ibbqeghoubwrj}
\end{align}
\eqref{eq:ibbqeghoubwrj} and $\|A\|_{\rm op}<1$ imply that $(\|x_n-\bar{x} \|_{\Z}^2 +\|\nu_n-\bar{\nu} \|_{\Hy}^2)_{n \in \mathbb{N}}$ decreases monotonically, i.e.,
 $(x_{n},\nu_n)_{n \in \mathbb{N}}$ is Fej\'{e}r monotone with respect to $\operatorname{Fix}(T_{\rm LAL})$, and $(\|x_n-\bar{x} \|_{\Z}^2 +\|\nu_n-\bar{\nu} \|_{\Hy}^2)_{n \in \mathbb{N}}$ converges to some $c\geq 0$. From this observation, we have
\begin{align*}
  & \ \sum_{n=0}^N\left[(1-\|A\|_{\rm op}^2)\|x_{n+1}  - x_n\|_{\cal X}^2 + \|Ax_n\|_{\Hy}^2\right] \\
  \leq & \ 
         \sum_{n=0}^N\left[(\|x_n-\bar{x} \|_{\Z}^2+\|\nu_n-\bar{\nu} \|_{\Hy}^2)
         -(\|x_{n+1} - \bar{x} \|_{\Z}^2
         +\|\nu_{n+1}-\bar{\nu}\|_{\Hy}^2)\right]  \\
=  & \ 
         (\|x_0-\bar{x} \|_{\Z}^2+\|\nu_0-\bar{\nu} \|_{\Hy}^2)
         -(\|x_{N+1} - \bar{x} \|_{\Z}^2
     +\|\nu_{N+1}-\bar{\nu}\|_{\Hy}^2) \\
  \to & \ (\|x_0-\bar{x} \|_{\Z}^2+\|\nu_0-\bar{\nu} \|_{\Hy}^2)
          - c < \infty \quad ( N \to \infty)
\end{align*}
and thus
\begin{align}
  \lim_{n\to \infty }\|x_{n+1}-x_n\|_{\cal X}= 0 \text{ and }  \lim_{n\to \infty } \|Ax_{n}\|_{\cal K}=  0.
  \label{eq:preconvergence}
\end{align}
By \cite[Theorem 9.12]{Frank2001}, the bounded sequence of $(x_n,\nu_n)_{n \in \mathbb{N}}$ has some subsequence $(x_{n_j},\nu_{n_j})_{j \in \mathbb{N}}$ which converges weakly to a some point, say $(x_{\star},\nu_{\star})$, in the Hilbert space $\mathcal{X} \times {\cal K}$.
Therefore,  by applying \cite[Theorem 9.1(iii)$\Leftrightarrow$(i)]{bauschke11} to $f \in \Gamma_0(\mathcal{X})$, we have
\begin{align}
  \label{eq:qefob}
  f(x_{\star}) \leq \underset{j \to \infty}{\operatorname{liminf}} \ f(x_{n_j})
\end{align}
and, by the Cauchy-Schwarz inequality and \eqref{eq:preconvergence},
\begin{align*}
  \|Ax_{\star}\|_{\cal K}^2 = & \  \langle Ax_{\star} - Ax_{n_j}, Ax_{\star} \rangle_{\cal K} +  \langle Ax_{n_j}, Ax_{\star} \rangle_{\cal K} \\
  \leq & \  \langle x_{\star} - x_{n_j}, A^*Ax_{\star} \rangle_{\cal X} +  \|Ax_{n_j}\|_{\cal K}\| Ax_{\star} \|_{\cal K} \to 0 \quad (j \to \infty),
\end{align*}
which implies $Ax_{\star}=0$.

Meanwhile, by \eqref{eq:Aupdate}, we have
\begin{align}\hspace{-5mm}
  & \ x_{n_j} = \operatorname{prox}_{f}(x_{n_j-1} - A^{*} Ax_{n_j-1} + A^{*} \nu_n)
= ({\rm I}+\partial f)^{-1}(x_{n_j-1} - A^{*} Ax_{n_j-1} + A^{*} \nu_{n_j-1})
 \nonumber \\ \hspace{-5mm}
  \Leftrightarrow & \ x_{n_j-1}-x_{n_j} - A^{*} Ax_{n_j-1} + A^{*} \nu_{n_j-1} \in \partial f(x_{n_j}) \nonumber \\ \hspace{-5mm}
  \Leftrightarrow & \ 
                    (\forall x \in \mathcal{X}) \ f(x_{n_j}) + \langle x_{n_j-1}-x_{n_j} - A^{*} Ax_{n_j-1} + A^{*} \nu_{n_j-1}, x - x_{n_j}\rangle_{\cal X}  \leq f(x), \label{eq:onaa}
\end{align}
where the inner product therein satisfies
\begin{align}
  \lim_{j \to \infty}\langle x_{n_j-1}-x_{n_j} - A^{*} Ax_{n_j-1} + A^{*} \nu_{n_j-1}, x - x_{n_j}\rangle_{\cal X}=\langle A^*\nu_{\star}, x - x_{\star} \rangle_{\cal X},     \label{eq:qeobqeg}
\end{align}
which is verified by $Ax_{\star}=0$, the triangle inequality, the Cauchy-Schwarz inequality, and \eqref{eq:preconvergence}, as follows:
\begin{align}
  & \ (\forall x \in \mathcal{X}) \nonumber \\
  & \ |\langle x_{n_j-1}-x_{n_j} - A^{*} Ax_{n_j-1} + A^{*} \nu_{n_j-1}, x - x_{n_j}\rangle_{\cal X} -\langle A^*\nu_{\star}, x - x_{\star} \rangle_{\cal X}| \nonumber \\
=  & \ |
\langle x_{n_j-1}-x_{n_j}, x - x_{n_j}\rangle_{\cal X}
-\langle Ax_{n_j-1}, A(x - x_{n_j})\rangle_{\cal K} \nonumber \\
  & \ \qquad \qquad \qquad \qquad  +\langle  \nu_{n_j-1}, A(x - x_{n_j})\rangle_{\cal K} -\langle \nu_{\star}, Ax \rangle_{\cal K}| \nonumber \\
=  & \ |
\langle x_{n_j-1}-x_{n_j}, x - x_{n_j}\rangle_{\cal X}
-\langle Ax_{n_j-1}, A(x - x_{n_j})\rangle_{\cal K} \nonumber \\
  & \ \qquad +\langle  \nu_{n_j-1},  - Ax_{n_j}\rangle_{\cal K}
    -\langle \nu_{n_j} - \nu_{n_j-1}, Ax  \rangle_{\cal K}
    -\langle \nu_{\star}-\nu_{n_j}, Ax  \rangle_{\cal K}| \nonumber \\
  \leq & \ (\| x_{n_j-1}-x_{n_j}\|_{\cal X}\| x - x_{n_j}\|_{\cal X}
         +\|Ax_{n_j-1}\|_{\cal K}\| A(x - x_{n_j})\|_{\cal K}  \nonumber \\         
  & \ \qquad
    +\| \nu_{n_j-1}\|_{\cal K}\| - Ax_{n_j}\|_{\cal K} + \|Ax_{n_j}\|_{\cal K}\| Ax\|_{\cal K} +|\langle \nu_{\star}-\nu_{n_j}, Ax  \rangle_{\cal K}|) \nonumber \\
  \to & 0 \quad ( j \to \infty).
        \nonumber
\end{align}
Now, by \eqref{eq:onaa}, \eqref{eq:qefob} and \eqref{eq:qeobqeg}, we have for any $ x \in \mathcal{X}$ 
\begin{align}   
  f(x) & \geq 
    f(x_{\star}) + \underset{j \to \infty}{\operatorname{liminf}} \langle x_{n_j-1}-x_{n_j} - A^{*} Ax_{n_j-1} + A^{*} \nu_{n_j-1}, x - x_{n_j}\rangle_{\cal X} \nonumber \\
       & = f(x_{\star}) +  \underset{j \to \infty}{\lim}\langle x_{n_j-1}-x_{n_j} - A^{*} Ax_{n_j-1} + A^{*} \nu_{n_j-1}, x - x_{n_j}\rangle_{\cal X} \nonumber \\
       & = f(x_{\star})+\langle A^{*} \nu_{\star}, x - x_{\star}\rangle_{\cal X},
         \nonumber
\end{align}
which implies
\begin{align}
   \label{eq:iubwrhr}
  A^{*} \nu_{\star} \in \partial f(x_{\star}).
\end{align}
By recalling \eqref{eq:u9bqevbh25}$\Leftrightarrow$\eqref{eq:u9bqevbh3},  
\eqref{eq:iubwrhr} and $Ax_{\star}=0$ prove $(x_{\star}, \nu_{\star}) \in \operatorname{Fix}(T_{\rm LAL})$.
The above discussion implies that every weak sequential cluster point \revb{(see Footnote \ref{fn:cluster} in Section~\ref{sec:Monotone})} of $(x_n,\nu_n)_{n \in \mathbb{N}}$, which is Fej\'{e}r monotone with respect to $\operatorname{Fix}(T_{\rm LAL})$, belongs to $\operatorname{Fix}(T_{\rm LAL})$. Therefore, \cite[Theorem 5.5]{bauschke11} guarantees that $(x_n,\nu_n)_{n \in \mathbb{N}}$ converges weakly to a point in $\operatorname{Fix}(T_{\rm LAL})$.
\hfill \hfill \qed

 \subsection*{C: Proof of Theorem \ref{th:DRS}}
 \label{asubsec:le:DRSa}
Now by recalling Proposition \ref{pr:DRS} in Section \ref{subsec:2.2} and Remark \ref{re:DRSbehind} in Section \ref{sec:3.1}, it is sufficient to prove Claim \ref{th:DRS}. Let $x_{\star} \in {\cal S}_p \not = \varnothing$. Then
the Fermat's rule, Fact \ref{fa:sc}(b) (applicable due to the qualification condition \eqref{eq:qualification0}) in Section \ref{subsec:2.1},
$\check{A}^*\colon \mathcal{K} \to \mathcal{X} \times \mathcal{K}\colon \nu \mapsto (A^*\nu, -\nu)$ for $\check{A}$ in \eqref{eq:checkA}, the property of $\iota_{\{0\}}$ in \eqref{eq:partialiota}, the straightforward calculations, and Fact \ref{fa:pdsc}(ii)$\Leftrightarrow$(i) (in Section \ref{subsec:2.1}) yield
\begin{align}
\hspace{0mm}
 x_{\star} \in {\cal S}_p  \Leftrightarrow & \ 
  0 \in \partial (f + g\circ A)(x_{\star})=\partial f(x_{\star}) + A^* \partial g(Ax_{\star})  \nonumber \\ \hspace{0mm} 
  \ \Leftrightarrow & \   \  y_{\star}=Ax_{\star} \text{ and }  0 \in \partial f(x_{\star}) + A^* \partial g({y}_{\star}) \nonumber \\ \hspace{0mm} 
  \ \Leftrightarrow & \ (\exists \nu_{\star} \in \mathcal{K})
 \ y_{\star}=Ax_{\star} \text{ and }  \left\{\begin{array}{l}
                                               A^*  \nu_{\star} \in \partial f(x_{\star}) \\
                                               -\nu_{\star} \in \partial g({y}_{\star})
                                             \end{array}
  \right. \nonumber \\ \hspace{0mm}     
  \ \Leftrightarrow & \ (\exists \nu_{\star} \in \mathcal{K})
 \ \check{A}(x_{\star},y_{\star})=0 \text{ and } \check{A}^*\nu_{\star} \in \partial F(x_{\star}, y_{\star}) \nonumber \\ \hspace{0mm}     
  \ \Leftrightarrow & \ (\exists \nu_{\star} \in \mathcal{K}) \
                      -\nu_{\star} \in \partial \iota_{\{0\}}(\check{A}(x_{\star},y_{\star}))
                      \text{ and } \check{A}^*\nu_{\star} \in \partial F(x_{\star}, y_{\star})
                      \nonumber \\ \hspace{0mm}     
  \ \Rightarrow & \ (\exists \nu_{\star} \in \mathcal{K}) \
                      -\check{A}^*\nu_{\star} \in \check{A}^*\partial \iota_{\{0\}}(\check{A}(x_{\star},y_{\star}))
                  \text{ and } \check{A}^*\nu_{\star} \in \partial F(x_{\star}, y_{\star})
                  \nonumber \\ \hspace{0mm}     
  \ \Rightarrow & \ (\exists \nu_{\star} \in \mathcal{K}) \
                      -\check{A}^*\nu_{\star} \in \partial (\iota_{\{0\}}\circ \check{A})(x_{\star},y_{\star})
                  \text{ and } \check{A}^*\nu_{\star} \in \partial F(x_{\star}, y_{\star})
  \nonumber \\ \hspace{0mm}     
  \ \Leftrightarrow & \ (\exists \nu_{\star} \in \mathcal{K}) \
                      -\check{A}^*\nu_{\star} \in \partial \iota_{\mathcal{N}(\check{A})}(x_{\star},y_{\star})
  \text{ and } \check{A}^*\nu_{\star} \in \partial F(x_{\star}, y_{\star})                  
\nonumber \\ \hspace{0mm}
  \  \Leftrightarrow & \
                       (\exists \nu_{\star} \in \mathcal{K}) \nonumber \\ \hspace{0mm}
  & \ 
                      \left\{\begin{array}{l}                       
                               (x_{\star}, y_{\star}) \in \operatorname{argmin}  (F + \iota_{\mathcal{N}(\check{A})})(\mathcal{X} \times \mathcal{K}) \\
                               \check{A}^*\nu_{\star} \in \operatorname{argmin}(F^*+\iota_{\mathcal{N}(\check{A})}^*\circ(-{\rm I}))(\mathcal{X} \times \mathcal{K}) \\
                               \min(F + \iota_{\mathcal{N}(\check{A})})(\mathcal{X} \times \mathcal{K})=- \min(F^*+\iota_{\mathcal{N}(\check{A})}^*\circ(-{\rm I}))(\mathcal{X} \times \mathcal{K}),
                                  \end{array} \right. \nonumber 
\end{align} 
which confirms Claim \ref{th:DRS}.
\hfill \hfill \qed

  \newcommand{\pdH}{\partial H}

  \subsection*{D: Proof of Theorem \ref{th:DRSave}}
  Now by recalling Proposition \ref{pr:DRS} in Section \ref{subsec:2.2} and Remark \ref{re:DRSavebehind} in Section \ref{sec:3.1}, it is sufficient to prove \eqref{eq:DRSIave2} by verifying Claim \ref{th:DRSave}. 
  We will use
  \begin{align}
    \label{eq:gasubdiff}
    \revb{A^* \circ \partial g \circ A=\sum_{i=1}^mA_i^* \circ \partial g_i \circ A_i = \sum_{i=1}^m\partial (g_i\circ A_i)}
  \end{align}
which is verified by \revb{$g=\bigoplus_{i=1}^mg_i$}, Fact \ref{fa:sc}(c) (see Section \ref{subsec:2.1}), and \rev{$\operatorname{ri}(\operatorname{dom}(g_j) - \operatorname{ran}(A_j))=\operatorname{ri}(\operatorname{dom}(g_j) - \mathbb{R})=\mathbb{R} \ni 0$ \revb{$(j=1,2,\ldots, m)$}.}
  Let $x_{\star}^{(m+1)} \in {\cal S}_p \not = \varnothing$.
  Then by using the Fermat's rule, Fact \ref{fa:sc}(b) (applicable due to \eqref{eq:qualification0}),
  \revb{\eqref{eq:gasubdiff}},
  $D$ in \eqref{eq:SDS}, and $H$ in \eqref{eq:H33}, we deduce the equivalence 
\begin{align}
  \hspace{-2mm}  & \  x_{\star}^{(m+1)} \in {\cal S}_p \nonumber \\ \hspace{-2mm}
  \ \Leftrightarrow & \ 0 \in \partial (f + g\circ A)(x_{\star}^{(m+1)}) = \partial f(x_{\star}^{(m+1)}) + A^* \partial g(Ax_{\star}^{(m+1)}) \nonumber \\ \hspace{-2mm}
   & \ \phantom{0 \in \partial (f + g\circ A)(x_{\star}^{(m+1)})} \revb{= \partial f(x_{\star}^{(m+1)}) + \sum_{i=1}^m \partial (g_i\circ A_i)(x_{\star}^{(m+1)})} \nonumber \\ \hspace{-2mm}  
         \Leftrightarrow & \ (j=1,\ldots,m) \ x_{\star}^{(j)}=x_{\star}^{(m+1)} \text{ and }  0 \in \partial f(x_{\star}^{(m+1)}) +\sum_{i=1}^m \partial (g_i\circ A_i)(x_{\star}^{(i)})
 \nonumber \\ \hspace{-2mm} 
  \ \Leftrightarrow & \  (\exists \nu^{(1)}, \ldots, \nu^{(m)} \in \mathcal{X}) (j =1,\ldots, m)
                      \left\{
                      \begin{array}{l}
                         x_{\star}^{(j)}=x_{\star}^{(m+1)} \\
                        \nu^{(j)} \in \partial (g_{j}\circ A_{j})(x_{\star}^{(j)}) \\
                        -\sum_{i=1}^m \nu^{(i)} \in \partial f(x_{\star}^{(m+1)}) 
                        \end{array}
  \right. \nonumber \\ \hspace{-2mm}
  \ \Leftrightarrow &
\  (\exists \nu^{(1)}, \ldots, \nu^{(m)} \in \mathcal{X}) \nonumber \\ \hspace{-2mm}
  & \hspace{-1mm} \! \left\{\!\!\!
                      \begin{array}{l}
                        (x_{\star}^{(1)},\ldots,x_{\star}^{(m+1)}) \in D \\
                        \begin{pmatrix}
                          \nu^{(1)},
                          \ldots,
                          \nu^{(m)},
                          -\sum_{i=1}^m \nu^{(i)}
                        \end{pmatrix} \! \in \!
                        \left[\varprod_{j=1}^{m}\partial (g_{j}\circ A_{j})(x_{\star}^{(j)})\right] \!\times \!\partial f(x_{\star}^{(m+1)}) \\
                        \phantom{
                        \begin{pmatrix}
                          \nu^{(1)},
                          \ldots,
                          \nu^{(m)},
                          -\sum_{i=1}^m \nu^{(i)}
                        \end{pmatrix} }
                        =\pdH(x_{\star}^{(1)},\ldots,x_{\star}^{(m+1)}).
                      \end{array}
  \right. 
  \olabel{eq:ouvcqefnoqeg}
\end{align}
Then by $-\begin{pmatrix}
                          \nu^{(1)},
                          \ldots,
                          \nu^{(m)},
                          -\sum_{i=1}^m \nu^{(i)}
                        \end{pmatrix} \in D^{\perp}=\partial \iota_{D}(x_{\star}^{(1)},\ldots,x_{\star}^{(m+1)})$
                        (see \eqref{eq:subdiffwgdobueqg3})
 and by Fact \ref{fa:pdsc}(ii)$\Leftrightarrow$(i) in Section \ref{subsec:2.1}, we have
                          \begin{align}
x_{\star}^{(m+1)} \in {\cal S}_p 
\ \Leftrightarrow &
\   (\exists \nu^{(1)}, \ldots, \nu^{(m)} \in \mathcal{X}) \nonumber \\
  & \ 
                      \left\{
                      \begin{array}{l}
                        -\begin{pmatrix}
                          \nu^{(1)},
                          \ldots,
                          \nu^{(m)},
                          -\sum_{i=1}^m \nu^{(i)}
                          \end{pmatrix}\in \partial \iota_{D}(x_{\star}^{(1)},\ldots,x_{\star}^{(m+1)}) \nonumber \\
                        \begin{pmatrix}
                          \nu^{(1)},
                          \ldots,
                          \nu^{(m)},
                          -\sum_{i=1}^m \nu^{(i)}
                          \end{pmatrix} \in \pdH(x_{\star}^{(1)},\ldots,x_{\star}^{(m+1)})   
                        \end{array}
                      \right.  \\
\ \Leftrightarrow &
                    \   (\exists \nu^{(1)}, \ldots, \nu^{(m)} \in \mathcal{X}) \nonumber \\
                            & \ 
                      \left\{\begin{array}{l}                       
                               (x_{\star}^{(1)},\ldots,x_{\star}^{(m+1)}) \in \operatorname{argmin}(H + \iota_{D})(\mathcal{X}^{m+1}) \\
                               \begin{pmatrix}
                          \nu^{(1)},
                          \ldots,
                          \nu^{(m)},
                          -\sum_{i=1}^m \nu^{(i)}
                          \end{pmatrix}  \in \operatorname{argmin}(H^*+\iota_{D}^*\circ(-{\rm I}))(\mathcal{X}^{m+1}) \\
                               \min(H + \iota_{D})(\mathcal{X}^{m+1})=- \min(H^*+\iota_{D}^*\circ(-{\rm I}))(\mathcal{X}^{m+1}),
                             \end{array} \right.
                            \nonumber
\end{align}
 which confirms Claim \ref{th:DRSave}.
\hfill \hfill \qed

\subsection*{E: Proof of Theorem \ref{th:HSDMLALI}}
Now by recalling Proposition \ref{pr:LAL} in Section \ref{subsec:2.2} and Remark \ref{re:LALIIbehind} in Section \ref{sec:3.2}, it is sufficient to prove Claim \ref{th:HSDMLALI}. 
Let $x_{\star} \in {\cal S}_p \not = \varnothing$. Then
the Fermat's rule, Fact \ref{fa:sc}(b) (applicable due to \eqref{eq:qualification0}) in Section \ref{subsec:2.1}, 
$\check{A}^*\colon \mathcal{K} \to \mathcal{X} \times \mathcal{K}\colon \nu \mapsto (A^*\nu, -\nu)$ for $\check{A}$ in \eqref{eq:checkA}, the property of $\iota_{\{0\}}$ in \eqref{eq:partialiota}, the straightforward calculations, and Fact \ref{fa:pdsc}(ii)$\Leftrightarrow$(i) (in Section \ref{subsec:2.1}) yield
\begin{align}
    x_{\star} \in {\cal S}_p
  \Leftrightarrow & \ 
  0 \in \partial (f + g\circ A)(x_{\star})=\partial f(x_{\star}) + A^* \partial g(Ax_{\star})  \nonumber \\ 
  \ \Leftrightarrow & \   \  y_{\star}=Ax_{\star} \text{ and }  0 \in \partial f(x_{\star}) + A^* \partial g({y}_{\star})  \nonumber \\ 
  \ \Leftrightarrow & \ (\exists \nu_{\star} \in \mathcal{K})
 \ y_{\star}=Ax_{\star} \text{ and }  \left\{\begin{array}{l}
                                               \mathfrak{u}A^*  \nu_{\star} \in \partial f(x_{\star}) \\
                                               -\mathfrak{u}\nu_{\star} \in \partial g({y}_{\star})
                                             \end{array}
  \right. \nonumber \\     
  \ \Leftrightarrow & \ (\exists \nu_{\star} \in \mathcal{K})
 \ (\sca\check{A})(x_{\star},y_{\star})=0 \text{ and } (\sca\check{A})^*\nu_{\star} \in \partial F(x_{\star}, y_{\star}) \nonumber \\     
  \ \Leftrightarrow & \ (\exists \nu_{\star} \in \mathcal{K}) \
                      -\nu_{\star} \in \partial \iota_{\{0\}}((\sca\check{A})(x_{\star},y_{\star}))
  \text{ and } (\sca\check{A})^*\nu_{\star} \in \partial F(x_{\star}, y_{\star}) 
\nonumber \\
  \  \Leftrightarrow & \
                       (\exists \nu_{\star} \in \mathcal{K}) \ \nonumber \\
                       & \quad 
                      \left\{\begin{array}{l}                       
                               (x_{\star}, y_{\star}) \in \operatorname{argmin} \ (F + \iota_{\{0\}}\circ (\sca\check{A}))(\mathcal{X} \times \mathcal{K}) \\
                               \nu_{\star} \in \operatorname{argmin}(F^*\circ (\sca\check{A})^*) (\mathcal{K}) \\
                               \min \ (F + \iota_{\{0\}}\circ (\sca\check{A}))(\mathcal{X} \times \mathcal{K})=- \min(F^*\circ (\sca\check{A})^*) (\mathcal{K}),
                                  \end{array} \right.
\nonumber 
\end{align}
which confirms Claim \ref{th:HSDMLALI}.
\hfill \hfill \qed

\subsection*{F: Proof of Theorem \ref{th:wrhwrh}}
\label{asubsec:th:wrhwrh}
(a) We have seen in \eqref{eq:DRScharacterization0} that, under the assumptions of Theorem \ref{th:wrhwrh}(a), for any vector $x_{\star} \in {\cal X}$,
\begin{align}
  \label{eq:DRSbounded1}
  \text{$x_{\star} \in {\cal S}_p\mbox{[in \eqref{viscosity-conv-prob}]}$ if and only if $(x_{\star}, y_{\star}) = {P}_{\mathcal{N}(\check{A})}({\mathbf \zeta}_{\star})$}
\end{align}
for some $y_{\star} \in {\cal X}$ and some ${\mathbf \zeta}_{\star}%
\in \operatorname{Fix}\left({\mathbf T}_{\rm DRS_{I}}\right)$, 
where $\check{A}\colon  \mathcal{X} \times \mathcal{K} \to \mathcal{K}\colon (x,y) \mapsto Ax-y$ (see \eqref{eq:checkA}),
$\mathcal{N}(\check{A})=\{(x,Ax) \in {\cal X} \times {\cal K}\mid x \in {\cal X} \}$,
and ${\mathbf T}_{\rm DRS_{I}}=(2\operatorname{prox}_F -{\rm I}) \circ (2{P}_{\mathcal{N}(\check{A})} -{\rm I})$ for $F\colon {\cal X} \times {\cal K} \to (-\infty,\infty]\colon (x,y)\mapsto f(x)+g(y)$
(see \eqref{eq:DRSIimplicit} and \eqref{eq:F31}).

Choose ${\mathbf \zeta}_{\star}:=(\zeta^x_{\star}, \zeta^y_{\star})  \in \operatorname{Fix}\left({\mathbf T}_{\rm DRS_{I}}\right)$ arbitrarily and let ${\mathbf z}_{\star}:=(x_{\star}, y_{\star}) := {P}_{\mathcal{N}(\check{A})}({\mathbf \zeta}_{\star})$. Then we have
\begin{align}
  & \ {\mathbf \zeta}_{\star}  \in \operatorname{Fix}\left({\mathbf T}_{\rm DRS_{I}}\right) \text{ and } {P}_{\mathcal{N}(\check{A})}({\mathbf \zeta}_{\star})={\mathbf z}_{\star} 
                      \nonumber  \\
  \Leftrightarrow & \ 
(2\operatorname{prox}_F-{\rm I})\circ (2{P}_{\mathcal{N}(\check{A})}-{\rm I})({\mathbf \zeta}_{\star})={\mathbf \zeta}_{\star}\text{ and } {P}_{\mathcal{N}(\check{A})}({\mathbf \zeta}_{\star})={\mathbf z}_{\star}
                    \label{eq:DRSbstart2} \\
  \Rightarrow & \ (2\operatorname{prox}_F-{\rm I})(2{\mathbf z}_{\star}-{\mathbf \zeta}_{\star})={\mathbf \zeta}_{\star}
                \ \Leftrightarrow  \ \operatorname{prox}_F(2{\mathbf z}_{\star}-{\mathbf \zeta}_{\star}) = {\mathbf z}_{\star}
    \nonumber \\
   \Leftrightarrow & \ ({\rm I}+\partial F)^{-1}(2{\mathbf z}_{\star}-{\mathbf \zeta}_{\star}) = {\mathbf z}_{\star} \
  \Leftrightarrow \ 2{\mathbf z}_{\star}-{\mathbf \zeta}_{\star} \in {\mathbf z}_{\star}+\partial F({\mathbf z}_{\star})
    \nonumber \\
  \Leftrightarrow & \ {\mathbf z}_{\star}-{\mathbf \zeta}_{\star} \in \partial F({\mathbf z}_{\star})=\partial f(x_{\star}) \times \partial g(y_{\star})
                    \label{eq:qeubqeDRS} \\
    \Leftrightarrow & \ 
                            x_{\star}-{\mathbf \zeta}_{\star}^x \in \partial f(x_{\star}) \text{ and }
                            y_{\star}-{\mathbf \zeta}_{\star}^y \in \partial g(y_{\star}).
\label{eq:jiuvqeqqhhrj}
\end{align}
Meanwhile, we have
\begin{align}
   {\mathbf z}_{\star}={P}_{\mathcal{N}(\check{A})}(\mathbf{\zeta}_{\star})  
 \ \Leftrightarrow & \ (\forall {\mathbf z}=(x,Ax) \in \mathcal{N}(\check{A})) \ \langle \zeta_{\star} - {\mathbf z}_{\star}, {\mathbf z} \rangle_{{\cal X} \times {\cal K}} = 0 \nonumber \\
  \Leftrightarrow &  \
                    (\forall x \in \mathcal{X}) \ \left\langle 
\zeta_{\star}^x 
- 
x_{\star}, x \right\rangle_{\mathcal{X}}
  +
  \left\langle 
\zeta_{\star}^y 
- 
y_{\star}
,   Ax \right\rangle_{\mathcal{K}}
  =0 \nonumber \\
  \Leftrightarrow & \
                    (\forall x \in \mathcal{X}) \ 
                    \left\langle 
                    (\zeta_{\star}^x - x_{\star}) +A^{*}(\zeta_{\star}^y -y_{\star}),
  x \right\rangle_{\mathcal{X}} =0  \nonumber \\
  \Leftrightarrow & \
                    A^{*}(\zeta_{\star}^y -y_{\star}) =-(\zeta_{\star}^x - x_{\star}).
 \label{eq:iyvqeroudgwbuoqey}
\end{align} 
Equations (\ref{eq:iyvqeroudgwbuoqey}) and \eqref{eq:jiuvqeqqhhrj} imply
\begin{align}
  \hspace{-3mm} & \  {\mathbf \zeta}_{\star}  \in \operatorname{Fix}\left({\mathbf T}_{\rm DRS_{I}}\right) \text{ and } {P}_{\mathcal{N}(\check{A})}({\mathbf \zeta}_{\star})={\mathbf z}_{\star} 
                      \nonumber  \\  \hspace{-3mm} 
\Rightarrow  
                & \
                  \revb{x_{\star}-\zeta_{\star}^x \in \partial f(x_{\star}) \text{ and } y_{\star}-\zeta_{\star}^y \in (-(A^{*})^{-1}(\partial f(x_{\star}))) \cap \partial g(y_{\star})}
\nonumber  \\
\hspace{-3mm}
\Rightarrow                                      
  & \ \zeta_{\star}\!=(\zeta_{\star}^x,\zeta_{\star}^y) \! \in\!  (x_{\star}, y_{\star}) \!-\! 
    \left(\partial f(x_{\star}) \times [(-(A^{*})^{-1}(\partial f(x_{\star}))) \cap \partial g(y_{\star})]\right).
\label{eq:char34}
\end{align}
Moreover, by noting that \eqref{eq:DRSbounded1} ensures $x_{\star} \in {\cal S}_p$ and $y_{\star}=Ax_{\star}$, we have from \eqref{eq:char34}
\begin{align}
  \hspace{-6mm}
  & \ {\mathbf \zeta}_{\star}  \in \operatorname{Fix}\left({\mathbf T}_{\rm DRS_{I}}\right) \text{ and } (x_{\star},Ax_{\star}) = {P}_{\mathcal{N}(\check{A})}({\mathbf \zeta}_{\star}) \nonumber \\ \hspace{-6mm}
   \Rightarrow                                      
  & \ \zeta_{\star}\in  (x_{\star},Ax_{\star})
 -  \left(\partial f(x_{\star}) \times [(-(A^{*})^{-1}(\partial f(x_{\star}))) \cap \partial g(Ax_{\star})]\right) \nonumber \\  \hspace{-6mm}
  \Rightarrow & \ \zeta_{\star} \in \bigcup_{x' \in {\cal S}_p }(x', Ax') - \bigcup_{x'' \in {\cal S}_p } \left(\partial f(x'') \times [(-(A^{*})^{-1}(\partial f(x''))) \cap \partial g(Ax'')]\right) \nonumber    \hspace{-6mm}
\end{align}
Since $\zeta_{\star}$ is chosen arbitrarily from $\operatorname{Fix}\left({\mathbf T}_{\rm DRS_{I}}\right)$, we have

\begin{align}
  & \ \operatorname{Fix}\left({\mathbf T}_{\rm DRS_{I}}\right) 
  \subset \bigcup_{x' \in {\cal S}_p }(x', Ax') - 
 \bigcup_{x'' \in {\cal S}_p } \left(\partial f(x'') \times [(-(A^{*})^{-1}(\partial f(x''))) \cap \partial g(Ax'')]\right),
    \olabel{eq:iuvbqeg2}
    \end{align}
    from which Theorem \ref{th:wrhwrh}(a) is confirmed.

\noindent {\bf \rm (b)} 
 We have seen in \eqref{eq:LALI2} that, under the assumptions of Theorem \ref{th:wrhwrh}(b), for any vector $x_{\star} \in {\cal X}$,
\begin{align}
  \label{eq:LALbounded1}
  \text{$x_{\star} \in {\cal S}_p\mbox{[in \eqref{viscosity-conv-prob}]}$ if and only if $(x_{\star}, y_{\star},\nu_{\star}) \in \operatorname{Fix}({\mathbf T}_{\rm LAL})$}
\end{align}
for some $(y_{\star},\nu_{\star}) \in {\cal K} \times {\cal K}$, 
where 
\begin{align}
  \hspace{-6mm} {\mathbf T}_{\rm LAL}&\colon {\cal X} \times {\cal K} \times {\cal K} \to  {\cal X} \times {\cal K} \times {\cal K} \nonumber \\ \hspace{-6mm} &\colon 
                                                                                                                                                       \begin{pmatrix}
x\\
y\\
\nu
\end{pmatrix}=
\begin{pmatrix}
{\mathbf z}\\
\nu
\end{pmatrix}
 \mapsto 
\begin{pmatrix}
  x_T\\
  y_T\\
  \nu_T
\end{pmatrix}
  =
\begin{pmatrix}
  {\mathbf z}_T \\
  \nu_T
\end{pmatrix}
 =
\left(
\begin{array}{l}
          \operatorname{prox}_F({\mathbf z}-(\sca\check{A})^*(\sca\check{A}){\mathbf z}+(\sca\check{A})^*\nu) \\
          \nu - \sca\check{A}{\mathbf z}_T
\end{array}
  \right)
  \nonumber
\end{align}
\revc{and $(\mathfrak{u} \check{A})^*\colon{\cal K} \to  {\cal X} \times {\cal K}\colon \nu \mapsto (\mathfrak{u}A^* \nu, -\mathfrak{u} \nu)$
(see %
\eqref{eq:LALexpress} and \eqref{eq:conjugateuAcheck}).}

Choose $({\mathbf z}_{\star}, \nu_{\star}) \in \operatorname{Fix}({\mathbf T}_{\rm LAL})$ arbitrarily and denote ${\mathbf z}_{\star}=(x_{\star},y_{\star}) \in \mathcal{X} \times \mathcal{K}$.
By passing similar steps in \eqref{eq:u9bqevbh3}$\Leftrightarrow$\eqref{eq:u9bqevbh25}, we deduce 
\begin{align}
   & \ ({\mathbf z}_{\star}, \nu_{\star}) \in \operatorname{Fix}({\mathbf T}_{\rm LAL})  \nonumber \\
 \ \Leftrightarrow & \ (\mathfrak{u}\check{A})^* \nu_{\star} \in \partial F({\mathbf z}_{\star})=\partial f(x_{\star}) \times \partial g(y_{\star}) \text{ and } \mathfrak{u}\check{A}({\mathbf z}_{\star}) = 0, 
\label{eq:oubqeg}
\end{align}
and then, from \eqref{eq:oubqeg}, straightforward calculations yield
\begin{align}
  & \ (x_{\star}, y_{\star}, \nu_{\star}) \in \operatorname{Fix}({\mathbf T}_{\rm LAL}) 
  \Leftrightarrow 
                     \left[\begin{array}{l}      
                                 \mathfrak{u}A^* \nu_{\star}=A^*(\mathfrak{u} \nu_{\star}) \in \partial f(x_{\star}) \\
                 -\mathfrak{u}\nu_{\star} \in \partial g( y_{\star}) \\
                                 Ax_{\star}=y_{\star}                                                                           \end{array} \right] \nonumber \\
  \Rightarrow & \
\revb{-\mathfrak{u} \nu_{\star} \in \left[-(A^*)^{-1}(\partial f(x_{\star}))\right] \cap \partial g( Ax_{\star}) \text{ and }  Ax_{\star}=y_{\star}}                
  \nonumber \\
\Leftrightarrow & \
                  -\mathfrak{u}(x_{\star}, y_{\star}, \nu_{\star}) \in   \{-\mathfrak{u}(x_{\star},Ax_{\star})\} \times [-(A^*)^{-1}(\partial f(x_{\star})) \cap \partial g( Ax_{\star})].
                  \label{eq:vebueh}
\end{align}
Moreover, from \eqref{eq:vebueh} and \eqref{eq:LALbounded1}, we have 
\begin{align}
  & \ (x_{\star},y_{\star}, \nu_{\star}) \in \operatorname{Fix}({\mathbf T}_{\rm LAL})  \nonumber \\
 \Rightarrow & \ -\mathfrak{u}(x_{\star}, y_{\star}, \nu_{\star}) \in  \bigcup_{x \in {\cal S}_p} \{-\mathfrak{u}(x,Ax)\} \times [-(A^*)^{-1}(\partial f(x)) \cap \partial g( Ax)]. \nonumber
\end{align}
Since $(x_{\star},y_{\star}, \nu_{\star})$ is chosen arbitrarily from $ \operatorname{Fix}({\mathbf T}_{\rm LAL})$, we have 
\begin{align*}
-\mathfrak{u} \operatorname{Fix}({\mathbf T}_{\rm LAL}) \subset \bigcup_{x \in {\cal S}_p} 
\ \{-\mathfrak{u}(x,Ax)\}  \times \left[-(A^*)^{-1}(\partial f(x)) \cap \partial g( Ax)\right]
\end{align*}
from which Theorem \ref{th:wrhwrh}(b) is confirmed.

\noindent (c)
 We have seen in \eqref{eq:DRSIave2b} that, under the assumptions of Theorem \ref{th:wrhwrh}(c), for any vector $x_{\star} \in {\cal X}$,
\begin{align}
  \label{eq:DRSavebounded1}
  \text{$x_{\star} \in {\cal S}_p\mbox{[in \eqref{viscosity-conv-prob}]}$ if and only if $(x_{\star}, x_{\star}, \ldots, x_{\star}) = {P}_{D}(\mathfrak{X}_{\star})$}
\end{align}
for some $
\mathfrak{X}_{\star}%
\in \operatorname{Fix}\left({\mathbf T}_{\rm DRS_{II}}\right)$, 
where 
$D=  \{(x^{(1)},\ldots,x^{(m+1)}) \in \mathcal{X}^{m+1} \mid x^{(i)}=x^{(j)} \ (i,j =1,2,\ldots, m+1)  \}$ (see \eqref{eq:SDS}),
$H\colon   \mathcal{X}^{m+1} \to (-\infty,\infty]\colon (x^{(1)},\ldots,x^{(m+1)}) \mapsto \sum_{i=1}^m g_i(A_ix^{(i)})+f(x^{(m+1)})$ (see \eqref{eq:H33}),
and ${\mathbf T}_{\rm DRS_{II}}=(2\operatorname{prox}_H -{\rm I}) \circ (2{P}_{D} -{\rm I})$ 
(see \eqref{eq:DRSIIexpress}) \revc{[For the availability of $\operatorname{prox}_H$ and ${P}_D$ as computational tools, see Remark \ref{re:DRSavebehind}(a)]}.

Choose $\mathfrak{X}_{\star}:=(\zeta_{\star}^{(1)},\ldots, \zeta_{\star}^{(m+1)})  \in \operatorname{Fix}\left({\mathbf T}_{\rm DRS_{II}}\right)$ arbitrarily, and let ${\mathbf X}_{\star}:=(x_{\star}, \ldots, x_{\star}) = {P}_D(\mathfrak{X}_{\star})$. Then we have
\begin{align}
  & \ \mathfrak{X}_{\star}  \in \operatorname{Fix}\left({\mathbf T}_{\rm DRS_{II}}\right) \text{ and } {P}_D(\mathfrak{X}_{\star})={\mathbf X}_{\star} \nonumber \\
  \Leftrightarrow & \   
                    (2\operatorname{prox}_H -{\rm I}) \circ (2{P}_D -{\rm I})(\mathfrak{X}_{\star})=\mathfrak{X}_{\star} \text{ and }
                       {P}_{D}(\mathfrak{X}_{\star}) ={\mathbf X}_{\star}.
                       \olabel{eq:F23}
\end{align}
Now, by passing similar steps for \eqref{eq:DRSbstart2}$\Rightarrow$\eqref{eq:qeubqeDRS},  we deduce that 
\begin{align}
  \hspace{-1mm} & \ \mathfrak{X}_{\star}  \in \operatorname{Fix}\left({\mathbf T}_{\rm DRS_{II}}\right) \text{ and } {P}_D(\mathfrak{X}_{\star})={\mathbf X}_{\star} \nonumber \\
  \hspace{-1mm}\Rightarrow &\ {\mathbf X}_{\star} - \mathfrak{X}_{\star} \in \partial H({\mathbf X}_{\star})= \left[\varprod_{j=1}^m \partial (g_j\circ A_j)(x_{\star})\right] \times \partial f(x_{\star})   \nonumber \\
  \hspace{-1mm}\Leftrightarrow & \ (j=1,2,\ldots, m) \ x_{\star} - \zeta_{\star}^{(j)}\in \partial (g_j\circ A_j)(x_{\star}) 
  \text{ and }
                    x_{\star} - \zeta_{\star}^{(m+1)}\in \partial f(x_{\star}) 
  \nonumber \\
  \hspace{-1mm}\Leftrightarrow & \ \revb{(j=1,2,\ldots, m) \ x_{\star} - \zeta_{\star}^{(i)}\in A_j^*\partial g_j(A_jx_{\star}) 
  \text{ and }
  x_{\star} - \zeta_{\star}^{(m+1)}\in \partial f(x_{\star}),}
\label{eq:jiuvqeqqhhrjave}
\end{align}
\revb{where the last equivalence follows from Fact \ref{fa:sc}(c) (applicable due to
  $\operatorname{ri}(\operatorname{dom}(g_j) - \operatorname{ran}(A_j))=\operatorname{ri}(\operatorname{dom}(g_j) - \mathbb{R})=\mathbb{R} \ni 0$). %
}
Meanwhile, we have 
\begin{align}
\hspace{-6mm}  {\mathbf X}_{\star} = {P}_D(\mathfrak{X}_{\star})\ \ \!\! \Leftrightarrow \!\! \  \  x_{\star}=\frac{1}{m+1}\sum_{i=1}^{m+1}\zeta_{\star}^{(i)} 
\ \ \! \Leftrightarrow \! \ \ x_{\star}- \zeta_{\star}^{(m+1)}=-\sum_{i=1}^{m}(x_{\star}-\zeta_{\star}^{(i)}).
  \label{eq:oubqege}
\end{align}
Equations \eqref{eq:oubqege} and \eqref{eq:jiuvqeqqhhrjave} imply
\begin{align} \hspace{-2mm}
  & \ \mathfrak{X}_{\star}  \in \operatorname{Fix}\left({\mathbf T}_{\rm DRS_{II}}\right) \text{ and } {P}_D(\mathfrak{X}_{\star})={\mathbf X}_{\star} \nonumber \\ \hspace{-2mm}
  \Rightarrow & \ 
\left\{
\begin{array}{l}
  (j=1,2,\ldots, m) \ x_{\star} - \zeta_{\star}^{(j)}\in A_j^* \partial g_j (A_jx_{\star}) \\
  x_{\star}- \zeta_{\star}^{(m+1)} \in  \partial f(x_{\star}) \cap \left[-\sum_{i=1}^{m} A_i^*\partial g_i( A_ix_{\star})\right]
\end{array}
\right. \nonumber \\ \hspace{-2mm}
\ \Rightarrow                                      
  & \ {\mathbf X}_{\star}-\mathfrak{X}_{\star} \! \in \!
    \left[\varprod_{j=1}^{m} A_j^* \partial g_j (A_jx_{\star})\right] \!\times \!
    \left[\partial f(x_{\star}) \cap \left(-\sum_{i=1}^{m}  A_i^*\partial g_i( A_ix_{\star})\right)\right].
    \label{eq:ubqegrh}
\end{align}
Moreover, by noting that \eqref{eq:DRSavebounded1} ensures $x_{\star} \in {\cal S}_p$, we have from \eqref{eq:ubqegrh}
\begin{align}
  & \ \mathfrak{X}_{\star}  \in \operatorname{Fix}\left({\mathbf T}_{\rm DRS_{II}}\right) \text{ and } {P}_D(\mathfrak{X}_{\star})={\mathbf X}_{\star}=(x_{\star},\ldots, x_{\star}) \nonumber \\ 
    \Rightarrow                                      
  & \ 
    \mathfrak{X}_{\star}
    \in
    {\cal S}_p^{m+1}- 
    \bigcup_{x \in {\cal S}_p }
    \left(
    \left[\varprod_{j=1}^{m} A_j^* \partial g_j (A_jx)\right] \times 
    \left[\partial f(x) \cap \left(-\sum_{i=1}^{m}  A_i^*\partial g_i( A_ix)\right)\right]
 \right).
    \nonumber
\end{align}
Since $\mathfrak{X}_{\star}$ is chosen arbitrarily from $\operatorname{Fix}({\mathbf T}_{\rm DRS_{II}})$, we have
\begin{align}
   \operatorname{Fix}\left({\mathbf T}_{\rm DRS_{II}}\right) 
  \subset 
  {\cal S}_p^{m+1} -
  \bigcup_{x \in {\cal S}_p }
  \left(
  \left[\varprod_{j=1}^{m} A_j^* \partial g_j (A_jx)\right] \times 
    \left[\partial f(x) \cap \left(-\sum_{i=1}^{m}  A_i^*\partial g_i( A_ix)\right)\right]
\right),
  \olabel{eq:DRSIboundedlast}
    \end{align}
from which Theorem \ref{th:wrhwrh}(c) is confirmed.
\hfill \hfill \qed

\subsection*{G: Proof of Lemma \ref{pr:TREXqualcond}}
Obviously, we have
from \eqref{def-gqj}
\begin{align}
  (j =1,2, \ldots, 2p) \quad \operatorname{dom}(g_{(j,q)})\supset
  \{ \eta \in \mathbb{R} \mid \eta > {\mathbf x}_j^{\top}{\mathbf z}\} \times \mathbb{R}^N.
             \label{eq:yivqehobuqrjoi2qh}
\end{align}
By recalling \revc{$0 \not = {\mathbf x}_{j} \in \mathbb{R}^{N}\text{ in \eqref{eq:def-xj}}$}
\revc{and ${\mathbf M}_{j} \in \mathbb{R}^{(N+1) \times p}$ in \eqref{def-Mj},}
we have  
  \begin{align*}
    (j =1,2, \ldots, 2p)
    \left[
    \begin{array}{l}
      \|{\mathbf X}^{\top}{\mathbf x}_j\|\geq \|{\mathbf x}_j\|^2>0 \\
      t
      (\|{\mathbf X}^{\top}{\mathbf x}_j\|)^{-2}
      {\mathbf M}_j{\mathbf X}^{\top}{\mathbf x}_j= \begin{pmatrix}
        t  \\
        t \frac{{\mathbf X}{\mathbf X}^{\top}{\mathbf x}_j}{\|{\mathbf X}^{\top}{\mathbf x}_j\|^2}
      \end{pmatrix}
              \quad 
      (\forall t \in \mathbb{R}),
    \end{array}
    \right.
  \end{align*}
and therefore 
\begin{align}
  \label{eq:iuvqetivqhpre24y}
  (j =1,2, \ldots, 2p) \quad {\mathbf M}_j\operatorname{dom}(\| \cdot \|_1) ={\mathbf M}_j(\mathbb{R}^p)  \supset \operatorname{span}\begin{pmatrix}
    1 \\
    \frac{{\mathbf X}{\mathbf X}^{\top}{\mathbf x}_j}{\|{\mathbf X}^{\top}{\mathbf x}_j\|^2}
\end{pmatrix}.
\end{align}
To prove $\operatorname{dom}(g_{(j,q)}) - {\mathbf M}_j\operatorname{dom}(\|\cdot\|_1)=\mathbb{R} \times \mathbb{R}^N$, choose arbitrarily $(\eta, {\mathbf y}) \in \mathbb{R} \times \mathbb{R}^N$.
Then \eqref{eq:yivqehobuqrjoi2qh} and \eqref{eq:iuvqetivqhpre24y} guarantee
\begin{align}
   \begin{pmatrix}
    \eta \\
    {\mathbf y}
\end{pmatrix}  
= & \
\begin{pmatrix}
  {\mathbf x}_j^{\top}{\mathbf z}+1 \\
{\mathbf y}+({\mathbf x}_j^{\top}{\mathbf z}+1-\eta)\frac{{\mathbf X}{\mathbf X}^{\top}{\mathbf x}_j}{\|{\mathbf X}^{\top}{\mathbf x}_j\|^2}
\end{pmatrix}
-
\begin{pmatrix}
  {\mathbf x}_j^{\top}{\mathbf z}+1-\eta \\
({\mathbf x}_j^{\top}{\mathbf z}+1-\eta)\frac{{\mathbf X}{\mathbf X}^{\top}{\mathbf x}_j}{\|{\mathbf X}^{\top}{\mathbf x}_j\|^2}
\end{pmatrix}
\nonumber \\
\in & \  \{ \tilde{\eta} \in \mathbb{R} \mid \tilde{\eta} > {\mathbf x}_j^{\top}{\mathbf z}\} \times \mathbb{R}^N
-\operatorname{span}\begin{pmatrix}
    1 \\
    \frac{{\mathbf X}{\mathbf X}^{\top}{\mathbf x}_j}{\|{\mathbf X}^{\top}{\mathbf x}_j\|^2}
\end{pmatrix} \nonumber \\
\subset & \ \operatorname{dom}(g_{(j,q)}) - {\mathbf M}_j\operatorname{dom}(\|\cdot\|_1), \nonumber
\end{align}
implying thus\rev{
\begin{align}
  \label{eq:qualTRXred}
  \operatorname{ri}(\operatorname{dom}(g_{(j,q)}) - {\mathbf M}_j\operatorname{dom}(\|\cdot\|_1))=\operatorname{ri}(\mathbb{R} \times \mathbb{R}^N)=\mathbb{R} \times \mathbb{R}^N \ni 0.
\end{align}
}
\hfill \hfill \qed

\subsection*{H: Proof of Theorem \ref{th:DRSTREX}}
\label{asubsec:th:DRSTREX}
By recalling Remark \ref{re:HSDMTREX} in Section \ref{subsec:5.2}, it is sufficient to prove Claim~\ref{th:DRSTREX}, for which we use the following inequality: for each $j=1,2,\ldots, 2p$, 
\begin{align}
  \label{eq:iuvqetivqhpre}
  (\forall (\eta,{\mathbf y}) \in \mathbb{R}\times \mathbb{R}^N) \quad  \left\|{\mathbf M}_j^{\top}  \begin{pmatrix}
    \eta \\
    {\mathbf y}
\end{pmatrix}\right\| \geq \left|\eta\|{\mathbf x}_j\|^2+ \langle {\mathbf x}_j, {\mathbf y} \rangle\right|,
\end{align}
where
${\mathbf x}_j \in \mathbb{R}^{N}$ in \eqref{eq:def-xj} and ${\mathbf M}_j \in \mathbb{R}^{(N+1) \times p}$ in \eqref{def-Mj}. Equation \eqref{eq:iuvqetivqhpre} is confirmed by
\begin{align}
  \olabel{eq:oubqehwrnj}
  (j=1,2,\ldots, 2p)(\forall (\eta,{\mathbf y}) \in \mathbb{R}\times \mathbb{R}^N) \quad 
  {\mathbf M}_j^{\top} 
  \begin{pmatrix}
    \eta \\
    {\mathbf y}
\end{pmatrix}
= \eta {\mathbf X}^{\top}{\mathbf x}_j  + {\mathbf X}^{\top}{\mathbf y}
\end{align}
and
\begin{align}\nonumber
  \begin{cases}
    \left[\eta {\mathbf X}^{\top}{\mathbf x}_j  + {\mathbf X}^{\top}{\mathbf y}\right]_j
    =\eta\|{\mathbf x}_j\|^2+ \langle {\mathbf x}_{j}, {\mathbf y} \rangle &  \text{if } j \in \{1,2, \ldots, p\} \\
    \left[\eta {\mathbf X}^{\top}{\mathbf x}_j  + {\mathbf X}^{\top}{\mathbf y}\right]_{j-p}=-\eta\|{\mathbf x}_j\|^2-\langle {\mathbf x}_{j}, {\mathbf y} \rangle &  \text{if }  j \in \{p+1,p+2, \ldots, 2p\}.
\end{cases}
\end{align}

Let
  $U_S:=\sup\{\| {\mathbf b}\| \mid {\mathbf b} \in S\}(<\infty)$.
  By supercoercivity of $\varphi$ and Example \ref{def-prop-perspective}, 
  the subdifferential %
  of its perspective $\widetilde{\varphi}$
  \revb{at each $(\eta, {\mathbf y}) \in \mathbb{R} \times \mathbb{R}^N$ can be expressed as \eqref{subdiff-perspective2}, and thus,}
to prove Claim~\ref{th:DRSTREX}, it is sufficient to show  
\begin{align*}
{\rm (i)}\quad & ({\mathbf M}_j^{\top})^{-1}(S) \cap \partial \widetilde{\varphi}(\mathbb{R}_{++} \times \mathbb{R}^{N}) \text{ is bounded;} \\
  {\rm (ii)}\quad & ({\mathbf M}_j^{\top})^{-1}(S) \cap \partial \widetilde{\varphi}(0,0) \text{ is bounded.}                    
\end{align*}

\noindent \underline{Proof of (i)} \  Choose $(\eta,{\mathbf y}) \in \mathbb{R}_{++} \times \mathbb{R}^{N}$ arbitrarily. Then, from \revb{\eqref{subdiff-perspective2}},
every ${\mathbf c}_{(\eta,{\mathbf y})} \in ({\mathbf M}_j^{\top})^{-1}(S) \cap  \partial \widetilde{\varphi}(\eta,{\mathbf y}) \in \mathbb{R} \times \mathbb{R}^{N}$ can be expressed with some ${\mathbf u} \in \partial \varphi({\mathbf y}/\eta)$ as 
\begin{align}
  \label{eq:basuobqeg}
{\mathbf c}_{(\eta,{\mathbf y})}=(\varphi({\mathbf y}/\eta) - \langle {\mathbf y}/\eta,{\mathbf u}\rangle, {\mathbf u})=(-\varphi^*({\mathbf u}), {\mathbf u}),
\end{align}
where the last equality follows from
$\varphi({\mathbf y}/\eta) +\varphi^*({\mathbf u})= \langle {\mathbf y}/\eta, {\mathbf u}\rangle$ due to
the Fenchel-Young identity \eqref{eq:FYidentity}.
By ${\mathbf M}_j^{\top}{\mathbf c}_{(\eta,{\mathbf y})} \in S$ and by applying the inequality \eqref{eq:iuvqetivqhpre} to \eqref{eq:basuobqeg}, we have
\begin{align}
  U_S \geq \|{\mathbf M}_j^{\top} {\mathbf c}_{(\eta,{\mathbf y})}\|
  =  \left\|{\mathbf M}_j^{\top}
  \begin{pmatrix}
    -\varphi^*({\mathbf u}) \\
    {\mathbf u}
  \end{pmatrix}
  \right\| 
  \geq & \   \left| (-\varphi^*({\mathbf u}))\|{\mathbf x}_j\|^2+ \left\langle {\mathbf x}_j , {\mathbf u}\right\rangle \right|
         \nonumber \\
  = & \ |\Upsilon({\mathbf u}) | 
  \geq \Upsilon_{+}({\mathbf u}), \label{eq:oiubqehwtjt}
\end{align}
where  
$\Upsilon\colon \mathbb{R}^N \to \mathbb{R}\colon {\mathbf v} \mapsto  \varphi^*({\mathbf v})\|{\mathbf x}_j\|^2- \left\langle {\mathbf x}_j , {\mathbf v}\right\rangle$
and $\Upsilon_{+}\colon \mathbb{R}^N \to \mathbb{R}\colon {\mathbf v} \mapsto \max\{\Upsilon({\mathbf v}), 0\}$ are coercive convex functions (see Section \ref{subsec:2.1}) and independent from the choice of $(\eta, {\mathbf y})$.
The coercivity of $\Upsilon_+$ ensures the existence of \revb{an open ball $B(0,\hat{U}_{\rm (i)})$ of radius $\hat{U}_{\rm (i)}>0$} such that $\operatorname{lev}_{\leq U_S} \Upsilon_+:=\{{\mathbf v} \in \mathbb{R}^N \mid \Upsilon_+({\mathbf v}) \leq U_S \} \subset B(0,\hat{U}_{\rm (i)} )$, and thus \eqref{eq:oiubqehwtjt} implies
\begin{align}
  \label{eq:qeoubvqrjjstj}
   \|{\mathbf u}\| \leq \hat{U}_{\rm (i)}.
\end{align}
Moreover, by ${\mathbf x}_j \not = 0$, the triangle inequality, the Cauchy-Schwarz inequality, \eqref{eq:oiubqehwtjt}, and \eqref{eq:qeoubvqrjjstj}, we have
 \begin{align}
   |\varphi^*({\mathbf u})|
   = & \ \left|\frac{\Upsilon({\mathbf u})}{\|{\mathbf x}_j\|^2} + \frac{\left\langle {\mathbf x}_j , {\mathbf u}\right\rangle}{\|{\mathbf x}_j\|^2}\right| 
       \leq  \left|\frac{\Upsilon({\mathbf u})}{\|{\mathbf x}_j\|^2}\right|+ \left|\frac{\left\langle {\mathbf x}_j , {\mathbf u}\right\rangle}{\|{\mathbf x}_j\|^2}\right| \nonumber \\
   \leq & \  \left|\frac{\Upsilon({\mathbf u})}{\|{\mathbf x}_j\|^2}\right|+ \frac{\| {\mathbf u}\|}{\|{\mathbf x}_j\|}
          \leq \frac{U_S}{\|{\mathbf x}_j\|^2} + \frac{\hat{U}_{\rm (i)} }{\|{\mathbf x}_j\|}=: U_{\rm (i)},
          \olabel{eq:qegjiubqegivye}
\end{align}
which yields ${\mathbf c}_{(\eta,{\mathbf y})}=(- \varphi^*({\mathbf u}), {\mathbf u}) \in [-U_{\rm (i)}, {U}_{\rm (i)}] \times B(0,\hat{U}_{\rm (i)} )$.
Since $(\eta,{\mathbf y})\in \mathbb{R}_{++} \times \mathbb{R}^{N}$ is chosen arbitrarily
and ${\mathbf c}_{(\eta,{\mathbf y})} \in ({\mathbf M}_j^{\top})^{-1}(S) \cap \partial \widetilde{\varphi}(\eta,{\mathbf y})$ is also chosen arbitrarily, we have
\begin{align*}
  ({\mathbf M}_j^{\top})^{-1}(S) \cap \partial \widetilde{\varphi}(\mathbb{R}_{++} \times \mathbb{R}^{N}) \subset [-U_{\rm (i)}, {U}_{\rm (i)}] \times B(0,\hat{U}_{\rm (i)} ),
\end{align*}
which confirms the statement (i).

\noindent \underline{Proof of (ii)} \
By introducing 
\begin{align}
  \label{eq:S1}
  \mathfrak{B}:= \left\{{\mathbf v} \in \mathbb{R}^N \left| \  \left|\left\langle  \frac{2}{\|{\mathbf x}_j\|^2}{\mathbf x}_j, {\mathbf v} \right\rangle\right| > |\varphi^*({\mathbf v})|\right. \right\},
\end{align}
we can decompose the set 
$({\mathbf M}_j^{\top})^{-1}(S) \cap \partial \widetilde{\varphi}(0,0)$ into 
\begin{align}
 \hspace{-5mm}
 ({\mathbf M}_j^{\top})^{-1}(S) \cap \partial \widetilde{\varphi}(0,0) \cap (\mathbb{R} \times \mathfrak{B}) \text{ and }
 ({\mathbf M}_j^{\top})^{-1}(S) \cap \partial \widetilde{\varphi}(0,0) \cap (\mathbb{R} \times \mathfrak{B}^c).
  \label{eq:twodecomp}
\end{align}
In the following, we show the boundedness of each set in \eqref{eq:twodecomp}.

First, we show the boundedness of $\mathfrak{B}$ by contradiction. Suppose that
$\mathfrak{B} \not \subset B(0,r)$
for all $r>0$. Then there exists a sequence $({\mathbf u}_k)_{k \in \mathbb{N}} \subset \mathbb{R}^N$ such that
\begin{align}
  \olabel{eq:cont2}
  (\forall k \in \mathbb{N}) \
  \frac{2}{\|{\mathbf x}_j\|} \geq  \left|\left\langle  \frac{2}{\|{\mathbf x}_j\|^2}{\mathbf x}_j, \frac{{\mathbf u}_k}{\|{\mathbf u}_k\|} \right\rangle\right| > \frac{|\varphi^*({\mathbf u}_k)|}{\|{\mathbf u}_k\|} \text{ and } \|{\mathbf  u}_k\| \geq k,
\end{align}
which contradicts the supercoercivity of $\varphi^*$, implying thus the existence of $r_* >0$ such that $\mathfrak{B} \subset B(0,r_*)$.

Next, we show the boundedness of the former set in \eqref{eq:twodecomp}.
Choose arbitrarily
\begin{align}
  \olabel{eq:wrjiqirbvuwrh2}
  (\mu,{\mathbf u}) \in ({\mathbf M}_j^{\top})^{-1}(S) \cap \partial \widetilde{\varphi}(0,0) \cap (\mathbb{R} \times \mathfrak{B}).
\end{align}
By $\mathbf{x}_j\not = 0$, ${\mathbf M}_j^{\top} (\mu,{\mathbf u}^{\top})^{\top} \in S \subset B(0,U_S)$, the inequality \eqref{eq:iuvqetivqhpre}, the triangle inequality, the Cauchy-Schwarz inequality, and ${\mathbf u} \in \mathfrak{B} \subset B(0,r_*)$, we have
  \begin{align}
  \frac{U_S}{\|{\mathbf x}_j\|^2} \geq & \ \frac{1}{\|{\mathbf x}_j\|^2}\left\|{\mathbf M}_j^{\top} \begin{pmatrix}
    \mu \\
    {\mathbf u}
    \end{pmatrix} \right\| 
  \geq \frac{1}{\|{\mathbf x}_j\|^2}|\mu\|{\mathbf x}_j\|^2+ \langle {\mathbf x}_j, {\mathbf u} \rangle| \nonumber \\
    \geq  & \ |\mu|- \left|\left\langle  \frac{{\mathbf x}_j}{\|{\mathbf x}_j\|^2}, {\mathbf u} \right\rangle\right|
            \geq |\mu|-  \frac{\| {\mathbf u}\|  }{\|{\mathbf x}_j\|} 
            \geq  |\mu|- \frac{r_*}{\|{\mathbf x}_j\|}
            \nonumber
\end{align}
which yields 
\begin{align}
  \nonumber
  \hat{U}_{\rm (iia)}  :=\frac{U_S}{\|{\mathbf x}_j\|^2}
+\frac{r_*}{\|{\mathbf x}_j\|} \geq |\mu|. 
\end{align}
Therefore, we have $(\mu,{\mathbf u}) \in [-\hat{U}_{\rm (iia)}, \hat{U}_{\rm (iia)}] \times B(0,r_{\star})$. Since $(\mu,{\mathbf u}) \in ({\mathbf M}_j^{\top})^{-1}(S) \cap \partial \widetilde{\varphi}(0,0) \cap (\mathbb{R} \times \mathfrak{B})$ is chosen arbitrarily, we have
\begin{align}
  ({\mathbf M}_j^{\top})^{-1}(S) \cap \partial \widetilde{\varphi}(0,0)\cap (\mathbb{R} \times \mathfrak{B}) \subset [-\hat{U}_{\rm (iia)}, \hat{U}_{\rm (iia)}] \times B(0,r_{\star}).
  \label{eq:iia}
\end{align}

Finally, we show the boundedness of the latter set in \eqref{eq:twodecomp}. Let
\begin{align}
  \label{eq:wrjiqirbvuwrh1}
  (\mu,{\mathbf u}) \in ({\mathbf M}_j^{\top})^{-1}(S) \cap \partial \widetilde{\varphi}(0,0) \cap (\mathbb{R} \times \mathfrak{B}^c).
\end{align}
\revb{From \eqref{subdiff-perspective2}}, we have
\begin{align}
  \label{eq:subdiffII}
  \partial \widetilde{\varphi}(0,0)= & \ \{(\mu',{\mathbf u}') \in \mathbb{R} \times \mathbb{R}^N \mid \mu' + \varphi^*({\mathbf u}') \leq 0 \}. 
\end{align}
Note that coercivity of $\varphi^*$ ($\Rightarrow \exists \min \varphi^*(\mathbb{R}^N) \in \mathbb{R}$, see Fact \ref{fa:uobqegbuqeh}) and \eqref{eq:subdiffII} yield
$\revb{\varphi^*({\mathbf u}) \in [\min \varphi^*(\mathbb{R}^N), -\mu]}$
and thus 
\begin{align}
  \label{eq:uviuyiqgborjh}
 |\varphi^*({\mathbf u})| \revb{\leq \max\{|\min \varphi^*(\mathbb{R}^N)|, |\mu| \}} \leq |\min \varphi^*(\mathbb{R}^N)|+ |\mu|.
  \end{align}
By $\mathbf{x}_j\not = 0$, ${\mathbf M}_j^{\top} (\mu,{\mathbf u}^{\top})^{\top} \in S \subset B(0,U_S)$ (see \eqref{eq:wrjiqirbvuwrh1}), the inequality \eqref{eq:iuvqetivqhpre}, the triangle inequality, ${\mathbf u} \in \mathfrak{B}^c$ (see \eqref{eq:wrjiqirbvuwrh1} and \eqref{eq:S1}), and \eqref{eq:uviuyiqgborjh}, we have
  \begin{align}
  \frac{2}{\|{\mathbf x}_j\|^2}U_S \geq & \ \frac{2}{\|{\mathbf x}_j\|^2}\left\|{\mathbf M}_j^{\top} \begin{pmatrix}
      \mu \\ 
      {\mathbf u}
    \end{pmatrix}\right\|
  \geq  \frac{2}{\|{\mathbf x}_j\|^2}|\mu\|{\mathbf x}_j\|^2+ \langle {\mathbf x}_j, {\mathbf u} \rangle| \nonumber \\
                                        \geq & \ 2|\mu|- \left|\left\langle  \frac{2}{\|{\mathbf x}_j\|^2}{\mathbf x}_j, {\mathbf u} \right\rangle\right| 
  \geq  2|\mu|- |\varphi^*({\mathbf u})| \nonumber \\
  \geq & \ 2|\mu|- |\min \varphi^*(\mathbb{R}^N)|- |\mu| 
   =  |\mu|- |\min \varphi^*(\mathbb{R}^N)| \nonumber
\end{align}
and thus, \revb{with \eqref{eq:uviuyiqgborjh},
\begin{align}
  \hat{U}_{\rm (iib)}  :=\frac{2}{\|{\mathbf x}_j\|^2}U_S  
+2|\min \varphi^*(\mathbb{R}^N)| \geq |\mu|+|\min \varphi^*(\mathbb{R}^N)| \geq |\varphi^*({\mathbf u})| \geq \varphi^*({\mathbf u}).
  \olabel{eq:qegobwrh}
\end{align}}
Hence, we have
\begin{align*}
 (\mu,{\mathbf u}) \in [-\hat{U}_{\rm (iib)}, \hat{U}_{\rm (iib)}] \times \operatorname{lev}_{\leq \hat{U}_{\rm (iib)}}(\varphi^*).
\end{align*}
Since $(\mu,{\mathbf u}) \in ({\mathbf M}_j^{\top})^{-1}(S) \cap \partial \widetilde{\varphi}(0,0)\cap (\mathbb{R} \times \mathfrak{B}^c)$ is chosen arbitrarily, we have 
\begin{align}
  \label{eq:iib}
  ({\mathbf M}_j^{\top})^{-1}(S) \cap \partial \widetilde{\varphi}(0,0)\cap (\mathbb{R} \times \mathfrak{B}^c) \subset [-\hat{U}_{\rm (iib)}, \hat{U}_{\rm (iib)}] \times \operatorname{lev}_{\leq \hat{U}_{\rm (iib)}}(\varphi^*).
\end{align}
Consequently, by using \eqref{eq:iia} and \eqref{eq:iib} and by letting $U_{\rm (ii)}:=\max\{\hat{U}_{\rm (iia)},\hat{U}_{\rm (iib)} \}$, we have 
\begin{align*}
  ({\mathbf M}_j^{\top})^{-1}(S) \cap \partial \widetilde{\varphi}(0,0) \subset [-U_{\rm (ii)}, U_{\rm (ii)}] \times [\operatorname{lev}_{\leq U_{\rm (ii)}}(\varphi^*) \cup B(0,r_{\star})],
\end{align*}
which guarantees the boundedness of $({\mathbf M}_j^{\top})^{-1}(S) \cap \partial \widetilde{\varphi}(0,0)$, due to the coercivity of $\varphi^*$, implying thus finally the statement (ii).

\hfill \hfill \qed

{\small
\bibliographystyle{YYspmpsci}
\bibliography{YY-2019-springer}
}

\end{document}